\documentclass[10pt, reqno]{amsart}

\usepackage{hyperref} 
\usepackage{amssymb, mathrsfs,bbold}
\usepackage{amsmath, amsthm}
\usepackage{enumerate}
\usepackage{dsfont}
\usepackage{color}
\usepackage[a4paper]{geometry}
\usepackage[utf8]{inputenc}   
\usepackage{stmaryrd}
\usepackage{titlesec, wasysym}
\usepackage{appendix}
\usepackage{todonotes, fancyhdr}
\usepackage{chngcntr}
\usepackage{cancel}
\usepackage{relsize,exscale}
\usepackage{bm}

\usepackage{tikz-cd}
\usetikzlibrary{calc}

\usepackage{setspace}
\renewcommand{\baselinestretch}{0.99}

\usepackage[all]{xy}
\numberwithin{subsection}{section}
\numberwithin{subsubsection}{subsection}
\numberwithin{equation}{section} 

\newenvironment{Dem}[1][\unskip]{%
    \begin{list}{\hspace{1.15cm}\textit{\textbf{Proof #1 --}}}{   
        \setlength{\topsep}{0pt}%
        \setlength{\leftmargin}{0pt}%
        \setlength{\rightmargin}{0pt}%
        \setlength{\listparindent}{0pt}%
        \setlength{\itemindent}{0pt}%
        \setlength{\parsep}{0pt}%
        \addtolength{\leftmargin}{0pt} 
        \addtolength{\rightmargin}{0pt}%
    } \item }{\hfill $\rhd$\end{list}\smallskip}

\newenvironment{Dem*}[1][\unskip]{%
    \begin{list}{\hspace{0cm}{\sf \textbf{{\small Proof #1 --}}}}{   %
        \setlength{\topsep}{0pt}%
        \setlength{\leftmargin}{0pt}%
        \setlength{\rightmargin}{0pt}%
        \setlength{\listparindent}{0pt}%
        \setlength{\itemindent}{0pt}%
        \setlength{\parsep}{0pt}%
        \addtolength{\leftmargin}{20pt}%
        \addtolength{\rightmargin}{0pt}%
    } \item }{\hfill $\rhd$\end{list}\smallskip}

\renewcommand\thesection       {\arabic{section}}
\renewcommand\thesubsection    {\thesection{\boldmath $.$}\arabic{subsection}}
\renewcommand\thesubsubsection    {\thesection{\boldmath $.$}\arabic{subsection}{\boldmath $.$}\arabic{subsubsection}} 

\titleformat{\section}[block] 
{\filcenter\normalfont\sffamily\bfseries\Large}  
{{\hspace{0cm}}\thesection \hspace{0.2em} --\vspace{0cm}}{0.5em}{} 

\titleformat{\subsection}[runin]
{\filcenter\normalfont\sffamily\bfseries\large}  
{{\hspace{0cm}}\thesubsection \hspace{0.15em} -- \vspace{0.1cm}}{.2em}{}   
\titlespacing{\subsection}{-0pc}{1.5ex plus .1ex minus .2ex}{0pc}   

\titleformat{\subsubsection}[runin]
{\filcenter\normalfont\sffamily\bfseries}   
{\filright\sffamily{\hspace{0cm}}\thesubsubsection\hspace{0.2em} --}{.5em}{}\titlespacing{\subsection}{-0pc}{1.5ex plus .1ex minus .2ex}{0pc}

\newtheoremstyle{mystyle}
{3pt}               
{3pt}               
{\it }                      
{}                      
{\bfseries}      
{}                      
{0.5em}                 
{\hspace{0cm}\textit{#2 --} {\hspace{-0.02cm}}\textit{#1}}   

\theoremstyle{mystyle}

\newtheorem{thm}{Theorem.}   
\newtheorem*{thm*}{Theorem.}

\newtheorem{cor}[thm]{{Corollary.} }
\newtheorem{lem}[thm]{{Lemma}. }
\newtheorem{prop}[thm]{{Proposition.}}
\newtheorem{defn}[thm]{{Definition.}}
\newtheorem*{rem*}{Remark.}

\newtheoremstyle{mystyle3}
{3pt}               
{3pt}               
{\it }                      
{}                      
{\bfseries}      
{}                      
{0.5em}                 
{\hspace{-0.8cm}{\textbf{\textit{#2}} --} {\hspace{-0.02cm}}{\textbf{\textit{#1}}}}

\theoremstyle{mystyle3}

\newtheoremstyle{mystyle2}
{3pt}               
{3pt}               
{\it }                      
{}                      
{\bfseries}    
{}                      
{0.5em}                 
{\llap{#2 }{\textbf{\textit{#1}}{\hspace{0.2cm}--}}}

\theoremstyle{mystyle2}

\newtheorem*{definition*}{Definition}
\newtheorem*{theorem*}{Theorem}
\newtheorem*{Remark*}{Remark}
\newtheorem*{lem*} {Lemma}
\newtheorem*{defn*} {Definition}
\newtheorem*{prop*} {Proposition}
\newtheorem*{cor*} {Corollary}

\newcommand{\ssk}{\smallskip}

\renewcommand{\epsilon}{\varepsilon}

\newcommand\bbN{\textbf{\textsf{N}}}
\newcommand\bbR{\textbf{\textsf{R}}}

\newcommand\bbZ{\mathbb{Z}}

\newcommand{\mcA}{\mathcal{A}}
\newcommand{\mcB}{\mathcal{B}} 
 
\newcommand{\mcD}{\mathcal{D}}
\newcommand{\mcE}{\mathcal{E}}

\newcommand{\mcL}{\mathcal{L}}

\newcommand{\mcP}{\mathcal{P}}

\newcommand\mcT{\mathcal T}
\newcommand\mcU{\mathcal U}
\newcommand\mcW{\mathcal W}

\newcommand{\bbfa}{{\color{blue} \bm{a}}}

\newcommand{\bbfb}{{\color{blue} \bm{b}}}

\newcommand{\bbfk}{{\color{blue} \bm{k}}}

\newcommand{\bbfl}{{\color{blue} \bm{\ell}}}

\newcommand{\bbfm}{{\color{blue} \bm{m}}}
\newcommand{\bbfmp}{{\color{blue} \bm{m}'}}
\newcommand{\bfN}{\mathbf{N}}

\newcommand{\bbfp}{{\color{blue} \bm{p}}}

\newcommand{\bbfq}{{\color{blue} \bm{q}}}
\newcommand{\bfR}{\mathbf{R}}
\newcommand{\bbfr}{{\color{blue} \bm{r}}}

\newcommand{\bbfv}{{\color{blue} \bm{v}}}

\makeatletter
\newcommand*{\defeq}{\mathrel{\rlap{%
                     \raisebox{0.3ex}{$\m@th\cdot$}}%
                     \raisebox{-0.3ex}{$\m@th\cdot$}}%
                     =}
\makeatother

\makeatletter
\newcommand*{\eqdef}{=\mathrel{\rlap{%
                     \raisebox{0.3ex}{$\m@th\cdot$}}%
                     \raisebox{-0.3ex}{$\m@th\cdot$}}%
                     }
\makeatother

\newcommand{\norme}[1]{\left\lVert #1 \right\lVert}

\newcommand{\dd}{\text{d}}
\renewcommand{\phi}{\varphi}

\begin{document}

\begin{center}
{\LARGE\sffamily{\textbf{Local expansion properties of paracontrolled systems}   \vspace{0.5cm}}}
\end{center}

\begin{center}
{\sf I. BAILLEUL and N. MOENCH}
\end{center}

\vspace{1cm}

\begin{center}
\begin{minipage}{0.8\textwidth}
\renewcommand\baselinestretch{0.7} \scriptsize \textbf{\textsf{\noindent Abstract.}} The concept of concrete regularity structure gives the algebraic backbone of the operations involved in the local expansions used in the regularity structure approach to singular stochastic partial differential equations. The spaces and the details of the structures depend on each equation. We introduce here a parameter-dependent universal algebraic regularity structure that can host all the regularity structures used in the study of singular stochastic partial differential equations. This is done by using the correspondence between the notions of model on a regularity structure and the notion of paracontrolled system. We prove that the iterated paraproducts that form the fundamental bricks of paracontrolled systems have some local expansion properties that are governed by this universal structure.
\end{minipage}
\end{center}

\vspace{0.4cm}

\section{Introduction}
\label{SectionIntro}

We work in the Euclidean space $\bbR^{d_0}$. The Besov-H\"older spaces $C^{\alpha_1}$ over $\bbR^{d_0}$ and their norms $\Vert\cdot\Vert_{\alpha_1}$ are defined as usual for any $\alpha_1\in\bbR$ from the Littlewood-Paley projectors $\Delta_i : \mcD'(\bbR^{d_0})\rightarrow C^\infty(\bbR^{d_0})$ and, setting 
$$
\Delta_{<j} \defeq \sum_{i\leq j-1}\Delta_i,
$$
the paraproduct ${\sf P}(f,g)$ of any two distributions $f,g$ is defined as
$$
{\sf P}(f,g) \defeq  \sum_{i\geq 1} \Delta_{<i-1}(f)\,\Delta_i(g).
$$
For $f\in C^{\alpha_1}$ and $g\in C^{\alpha_2}$ we have the optimal continuity estimates
$$
\Vert {\sf P}(f,g)\Vert_{\alpha_2} \lesssim \Vert f\Vert_\infty \Vert g\Vert_{\alpha_2}   \qquad   \textrm{if $\alpha_1\geq 0$}
$$
and 
$$
\Vert {\sf P}(f,g)\Vert_{\alpha_1+\alpha_2} \lesssim \Vert f\Vert_{\alpha_1} \Vert g\Vert_{\alpha_2}   \qquad   \textrm{if $\alpha_1<0$}.
$$

\medskip

\subsection{Local expansion properties of iterated paraproducts$\boldmath{.}$ \hspace{0.1cm}}
\label{IntroExpansionIPP}

We define inductively the iterated paraproduct operator by setting ${\sf P}(f)=f$, for any distribution $f\in\mcD'(\bbR^{d_0})$, and 
$$
{\sf P}(f_1,\dots,f_n) = {\sf P}\big({\sf P}(f_1,\dots,f_{n-1}),f_n\big)
$$
for any $n\geq 2$ and any distributions $f_1,\dots, f_n$ in $\mcD'(\bbR^{d_0})$. If $f_n\in C^{\alpha_n}$, the above continuity estimates on the paraproduct imply that the iterated paraproduct ${\sf P}(f_1,\dots,f_n)$ is in some $C^\gamma$ space where $\gamma\leq \alpha_n$. Such distributions may nonetheless have some local descriptions to an accuracy strictly larger that $\alpha_n$ around an arbitrary point. This was noticed for instance in Corollary 1 of Bailleul \& Bernicot's work \cite{BailleulBernicotHighOrder} in the case of a paraproduct ${\sf P}(f,g)$ where $f\in C^{\alpha_1}, g\in C^{\alpha_1}$ and $0<\alpha_1<1/2$. One has indeed in that case
\begin{equation} \label{EqFirstLocalExpansionParaproduct}
\big\vert {\sf P}(f,g)(y) - {\sf P}(f,g)(x) - f(x) \big(g(y)-g(x)\big)  \big\vert \lesssim \vert y-x\vert^{2\alpha_1},
\end{equation}
so one can give in that setting a local description of the behaviour of ${\sf P}(f,g)$ around an arbitrary point $x$ up to a precision $\vert y-x\vert^{2\alpha_1}$. For $g\in C^{\alpha_2}$ with $\alpha_2>0$ we write
$$
R^{\alpha_2}(g)(y,x) \defeq g(y) - \sum_{\vert k\vert<\alpha_2} \partial^kg(x)\,\frac{(y-x)^k}{k!}
$$
for the Taylor remainder function of $g$ at order $\alpha_2$. We use here the convention that for $z=(z^1,\dots,z^{d_0})\in\bbR^{d_0}$ and $k\in\bbN^{d_0}$ one sets $z^k=\prod_{1\leq i\leq d_0} (z^i)^{k_i}$. M. Hoshino extended in \cite{Hoshino2024} the expansion result \eqref{EqFirstLocalExpansionParaproduct} for ${\sf P}(f,g)$ to any $f\in C^{\alpha_1}, g\in C^{\alpha_2}$ for $\alpha_1,\alpha_2>0$ by proving amongst other things that
\begin{equation} \label{EqLocalExpansionParaproduct}
\bigg\vert {\sf P}(f,g)(y) - \sum_{\vert k\vert<\alpha_1+\alpha_2} \partial_\star^k{\sf P}(f,g)(x)\,\frac{(y-x)^k}{k!} - \sum_{\vert k\vert<\alpha_1} \partial^kf(x) R^{\alpha_2}(g)(y,x) \bigg\vert \lesssim \vert y-x\vert^{\alpha_1+\alpha_2}, 
\end{equation}
where the {\it generalized derivative}
$$
\partial_\star^k{\sf P}(f,g) \defeq \partial^k{\sf P}(f,g) - \underset{\vert k_1\vert<\alpha_1, \vert k_2\vert\geq\alpha_2}{\sum_{k_1+k_2=k}} \binom{k}{k_1} (\partial^{k_1}f)(\partial^{k_2}g)
$$
is indeed well-defined pointwise. The inequality \eqref{EqLocalExpansionParaproduct} provides a local description of the behaviour of ${\sf P}(f,g)$ around an arbitrary point $x$ to a precision $\vert y-x\vert^{\alpha_1+\alpha_2}$ when $\alpha_1,\alpha_2>0$. Hoshino was able to prove in \cite{Hoshino2024} a local expansion result for ${\sf P}(f_1,f_2,f_3)$ when $\alpha_1, \alpha_2,\alpha_3$ are all three positive. Theorem \ref{ThmMain1} below provides the most general extension of this result for some arbitrary iterated paraproducts ${\sf P}(f_1,\dots,f_n)$. In the particular case where the $f_k\in C^{\alpha_k}$ with $\alpha_k>0$ for all $1\leq k\leq n$, it implies that the function ${\sf P}(f_1,\dots,f_n)$ has a local description around an arbitrary point $x$ up to a precision $\vert y-x\vert^{\alpha_1+\cdots+\alpha_n}$. The statement of Theorem \ref{ThmMain1} does not require that all the $\alpha_k$ be positive and takes a very precise form. Not only does ${\sf P}(f_1,\dots,f_n)$ have a local expansion around any point $x$, but the functions whose values at $x$ give the coefficients of the expansion of ${\sf P}(f_1,\dots,f_n)$ also have some local expansion, to a lower precision though. The coefficients that appear in the latter expansion can also be expanded, to an even lower precision, and so on. A reader acquainted with regularity structures will recognize here the verbal description of a modelled distribution over a regularity structure. Theorem \ref{ThmMain1} states that a certain family of functions and distributions is a model over a particular regularity structure.

\subsection{Regularity structures associated with iterated paraproducts$\boldmath{.}$}
\label{IntroExpansionPPSystems}

\hspace{0.1cm} The reader will find in Appendix \ref{SectionBasicsRS} some basic facts about regularity structures. It suffices to mention here that they involve some pairs of vector spaces $(T,T^+)$ equipped with some algebraic structures 
$$
\Delta : T \rightarrow T\otimes T^+
$$
and 
$$
\Delta^+ : T^+ \rightarrow T^+\otimes T^+.
$$ 
{\it 1. The regularity structure.} We need some notations to introduce the structure that is involved in Theorem \ref{ThmMain1}. This structure involves a positive integer $n$ and a fixed tuple of real numbers 
$$
\alpha = (\alpha_1,\dots,\alpha_n).
$$ 
We use some {\color{blue} \textbf{blue bold}} letters $\bbfk=(k_1,\dots,k_c)\in(\bbN^{d_0})^c$ to denote some tuples of multi-indices $k_i\in\bbN^{d_0}$ of arbitrary length $c$. Denote by $\vert k\vert=k^1+\cdots+k^d$ the $\ell^1(\bbN)$-norm of an arbitrary $k=(k^1,\dots,k^d)\in\bbN^{d_0}$, and set for $\bbfk=(k_1,\dots,k_c)\in(\bbN^{d_0})^c$
\begin{equation*} \begin{split}
\vert\bbfk\vert &\defeq \big(\vert k_1\vert, \cdots, \vert k_c\vert\big) \in\bbN^c   \\
\Vert\bbfk\Vert &\defeq \vert k_1\vert+ \cdots + \vert k_c\vert \in\bbN.
\end{split} \end{equation*}
For $k\in\bbN^{d_0}$ and a non-null integer $c$ we define the set $\mcP_c(k)$ of partitions of $k$ into $c$ sub-mutli-indices as
$$
\mcP_c(k) \defeq \Big\{(k_1,\dots,k_c)\in(\bbN^{d_0})^c \,;\,  k=k_1+\cdots+k_c \Big\}.
$$
One has
$$
\Vert\bbfk\Vert = \vert k\vert
$$
for any $\bbfk\in\mcP_c(k)$ where $k\in\bbN^{d_0}, c\geq 1$. For some integers $a<b$ we write $\llbracket a,b\rrbracket$ for the set of integers in the closed interval $[a,b]$. Let $X=(X^1,\dots,X^{d_0})$ stand for an abstract $d_0$-dimensional monomial with  commutative symbol coordinates. For $p\in\bbN^{d_0}$ we set 
$$
X^p \defeq (X^1)^{p_1}\cdots (X^{d_0})^{p_{d_0}}.
$$ 
Denote by $(\epsilon_1,\dots,\epsilon_{d_0})$ the canonical basis of $\bbN^{d_0}$, so $X^{\epsilon_i} = X^i$. The following symbols
$$
\mcB \defeq \Big\{  \llbracket a,b\rrbracket_{\bbfl} \, X^p \Big\}_{1\leq a<b\leq n ,\, \bbfl\in\mcP_{b-a}(\ell), \,\ell\in\bbN^{d_0},\, p\in \bbN^{d_0}} \cup \big\{X^p\big\}_{p\in\bbN^{d_0}}
$$
form the basis of a vector space denoted by $T$. Similarly the following symbols
$$
\mcB^+ \defeq \Big\{ \llbracket a,b\rrbracket_{\bbfl}^\bbfk \Big\}_{\textrm{\textbf{\textsf{condition}}}(a,b,\bbfk,\bbfl)}
\cup \big\{ X^{\epsilon_i} \big\}_{1\leq i\leq d} 
$$
generate freely an algebra with unit ${\bf 1}^+$ that we denote by $T^+$. One says that $(a,b,\bbfk,\bbfl)$ satisfies $\textrm{\textbf{\textsf{condition}}}(a,b,\bbfk,\bbfl)$ if $1\leq a<b\leq n, \bbfk=(k_a,\dots,k_b)\in\mcP_{b-a+1}(k)$ for some $k\in\bbN^{d_0}$, and $\bbfl\in\mcP_{b-a}(\ell)$ for some $\ell\in\bbN^{d_0}$, and we have
$$
\max(\vert k\vert, \vert\ell\vert) < \sum_{1\leq j\leq n} \vert \alpha_j\vert
$$
and
\begin{equation} \label{EqCondition}
\vert\ell\vert + \sum_{a\leq j\leq b} \alpha_j > \vert k\vert.
\end{equation}
We emphasize that the tuples $\bbfk=(k_a,\dots,k_b)\in\mcP_{b-a+1}(k)$ have $b-a+1$ components while the tuples $\bbfl\in\mcP_{b-a}(\ell)$ have $b-a$ components. (To have a unified picture in mind one can think of $\bbfl=(\ell_a,\dots,\ell_{b-1})$ as the tuple $(\ell_a,\dots,\ell_{b-1},0)$ with $b-a+1$ components.) The $k_i$ in $\bbfk$ will represent below some derivatives in some analytic expressions like \eqref{EqGMap} below. The $\ell_j$ in $\bbfl$ will represent some polynomial weights in some analytical expressions like \eqref{EqDefnDeltaiell} below. The symbols of $\mcB$ and $\mcB^+$ index some analytic quantities that will be described below. We define an $\alpha$-dependent grading on $T$ and $T^+$ by defining the degree of $ \llbracket a,b\rrbracket_{\bbfl} X^p \in\mcB$ as

$$
\big\vert \llbracket a,b\rrbracket_{\bbfl} \, X^p \big\vert_{\alpha} \defeq \Vert\bbfl\Vert + \sum_{a\leq j\leq b} \alpha_j + \vert p\vert,
$$
and, requiring that the degree map is multiplicative on $T^+$, we set $\vert \epsilon_i\vert_{\bm{\alpha}}=1$ and define the degree of $\llbracket a,b\rrbracket^{\bbfk}_{\bbfl} \in\mcB^+$ as 

$$
\big\vert \llbracket a,b\rrbracket^{\bbfk}_{\bbfl} \big\vert_{\alpha} \defeq \Vert\bbfl\Vert + \sum_{a\leq j\leq b} \alpha_j - \Vert \bbfk\Vert.
$$
We read on the condition \eqref{EqCondition} that the elements of $\mcB^+$ have a positive degree. We will see in Section \ref{RSIPP} that there are some particular splitting maps $\Delta^+$ and $\Delta$ that turn the pair
$$
\mathscr{T}_\alpha \defeq \big((T,\Delta) , (T^+,\Delta^+)\big)
$$
into a concrete regularity structure.

\ssk

\noindent {\it 2. A model on the regularity structure.} We now define the analytic objects $\sf \Pi$ and $\sf g$ that we associate to the symbols of the regularity structure. Jointly, they define a model ${\sf (\Pi,g)}$ over a truncated version of $\mathscr{T}_\alpha$ that is parametrized by some non-null integer $n$ and some distributions $(f_1,\dots, f_n)$, where $f_i\in C^{\alpha_i}$ for some regularity exponents $\alpha_i\in\bbR$. We make the following assumption on these exponents.

\medskip

\noindent \textit{\textbf{Assumption (A) --} One has $\sum_{a\leq j\leq b} \alpha_j\notin \bbZ$ for all $1\leq a\leq b\leq n$.}

\medskip

For $\ell\in\bbN^{d_0}$ and $i\geq -1$ we define the modified Littlewood-Paley projector $\Delta_i^\ell$ by setting
\begin{equation} \label{EqDefnDeltaiell}
(\Delta_i^\ell f)(x) \defeq f\big((\cdot-x)^\ell K_i(\cdot-x)\big)
\end{equation}
for all $f\in\mcD'(\bbR^d)$ and $x\in\bbR^d$, where $\Delta_i^0=\Delta_i$. For $j\geq 0$ we define 
$$
\Delta_{<j}^\ell \defeq \sum_{-1\leq j'\leq j-1}\Delta_{j'}^\ell
$$ 
and set
$$
{\sf P}_\ell(f,g) \defeq \sum_{i\geq 1} \big(\Delta_{<i-1}^\ell f\big)\,(\Delta_ig)
$$
for any $f,g\in\mcD'(\bbR^d)$. For $\bbfl=(\ell_1,\dots,\ell_c)\in(\bbN^{d_0})^c$ and $\bbfl_{\leq c-1}=(\ell_1,\dots,\ell_{c-1})\in(\bbN^{d_0})^{c-1}$ we define recursively
$$
{\sf P}_{\bbfl}(f_1,\dots,f_c) \defeq {\sf P}_{\ell_c}\big({\sf P}_{\bbfl_{\leq c-1}}(f_1,\dots, f_{c-1}) , f_c\big).
$$
For $ \llbracket a,b\rrbracket_{\bbfl} X^p \in\mcB$ we define the distribution ${\sf \Pi}\big( \llbracket a,b\rrbracket_{\bbfl} X^p \big)$ by its action on a test function $\varphi$
\begin{equation} \label{EqPiMap1}
{\sf \Pi} \big(\llbracket a,b\rrbracket_{\bbfl} X^p\big) (\varphi) = {\sf \Pi}\big(\llbracket a,b\rrbracket_{\bbfl}\big)(\cdot^p\varphi)
\end{equation}
with $(\cdot^p\varphi)(y)=y^p\varphi(y)$ and
\begin{equation} \label{EqPiMap2}
{\sf \Pi}\big(\llbracket a,b\rrbracket_{\bbfl}\big) \defeq {\sf P}_{\bbfl}\big(f_a,\dots,f_b\big).
\end{equation}
The definition of the character $\sf g$ on $T^+$ requires a notation. For a tuple $\beta=(\beta_1,\dots,\beta_c)\in\bbR^c$ of regularity exponents and $\bbfl=(\ell_1,\dots,\ell_{c-1})\in(\bbN^{d_0})^{c-1}$ we set $\ell_c=0$ and define the set of $\bbfl$-admissible cuts of $\beta$ as
\begin{equation} \label{EqDefnEllCut}
\bbfl-\textsf{Cut}(\beta) \defeq \bigg\{1\leq d\leq c-1 \,;\, \ell_d=0, \, \sum_{1\leq e\leq d}\big(\beta_e + \vert \ell_e\vert\big) > 0,\; \sum_{d+1\leq e\leq c}\big(\beta_e + \vert \ell_e\vert\big) < 0 \bigg\}
\end{equation}
and for $d\in \bbfl-\textsf{Cut}(\beta)$ we set
$$
r_d = r_d(\beta,\bbfl) \defeq \min \bigg\{\sum_{1\leq e\leq d}\big(\beta_e + \vert \ell_e\vert\big) , \; - \hspace{-0.2cm}\sum_{d+1\leq e\leq c}\big(\beta_e + \vert \ell_e\vert\big) \bigg\}.
$$
Set 
$$
\beta_{\leq e} \defeq (\beta_1,\dots,\beta_e),   \qquad   \beta_{>e} \defeq (\beta_{e+1},\dots,\beta_c), \qquad \beta_{\llbracket a,b\rrbracket} \defeq (\beta_a,\dots,\beta_b)
$$ 
for any $1\leq e\leq c$ and $e\leq c-1$ and $1\leq a\leq b\leq c$, respectively. We define recursively
\begin{equation} \begin{split}
\widetilde{\sf P}^{\beta}_{\bbfl} &(g_1,\dots,g_c)   \\
&\defeq {\sf P}_{\bbfl}(g_1,\dots,g_c)   \\
&\qquad- \underset{\vert m\vert<r_d}{\sum_{d\in \bbfl-\textsf{Cut}(\bm{\beta})}} \underset{\bbfm'\in\mcP_{c-d}(m)}{\sum_{\bbfm\in \mcP_d(m)}}  \frac{m!}{\bbfm!\, \bbfm'!} \, \widetilde{\sf P}_{\bbfl_{\leq d}}^{\beta_{\leq d}-\bbfm}\big(\partial^{m_1}g_1,\dots, \partial^{m_d}g_d\big) \, \widetilde{\sf P}_{\bbfl_{>d}+\bbfmp}^{\beta_{>d}}\big(g_{d+1},\dots,g_c\big)   
\end{split} \end{equation}
where $m\in\bbN^{d_0}$ and with the convention that $\widetilde{\sf P}_m^{\beta_c}\big(g_c\big)=g_c$. For any $\beta_i\in\bbR$ we denote by $C^{\beta_i}_\circ$ the closure of $C^\infty\cap C^{\beta_i}$ in $C^{\beta_i}$. In the course of proving Theorem \ref{ThmMain1} below we will prove that $\widetilde{\sf P}^{\beta}_{\bbfl} (g_1,\dots,g_c)\in L^\infty$ if $g_i\in C^{\beta_i}_\circ$ for all $1\leq i\leq c$ and $\Vert\bbfl\Vert + \sum_{1\leq i\leq c} \beta_i >0$. Given a tuple $\alpha=(\alpha_1,\dots,\alpha_n)\in\bbR^n$ of regularity exponents and $f_i\in C^{\alpha_i}_\circ$, for all $1\leq i\leq n$, we can then define for $\llbracket a,b\rrbracket_{\bbfk}^{\bbfl} \in \mcB^+$ with $\bbfk=(k_a,\dots,k_b)$
\begin{equation} \label{EqGMap}
{\sf g}\big(\llbracket a,b\rrbracket^{\bbfk}_{\bbfl}\big) \defeq \widetilde{\sf P}^{\alpha_{\llbracket a,b\rrbracket}-\vert\bbfk\vert}_{\bbfl}\big(\partial^{k_a}f_a,\dots, \partial^{k_b}f_b\big).
\end{equation}

\ssk

\begin{thm} \label{ThmMain1}
The pair ${\sf (\Pi,g)}$ is a model on the regularity structure $\mathscr{T}_{\alpha}$. It depends continuously on $(f_1,\dots,f_n)\in \prod_{i=1}^n C^{\alpha_i}_\circ$.
\end{thm}

\ssk

For $\sf g$ to be part of a model we need to prove that each function ${\sf g}\big(\llbracket a,b\rrbracket^{\bbfk}_{\bbfl}\big)$ has a local expansion to accuracy $\vert y-x\vert^{\vert \llbracket a,b\rrbracket^{\bbfk}_{\bbfl} \vert_\alpha}$ around any point $x$, with the different terms in the expansion indexed by the algebraic structure of the Hopf algebra $(T^+,\Delta^+)$. For $\sf \Pi$ to be part of a model it also needs to satisfy some local expansion property that involves $\sf g$ as well. The strategy that we adopt to prove Theorem \ref{ThmMain1} is first to prove a statement of a similar flavor for some distributions and functions that are built from a simplified version of the iterated paraproducts. The algebra involved in the analysis of these operators is simpler than that of the true iterated paraproducts, and their analytical properties are more flexible. At the same time, we will see in Proposition \ref{PropoRepresentationProposition} of Section \ref{SectionRepresentationFormula} that ${\sf P}(f_1,\dots,f_n)$ can be written as a sum of simplified paraproducts evaluated on some other functions/distributions built from the $f_i$. This fact will play a crucial role in transfering the local expansion properties of the simplified iterated paraproducts to the true iterated paraproducts.

\medskip

\subsection{Local expansion properties of paracontrolled systems$\boldmath{.}$ \hspace{0.1cm}}
\label{IntroUniversalStructure}

\hspace{0.1cm} Paracontrolled systems are, within paracontrolled calculus, the equivalent of modelled distributions. Assume we are given a finite set of letters $\mcL=\{l_1,\dots,l_{\vert \mcL\vert}\}$ and a family $[l]\in C^{r_l}$ of distributions on $\bbR^{d_0}$ indexed by $\mcL$. We denote by $w=l_{i_1}\dots l_{i_w}$ a generic word with letters from $\mcL$. The concatenation of two words $w_1$ and $w_2$ is denoted by $w_1w_2$. If $w=w_1w_2$ we say that $w_1$ the a {\it begining of the word} $w$. We assume that the letters come with a notion of size $\vert l_i\vert\in\bbR$ and set 
$$
\vert w\vert \defeq \vert l_{i_1}\vert+\cdots+\vert l_{i_w}\vert.
$$ 
The empty word $w_\emptyset$ has size $0$. For a positive real number $r$ we denote by $\mathcal{W}_{<r}$ the set of words of size less than $r$, including the empty word. An $r$-\textbf{\textit{paracontrolled system}} is a family $(u_w)_{w\in \mcU_{<r}}$ of functions/distributions on $\bbR^{d_0}$ indexed by a subset $\mcU_{<r}$ of $\mathcal{W}_{<r}$ that contains the empty word $w_\emptyset$ and which has the following properties. 
\begin{itemize}
	\item[(1)] There is a {\it finite} subset $\mcU_{<r}^f$ of $\mcU_{<r}$ made up of words of positive size and such that every word of $\mcU_{<r}$ is the begining of one of the words of $\mcU_{<r}^f$. (The exponent $f$ in $\mcU_{<r}^f$ stands for `{\sl final}'.)   \vspace{0.1cm}
	
	\item[(2)] For all $w\in\mcU_{<r}$ one has
\begin{equation} \label{EqPCSystem}
u_w = \sum_{l\in \mcL} {\sf P}_{u_{wl}}[l] + u_w^\sharp
\end{equation}
with $u_w^\sharp\in C^{r-\vert w\vert}$.
\end{itemize}

Condition (1) ensures that the family $(u_w)_{w\in \mcU_{<r}}$ is finite even if some of the sizes $\vert l\vert$ are non-positive. This condition is automatically satisfied if all the $\vert l\vert$ are positive. We talk of the $[l]$ as the {\sl reference functions/distributions}. We see from this description that each $u_w$ writes as a finite sum of iterated paraproducts of the form ${\sf P}\big(u_{w_\bullet^\sharp}, [l_{i_1}], \dots, [l_{i_\bullet}]\big)$, including $u_w^\sharp={\sf P}(u_w^\sharp)$. 

The notion of paracontrolled system is useful even for the study of regularity structures. Bailleul \& Hoshino proved for instance in \cite{PCandRS1} that, for a model $\sf M=(\Pi,g)$ on a fixed regularity structure, the distributions/functions ${\sf \Pi}(\tau)$ and ${\sf g}(\mu)$ can be described by some paracontrolled systems 
\begin{equation} \label{EqPCSForModel} \begin{split}
{\sf \Pi}(\tau) &= \sum_{\sigma<\tau} {\sf P}_{{\sf g}(\tau/\sigma)}[ \sigma]^{\sf M} + [ \tau] ^{\sf M}  \\
{\sf g}(\mu) &= \sum_{1^+<^+\nu<^+\mu} {\sf P}_{{\sf g}(\mu/\nu)}[\nu]^{\sf g} + [ \mu]^{\sf g}
\end{split} \end{equation}
for some reference functions/distributions $[ \tau]^{\sf M} \in C^{\vert\tau\vert}, [ \mu]^{\sf g}\in C^{\vert\mu\vert}$ built from $\sf M$, and that for any modelled distribution $\bm{v}=\sum_\tau v_\tau \tau$ of positive regularity $r$ the family $\big({\sf R}^{\sf M}(\bm{v}),(v_\tau)_\tau\big)$ is an $r$-paracontrolled system 
\begin{equation} \label{EqPCSForModelledDistributions} \begin{split}
{\sf R}^{\sf M}(\bm{v}) &= \sum_{\vert\tau\vert<r} {\sf P}_{v_\tau} [\tau]^{\sf M} + [ \bm{v}]   \\
v_\tau &= \sum_{\tau<\sigma, \vert\sigma\vert<r} {\sf P}_{v_\sigma} [\sigma/\tau]^{\sf g} + [ v_\tau]
\end{split} \end{equation}
with reference functions/distributions the family of brackets $[ \tau]^{\sf M}, [\mu]^{\sf g}$. (This is Proposition 12 and Theorem 1 in \cite{PCandRS1}.)  Bailleul \& Hoshino further proved in Theorem 1 of \cite{PCandRS2} that a sub-family of these `brackets' $[ \tau] ^{\sf M}, [ \mu]^{\sf g}$ parametrizes the set of models, providing in particular a linear parametrization of the nonlinear space of models. These results hold for any reasonable regularity structure. For a particular class of regularity structures $\mathscr{T}$ including the BHZ regularity structures used for the study of subcritical singular stochastic PDEs, they proved that for a given model on $\mathscr{T}$ the set of admissible modelled distributions with regularity $r$ is parametrized by the family of functions/distributions $\big\{[ \bm{v}]\in C^r\big\} \cup \big\{[ v_\tau] \in C^{r-\vert\tau\vert},\, \textrm{$\tau$ in a linear basis of }T\big\}_{\vert\tau\vert<r}$ -- this is Theorem 5 and Theorem 7 in \cite{PCandRS2}. In all these results the regularity structure is fixed. In particular, if we are given some placeholders for $[ \bm{v}]$ and the $[ v_\tau]$ there is a unique admissible modelled distribution {\it over the given regularity structure} that has these functions/distributions as its brackets.

In the more general situation of an arbitrary paracontrolled system there is no a priori regularity structure. Our second main result means informally that we can lift a paracontrolled system into a modelled distribution on some universal regularity structure and for some system-dependent model. Recall we assume $[l]\in C^{r_l}$.

\ssk

\begin{thm} \label{ThmMain2}
Pick $r>0$. Given an $r$-paracontrolled system $(u_w)_{w\in\mathcal{U}_{<r}}$ as in \eqref{EqPCSystem} there is an explicit regularity structure $\mathscr{T}_\mcL$ that depends {\it only} on $\vert\mcL\vert, \,r$ and the $r_l$, a model $\sf M$ on $\mathscr{T}_\mcL$ and a modelled distribution $\bm{u}$ of regularity $r$ such that $u_{w_{\emptyset}} = {\sf R}^{\sf M}(\bm{u})$.
\end{thm}

\medskip

Regularity structures were first introduced by M. Hairer as a setting adapted to give sense to, and study, a large class of stochastic partial differential equations that are outside of the scope of classical stochastic calculus. Each equation in this class can be formulated as a fixed point problem in a random space of modelled distributions over a deterministic, equation-dependent, regularity structure. A solution to a singular stochastic partial differential equation then comes under the form of a local expansion around each state space point. The choice of regularity structures as a language to make sense of solutions to such equations is not the only choice possible. Gubinelli, Imkeller \& Perkowski \cite{GIP} introduced the notion of paracontrolled distribution at the same time that Hairer introduced regularity structures. A number of singular equations can be dealt with using the setting of the high order paracontrolled calculus of Bailleul \& Bernicot \cite{BailleulBernicotHighOrder}. In that language a solution to an equation comes under the form of a paracontrolled system. It was thus natural to ask if there is any correspondence between the fundamental notions of models and modelled distributions in regularity structures and the basic notion of paracontrolled system in paracontrolled calculus. Bailleul \& Hoshino's works \cite{PCandRS1, PCandRS2} establish such a correspondence building on \eqref{EqPCSForModel} and \eqref{EqPCSForModelledDistributions}. One associates to an equation a regularity structure $\mathscr{T}$, to a model $\sf M$ the paracontrolled system \eqref{EqPCSForModel} and to a model distribution $\bm{v}$ defined of $\sf M$ the paracontrolled system \eqref{EqPCSForModelledDistributions}. The inverse map consists in getting back the model $\sf M$ over $\mathscr{T}$ from \eqref{EqPCSForModel} and the modelled distribution $\bm{v}$ from \eqref{EqPCSForModelledDistributions}. In that context, Theorem \ref{ThmMain2} does something of a different nature. Starting from the paracontrolled systems \eqref{EqPCSForModel} and \eqref{EqPCSForModelledDistributions}, it introduces 
\begin{itemize}
	\item[--] another regularity structure $\mathscr{T}_\mcL$ that retains little information about the initial regularity structure $\mathscr{T}$,
	\item[--] a model and a modelled distribution on $\mathscr{T}_\mcL$,
\end{itemize}
whose associated paracontrolled systems are also given by \eqref{EqPCSForModel} and \eqref{EqPCSForModelledDistributions}. This situation is somewhat reminiscent of the study by Hairer \& Kelly \cite{HK} of the links between the notions of geometric and branched rough paths.
	
\medskip

\textbf{\textit{Organisation of the article.}} A `simplified' iterated paraproduct operator ${\sf P}_<(f_1,\dots,f_n)$ is introduced in Section \ref{SectionSimplifiedIPP}, and we provide in Section \ref{SectionExpansionSimplifiedIPP} its local expansion properties. The latter involve some functions $\partial_\star^k{\sf P}(f_1,\dots,f_n)$ that are introduced in Section \ref{SectionDStarK}. These functions also have some local expansion properties which we investigate in Section \ref{SectionLocalExpansionDStarK}. We leave aside the simplified iterated paraproducts in Section \ref{RSIPP} and describe in this section the regularity structure $\mathscr{T}_\alpha$ that is involved in the statement of Theorem \ref{ThmMain1}. This statement is proved in Section \ref{SectionProofThmsMain}. We build in Section \ref{SectionBuildingBlocks} a number of functions/distributions that will be used to represent an iterated paraproduct ${\sf P}(f_1,\dots,f_n)$ as a sum of simplified ${\sf P}_<$ iterated paraproducts. The representation formula itself is proved in Section \ref{SectionRepresentationFormula}. We prove Theorem \ref{ThmMain1} in Section \ref{SectionProofTheorems}. Section \ref{SectionPCSystems} is dedicated to proving Theorem \ref{ThmMain2}. We describe the universal regularity structure involved in this statement in Section \ref{SectionRSPCSystems} and prove Theorem \ref{ThmMain2} in Section \ref{SectionPCSystemsModelledDistributions}. A number of technical lemmas are deferred to some appendices. The proof of the local expansion property of the $\partial_\star^k{\sf P}(f_1,\dots,f_n)$ involves in particular some algebraic results that are proved in Appendix \ref{SectionAppendixAlgebraicLemmas}. So is the proof of some algebraic identities that play a crucial role in our proof of Theorem \ref{ThmMain1}. Appendix \ref{SectionBasicsRS} gives some background on regularity structures and Appendix \ref{SectionBasicsAnalysis} gives some general and particular analysis results.

\bigskip

\textit{\textbf{Notations.}} {\it We collect here a number of notations that are used throughout the text.

\begin{itemize}
	\item[--] The letters $i,j$ and $a,b,c,d,e$ will exclusively be used to denote some integers.
	
	\item[--] The letters $k,\ell,m$ will denote exclusively some elements of $\bbN^{d_0}$.
	
	\item[--] We denote by $\alpha=(\alpha_1,\alpha_2,\dots)$ or $\beta=(\beta_1,\beta_2,\dots)$ some finite tuples of regularity exponents $\alpha_i,\beta_j$ in $\bbR$.
	
	\item[--] For $z=(z^1,\dots,z^{d_0})\in\bbR^{d_0}$ and $k\in\bbN^{d_0}$ we write $z^k=\prod_{1\leq i\leq d_0} (z^i)^{k_i}$.
	
	\item[--] For $k=(k^1,\dots,k^{d_0})\in\bbN^{d_0}$ we write $k!=\prod_{i=1}^{d_0}k^i!$. For $m,m'_1,\dots,m'_r$ in $\bbN^{d_0}$ we set
	$$
	\binom{m}{(m'_1,\dots,m'_r)} \defeq \frac{m!}{\prod_{1\leq i\leq r}m'_i!}.
	$$
	
	\item[--] We write $\lesssim_p$ for an inequality that holds up to a multiplicative positive constant that only depends on a parameter $p$.
	
	\item[--] We work here in the Euclidean space $\bbR^{d_0}$. All that follows has a direct counterpart in an anisotropic version of $\bbR^{d_0}$. We stick to the Euclidean setting not to distract the reader from the main points of this work. 
\end{itemize}
}

\medskip

\section{Simplified iterated paraproducts and their local expansion properties} 
\label{SectionSimplifiedIPP}

We introduce in this section some simplified iterated paraproducts defined inductively from their Littlewood-Paley description
$$
\Delta_i\big({\sf P}_<(f_1,\dots,f_n)\big) = \Delta_{<i-1}\big({\sf P}_<(f_1,\dots,f_{n-1})\big)\,\Delta_i(f_n).
$$
It turns out to be convenient to define these functions/distributions on a slightly larger class of objects than the usual $C^\alpha$ spaces, $\alpha\in\bbR$. The setting is described in Section \ref{SectionSimplifiedPP}. The description of the local expansion properties of the simplified iterated paraproducts involves some generalized derivative operators $\partial_\star^k$ that we introduce in Section \ref{SectionDStarK}. We state and prove the local expansion property of the ${\sf P}_<({\sf f}_1,\dots, {\sf f}_n)$ in Section \ref{SectionExpansionSimplifiedIPP}, for ${\sf f}_j\in{\sf C}^{\alpha_j}$ with $\alpha_j>0$ for all $1\leq j\leq n$.

\subsection{Simplified iterated paraproducts.}
\label{SectionSimplifiedPP}

\hspace{0.1cm} We will work through part of this document with the following extension of the H\"older-Besov spaces.

\ssk

\begin{defn*}
For $r\in\bbR$ we define ${\sf C}^r$ as the vector space of sequences ${\sf{f}} = (f_i)_{i\geq -1}$ of smooth functions to which one can associate a ball $B\subset(\bbR^{d_0})'$ such that each $f_i$ is spectrally supported in $2^i B$ and
$$
\norme{ \sf{f}}_r \defeq \sup_{i\geq -1} 2^{ir}\norme{f_i}_{L^\infty} < \infty.
$$
This formula defines a norm on ${\sf C}^r$. An element of 
$$
{\sf C}^\infty \defeq \bigcap_{r>0} {\sf C}^r
$$
is said to be {\sl smooth}, and we set 
$$
{\sf C}^{-\infty} \defeq \bigcup_{r\in\bbR} {\sf C}^r, \qquad {\sf C}^{0+} \defeq \bigcup_{r>0}{\sf C}^r.
$$
We write ${\sf C}^r_\circ$ for the closure of ${\sf C}^\infty $ in ${\sf C}^r$.\end{defn*}

\ssk

For $r>0$ there is a canonical continuous non-injective surjection from ${\sf{C}}^r$ onto $ C^r$ sending $\sf f$ to $\sum_{i\geq -1} f_i$. The Paley-Littlewood projectors give a continuous injection from $C^r$ into $\sf C^r$ for any $r\in\bbR$. We define for any distribution $f$ on $\bbR^{d_0}$ and $o \in \bbR_+$ its Taylor polynomial $T_h^o f$ of order $o$ in the direction $h\in\bbR^d$ as the distribution
$$
(T_h^o f)(\cdot) \defeq \sum_{|k|<o} \frac{h^k}{k!} \, (\partial^k f)(\cdot).
$$
Its associated Taylor remainder $R_h^o f$ is defined from the relation
$$
|h|^o (R_h^o f)(\cdot) \defeq f(\cdot+h) - (T_h^o f)(\cdot).
$$
The derivation operator $\partial^k$, the Taylor expansion and remainder maps $T^o_h, R^o_h$, can be applied  to any ${\sf f}= (f_i)_{i\geq -1}\in{\sf C}^r$ by applying the corresponding classical operators to every $f_i$. These operations behave well in this context. For any $r\in\bbR$ and $k\in\bbN^{d_0}$, Bernstein inequalities ensures that $\partial^k  : {\sf C}^r \rightarrow {\sf{C}}^{r-|k|}$ defines a continuous operator. We give the proof of the following elementary fact in Appendix \ref{SectionBasicsAnalysis}.

\ssk

\begin{lem} \label{LemResteTaylor}
For $r\in\bbR, {\sf f}\in {\sf C}^r$ with $f_i$ is spectrally supported in $2^iB$, and $o\in \bbR_+$ we have
$$
{\sf f}(\cdot+h) - \sum_{|k|<o}\frac{h^k}{k!} (\partial^k {\sf f})(\cdot) = |h|^o \, (R_h^o{\sf f})(\cdot)
$$
with 
$$ 
\norme{R_h^o{\sf f}}_{r-o} \lesssim_B \norme{{\sf f}}_r,
$$
uniformly over $\vert h\vert\leq 1$.
\end{lem}

\ssk

For ${\sf{f}}_1,\dots, {\sf{f}}_n$ in ${\sf C}^{-\infty}$ we define iteratively the simplified iterated paraproducts 
$$
\textbf{\textsf{P}}_<({\sf{f}}_1,\dots,{\sf{f}}_n) = \big( {\sf{P}_<}({\sf{f}}_1,\dots,{\sf{f}}_n)_i \big)_{i\geq -1}
$$ 
as the element of ${\sf C}^{-\infty}$ given by ${\sf{P}_<}({\sf f}_1) =  {\sf f}_1$ and with ${\sf f}_n=(f_{ni})_{i\geq -1}$
$$
{\sf{P}_<}({\sf{f}}_1,\dots,{\sf{f}}_n)_i \defeq \Delta_{<i-1}\big({\sf {P}_<}({\sf{f}}_1,\dots,{\sf{f}}_{n-1})\big) \, f_{ni}.
$$
We write
$$
{\sf{P}_<}({\sf f}_1,\dots,{\sf f}_n) \defeq \sum_{i\geq -1} {\sf{P}}_<({\sf f}_1,\dots,{\sf f}_{n})_i
$$
for its corresponding distribution. Recall from \S 2 of Section \ref{IntroExpansionPPSystems} the statement of \textbf{\textit{Assumption (A)}}. From now on

\medskip

 {\it all our tuples $\alpha=(\alpha_1,\dots,\alpha_n), \beta=(\beta_1,\dots,\beta_n)$ in $\bbR^n$ will satisfy \textbf{\textit{Assumption (A)}}. }

\medskip

\subsection{Generalized derivative operator $\partial_{\star}^k$.}
\label{SectionDStarK}

\hspace{0.1cm} Recall from \eqref{EqDefnEllCut} the definition of the set of the $\bbfl$-admissible cuts of a tuple $\beta=(\beta_1,\dots,\beta_n)\in\bbR^n$. We define here the set of \textbf{\textit{cuts of $\beta$}} as 
$$
{\sf Cut}(\beta) \defeq {\color{blue} \bm{0}}-{\sf Cut}(\beta) = \bigg\{ d \in [\![1,n-1]\!],\quad \sum_{j=1}^d \beta_j >0 \enskip \text{  and } \sum_{j=d+1}^n \beta_j <0 \bigg\}.
$$
We also define the following set of \textbf{\textit{multi-cuts of $\beta$}}
$$
{\sf MultiCut}(\beta) \defeq \Big\{ \bm{d} = \Big(0=d_0<d_1<\cdots<d_{n(\bm{d})} = n\Big) \,;\, \forall e\in\llbracket1,n(\bm{d})-1\rrbracket, \; d_e\in {\sf Cut}(\beta)  \Big\}.
$$

\ssk

\begin{defn}
For $\beta\in\bbR^n$ and ${\sf f}_1,\dots,{\sf f}_n \in {\sf C}^{0+}$ we set
\begin{align*}
\widetilde{{\sf P}}_<^{\beta}({\sf f}_1,\dots,{\sf f}_n) \defeq \sum_{\bm{d} \in {\sf MultiCut}(\beta)}(-1)^{n(\bm{d})+1} \prod_{e=1}^{n(\bm{d})}{\sf P}_<\big({\sf f}_{d_{e-1}+1},\dots,{\sf f}_{d_e}\big).
\end{align*} 
We also set for $i\geq -1$
\begin{equation*} \begin{split}
\widetilde{{\sf P}}_<^{\beta}({\sf h}_1,&\dots,{\sf h}_n)\{i\}   \\
&\defeq \sum_{\bm{d} \in {\sf MultiCut}(\beta)}(-1)^{n(\bm{d})+1} \bigg\{ \prod_{c=1}^{n(\bm{d})-1}{\sf P}_<\big({\sf h}_{d_{c-1}+1},\dots,{\sf h}_{d_c}\big) \bigg\} {\sf P}_<\big({\sf h}_{d_{n(\bm{d})-1}+1},\dots,{\sf h}_n\big)_{i}.
\end{split} \end{equation*}
\end{defn}

\ssk

\noindent One has for instance $ \widetilde{{\sf P}}_<^{\beta}({\sf f})= {\sf P}_<({\sf f})= \sum_{i\geq -1} f_i$ for all ${\sf f}=(f_i)_{i\geq -1}\in{\sf C}^{-\infty}$, and 
\begin{equation*} \begin{split}
&\widetilde{{\sf P}}_<^{(2,1)}({\sf f}_1,{\sf f}_2) = {\sf P}_<({\sf f}_1,{\sf f}_2),   \\
&\widetilde{{\sf P}}_<^{(1,-2)}({\sf f}_1,{\sf f}_2) = {\sf P}_<({\sf f}_1,{\sf f}_2) - {\sf f}_1 {\sf f}_2,   \\
&\widetilde{{\sf P}}_<^{(1,-2,2,-1)}({\sf f}_1,{\sf f}_2,{\sf f}_3,{\sf f}_4) = {\sf P}_<({\sf f}_1,{\sf f}_2,{\sf f}_3,{\sf f}_4) - {\sf f}_1 {\sf P}_<({\sf f}_2,{\sf f}_3,{\sf f}_4) - {\sf P}_<({\sf f}_1,{\sf f}_2,{\sf f}_3) {\sf f}_4+{\sf f}_1{\sf P}_<({\sf f}_2,{\sf f}_3) {\sf f}_4,   \\
&\widetilde{{\sf P}}_<^{(1,-1,3/2)}({\sf f}_1,{\sf f}_2,{\sf f}_3) = {\sf P}_<({\sf f}_1,{\sf f}_2,{\sf f}_3).
\end{split} \end{equation*}

One has the relation $\widetilde{{\sf P}}_<^{\beta}({\sf h}_1,\dots,{\sf h}_n)=\sum_{i\geq -1} \widetilde{{\sf P}}_<^{\beta}({\sf h}_1,\dots,{\sf h}_n)\{i\}$, however $\widetilde{{\sf P}}_<^{\beta}({\sf h}_1,\dots,{\sf h}_n)\{i\}$ does not represent the Paley-Littlewood projection of $\widetilde{{\sf P}}_<^{\beta}({\sf h}_1,\dots,{\sf h}_n)$ as it is not spectrally supported in a ball, we introduce it as it appears naturally in the algebraic manipulations involving $\widetilde{\sf P }_<^\beta$ operators.

\ssk

\begin{prop} \label{EstiContinuityPTilde}
For any $\beta\in\bbR^n$, setting
$$
\mcE_\beta \defeq \bigg\{ c\in\llbracket 1,n\rrbracket \,;\, \sum_{j=1}^{c-1}\beta_j>0 \;\text{ and }\; \sum_{j=c}^{n}\beta_j>0 \bigg\}
$$
and 
$$
m_0\defeq 
\left\{ \begin{array}{ll}
\max \mcE_\beta, & \text{if } \mcE_\beta\neq\emptyset,   \\
1,& \text{otherwise}.
\end{array}\right.,
$$
one has for every $\big({\sf f}_1,\dots,{\sf f}_n\big) \in ({\sf C}^{0+})^n$ the estimate
$$
\big| \widetilde{{\sf P}}_<^{\beta} \big( {\sf f}_1,\dots,{\sf f}_n\big)\{ i\}\big| \lesssim 2^{-i \sum_{j=m_0}^n\beta_j}\prod_{j=1}^n\norme{{\sf f}_j}_{\beta_j}.
$$
\end{prop}

\ssk

The exponent $\beta$ in the notation $\widetilde{\sf P}^\beta_<$ records this somewhat optimal extension result.

\ssk

\begin{cor} \label{PropContinuityPTilde}
For any $\beta\in\bbR^n$ such that $\sum_{j=1}^n \beta_j > 0$, the map
$$
\big({\sf f}_1,\dots,{\sf f}_n\big) \in ({\sf C}^{\infty})^n \mapsto \widetilde{{\sf P}}_<^{\beta} \big( {\sf f}_1,\dots,{\sf f}_n\big) \in L^\infty
$$ 
has a continuous extension as a map from $\prod_{j=1}^n {\sf{C}}^{\beta_j}_\circ$ into $L^\infty$.
\end{cor}

\ssk

The remaining of this section is dedicated to proving Proposition \ref{EstiContinuityPTilde}. We will use for that purpose the following algebraic result. We state it and prove it first before giving the proof of Proposition \ref{PropContinuityPTilde}. In the following statement any constant in the open interval $(1,2)$ could be used in place of $3/2$.

\ssk

\begin{lem}\label{lemsomitere}
Given $\beta\in\bbR^n$ we define
\begin{equation*} \begin{split}
\rho_c &\defeq 
\left\{ 
\begin{array}{ll} 
+1 & \emph{if }(n-c)\in {\sf Cut}(\beta)  \\
-1 & \emph{otherwise   }
\end{array} 
\right.,   \qquad 
\rho \defeq \prod_{c=1}^{n-1}(-\rho_c).
\end{split} \end{equation*}
For any ${\sf h}_1=(h_{1i})_{i\geq -1},\dots, {\sf h}_n=(h_{ni})_{i\geq -1}$ in ${\sf C}^{\infty}$ we have
$$
\widetilde{\sf P}_<^{\beta}({\sf h}_1,\dots,{\sf h}_n)\{ i_1\} = \rho \sum_{\rho_1(i_2-i_1+3/2)>0} \dots \sum_{\rho_{n-1} (i_{n}-i_{n-1}+3/2)>0} \enskip\prod_{c=1}^n h_{ci_{n-c+1}}.
$$
\end{lem}

\ssk

\begin{Dem}
We prove the identity by induction on $n$. The result holds for $\widetilde{\sf P}_<^{\beta}({\sf h}_1)$. Suppose now that it holds for $(n-1)$ functions and consider first the case that $(n-1)\notin {\sf Cut}(\beta)$, so $\rho_1=-1$ and the condition $\rho_1(i_2-i_1+3/2)>0$ reads $i_2<i_1-1$. Then $\widetilde{{\sf P}}_<^{\beta}\big({\sf h}_1,\dots,{\sf h}_n\big)\{i_1\}$ is equal to
$$
\sum_{\bm{d} \in {\sf MultiCut}(\beta)}  \hspace{-0.1cm} (-1)^{n(\bm{d})+1}  \bigg\{ \prod_{c=1}^{n(\bm{d})-1}{\sf P}_<\big({\sf h}_{d_{c-1}+1},\dots,{\sf h}_{d_c}\big) \bigg\} {\sf P}_<\big({\sf h}_{d_{n(\bm{d})-1}+1},\dots,{\sf h}_{n}\big)_{i_1}
$$

\begin{equation*} \begin{split}
&=\sum_{\bm{d} \in {\sf MultiCut}(\beta)} \hspace{-0.1cm} (-1)^{n(\bm{d})+1}  \bigg\{ \prod_{c=1}^{n(\bm{d})-1} {\sf P}_<\big({\sf h}_{d_{c-1}+1},\dots,{\sf h}_{d_c}\big) \bigg\} {\sf P}_<\big({\sf h}_{d_{n(\bm{d})-1}+1},\dots,{\sf h}_{n-1}\big)_{<i_1-1} \, h_{ni_1}
\\
&= \sum_{i_2<i_1-1} \widetilde{{\sf P}}_<^{\beta^*}\big({\sf h}_1,\dots,{\sf h}_{n-1}\big)\{i_2\} \, h_{ni_1}
\end{split} \end{equation*}
where
$$
\beta^*\defeq \Big(\beta_1, \, \dots,\,\beta_{n-2}, \, \beta_{n-1}+\beta_n\Big).
$$
From the induction hypothesis we have 
$$
\widetilde{{\sf P}}_<^{\beta^*}\big({\sf h}_1,\dots,{\sf h}_{n-1}\big)\{i_2\} =  \rho \sum_{\rho_2(i_3-i_2+3/2)>0} \cdots \sum_{\rho_{n-1} (i_{n}-i_{n-1}+3/2)>0} \enskip\prod_{c=1}^{n-1} h_{ci_{n-c+1}},
$$
so we can conclude the induction in that case. If now $(n-1)\in {\sf Cut}(\beta)$ we have $\rho_1=1$ and the condition $\rho_1(i_2-i_1+3/2)>0$ reads $i_2\geq i_1-1$. We have in that case
\begin{align*}
\widetilde{{\sf P}}_<^{\beta}&({\sf h}_1,\dots,{\sf h}_n)\{i_1\}   \\
&=\sum_{\substack{\bm{d} \in {\sf MultiCut}(\beta) \\  (n-1)\in \bm{d} }}(-1)^{n(\bm{d})+1} \prod_{c=1}^{n(\bm{d})-1}{\sf P}_<\big({\sf h}_{d_{c-1}+1},\dots,{\sf h}_{d_c}\big) h_{ni_1}
\\
&\enskip+\sum_{\substack{\bm{d} \in {\sf MultiCut}(\beta) \\  (n-1)\notin \bm{d} }}(-1)^{n(\bm{d})+1} \prod_{c=1}^{n(\bm{d})-1}{\sf P}_<\big({\sf h}_{d_{c-1}+1},\dots,{\sf h}_{d_c}\big){\sf P}_<\big({\sf h}_{d_{n(\bm{d})-1}+1},\dots,{\sf h}_{{n}}\big)_{i_1}
\\
&=\sum_{\substack{\bm{d} \in {\sf MultiCut}(\beta) \\  (n-1)\in\bm{d} }}(-1)^{n(\bm{d})+1} \prod_{c=1}^{n(\bm{d})-1}{\sf P}_<\big({\sf h}_{d_{c-1}+1},\dots,{\sf h}_{d_c}\big) h_{ni_1}
\\
&\enskip+\sum_{\substack{\bm{d} \in {\sf MultiCut}(\beta) \\  (n-1)\notin \bm{d} }}(-1)^{n(\bm{d})+1} \prod_{c=1}^{n(\bm{d})-1}{\sf P}_<\big({\sf h}_{d_{c-1}+1},\dots,{\sf h}_{d_c}\big){\sf P}_<\big({\sf h}_{d_{n(\bm{d})-1}+1},\dots,{\sf h}_{n-1}\big)_{<i_1-1} h_{ni_1}
\\
&= \sum_{\bm{d}\in {\sf MultiCut}(\beta^*)} (-1)^{n(\bm{d})+1} \bigg\{  \prod_{c=1}^{n(\bm{d})-1}{\sf P}_<\big({\sf h}_{d_{c-1}+1},\dots,{\sf h}_{d_c}\big)
\\
&\hspace{3cm} - \prod_{c=1}^{n(\bm{d})-2}{\sf P}_<\big({\sf h}_{d_{c-1}+1},\dots,{\sf h}_{d_c}\big){\sf P}_<\big({\sf h}_{{ d_{n(\bm{d})-2}+1}},\dots,{\sf h}_{n-1}\big)_{<i_1-1}   \bigg\}  h_{ni_1}
\\
&= - \sum_{i_2>i_1-2} \widetilde{{\sf P}}_<^{\beta^*}\big({\sf h}_1,\dots,{\sf h}_{n-1}\big)\{i_2\} \, h_{ni_1}.
\end{align*}
We conclude from the induction hypothesis that
$$
\widetilde{{\sf P}}_<^{\beta}\big({\sf h}_1,\dots,{\sf h}_n\big)\{i_1\} = -\sum_{i_2 \geq i_1-1}  (-\rho) \sum_{\rho_2(i_3-i_2+3/2)>0} \cdots \hspace{-0.15cm} \sum_{\rho_{n-1} (i_{n}-i_{n-1}+3/2)>0}  \bigg\{ \prod_{c=1}^{n-1} h_{ci_{n-c+1}} \bigg\} \, h_{ni_1}
$$
which allows us to close the induction in that case.
\end{Dem}

\medskip

\begin{Dem}[of Proposition \ref{PropContinuityPTilde}] 
For ${\sf h}_1,\dots {\sf h}_n \in \sf C^{+\infty}$ we have from Lemma \ref{lemsomitere} the bound 
$$
\big|\widetilde{\sf P}_<^{\beta}\big({\sf h}_1,\dots,{\sf h}_n\big)\{i\}\big| \lesssim C_\beta(i) \prod_{j=1}^n \norme{{\sf h}_j}_{\beta_j},
$$
where 

$$
C_\beta(i_1) \defeq \sum_{\rho_1(i_2-i_1+3/2)>0} \cdots \hspace{-0.2cm}\sum_{\rho_{n-1} (i_{n}-i_{n-1}+3/2)>0} \; \prod_{c=1}^n 2^{-i_{n+c-1}\beta_c}.
$$
We prove by induction that
\begin{equation} \label{estimeesommeitere}
C_\beta(i) \lesssim 2^{-i \sum_{j=m_0}^n\beta_j}.
\end{equation}
\begin{itemize}
	\item[--] If $\beta_1<0$ we have $\rho_{n-1}=-1$ and 
$$
\sum_{\substack{i_n;\enskip \rho_{n-1} (i_{n}-i_{n-1}+3/2)>0}} \enskip 2^{-i_{n}\beta_1} \simeq 2^{-i_{n-1}\beta_1}.
$$
We have in that case
$$
C_{(\beta_1,\dots,\beta_n)}(i) \simeq C_{(\beta_1+\beta_2,\beta_3\dots,\beta_n)}(i).
$$
	\item[--] If now $\beta_1>0$ and $\sum_{j=2}^n\beta_j<0$, then $\rho_{n-1}=+1$ and we have
$$
\sum_{\substack{i_n;\enskip \rho_{n-1} (i_{n}-i_{n-1}+3/2)>0}} \enskip 2^{-i_{n}\beta_1} \simeq 2^{-i_{n-1}\beta_1},
$$
so we have again
$$
C_{(\beta_1,\dots,\beta_n)}(i) \simeq C_{(\beta_1+\beta_2,\beta_3\dots,\beta_n)}(i).
$$
	\item[--] If finally $\beta_1>0$ and $\sum_{j=2}^n\beta_j>0$, we have this time
$$
\sum_{\substack{i_n;\enskip \rho_{n-1} (i_{n}-i_{n-1}+3/2)>0}} \enskip 2^{-i_{n}\beta_1} \simeq 1,
$$
so
$$
C_{(\beta_1,\dots,\beta_n)}(i) \simeq C_{(\beta_2,\beta_3\dots,\beta_n)}(i).
$$
\end{itemize}
In all the cases the inequality \eqref{estimeesommeitere} follows by induction since $C_{(\beta_1)}(i)=2^{-i\beta_1}$.
\end{Dem}

\ssk

\begin{defn*}
Pick some integers $1\leq a\leq b\leq n$ and $\alpha=(\alpha_a,\dots,\alpha_b)\in\bbR^{b-a+1}$. For $\bbfk=(k_a,\dots,k_b)\in(\bbN^{d_0})^{b-a+1}$ and ${\sf f}_a,\dots,{\sf f}_b$ in ${\sf C}^\infty$ we define
$$
\partial^{\bbfk}_{\star\alpha}{\sf P}_<\big({\sf f}_a,\dots,{\sf f}_b\big) \defeq \widetilde{{\sf P}}_<^{\alpha_{[a,b]}-\vert\bbfk\vert} \Big(\partial^{k_a}{\sf f}_a,\dots, \partial^{k_b}{\sf f}_b \Big).
$$
and
$$
\partial_{\star\alpha}^k {\sf P}_<\big({\sf f}_a,\dots,{\sf f}_b\big) \defeq \sum_{\bbfk \in \mcP_{b-a+1}(k)} \binom{k}{\bbfk} \, \partial_{\star\alpha}^{\bbfk}{\sf P}_<\big({\sf f}_a,\dots,{\sf f}_b\big),
$$
\end{defn*}

\ssk

As a consequence of Corollary \ref{PropContinuityPTilde} the map $\partial_{\star\alpha}^k{\sf P}_<$ is continuous from $\prod_{j=a}^b {\sf C}_\circ^{\alpha_j}$ into $L^\infty$ if $\vert k\vert<\sum_{j=a}^b \alpha_j$. It makes sense in that setting so simply write $\partial_\star^k$ rather than $\partial_{\star\alpha}^k$, as the information on $\alpha$ is already recorded in the domain $\prod_{j=a}^b {\sf C}_\circ^{\alpha_j}$ of the extension. 

The following lemma gives a recursive definition of the $\widetilde{\sf P}_<^{\beta}({\sf h}_1,\dots,{\sf h}_n)\{i\}$.

\ssk

\begin{lem} \label{rectilde1}
For any $\beta = (\beta_1,\dots,\beta_n)\in\bbR^n$ and any ${\sf h}_1,\dots,{\sf h}_n$ in ${\sf C}^\infty$ we have
\begin{align} \label{EqInductiveTildePMinus}
\widetilde{{\sf P}}_<^{\beta}\big({\sf h}_1,\dots {\sf h}_n\big)\{i\} &= {\sf P}_<\big({\sf h}_1,\dots,{\sf h}_n\big)_i 
\\
&\qquad- \sum_{d\in {\sf Cut}(\beta)}  \widetilde{\sf P}_<^{(\beta_1,\dots,\beta_d)}\big({\sf h}_1,\dots, {\sf h}_d\big) \, \widetilde{{\sf P}}_<^{(\beta_{d+1},\dots,\beta_n)}\big({\sf h}_{d+1},\dots, {\sf h}_n\big)\{i\}.
\end{align}
\end{lem}

\begin{Dem}
\textit{\textbf{Assumption (A)}} implies in particular that the $\sum_{c=1}^j\beta_c$ are all distinct for different $j\in \llbracket 1,n-1 \rrbracket$. We then have the following partition of ${\sf MultiCut}(\beta)$
$$
{\sf MultiCut}(\beta) = \big\{(0,n)\big\} \sqcup \bigsqcup_{d\in {\sf Cut}(\beta)} {\sf MultiCut}(\beta)[d],
$$
with
$$
{\sf MultiCut}(\beta)[d] \defeq \bigg\{ \bm{d} \in {\sf MultiCut}(\beta) \,;\, d\in \bm{d}, \quad \sum_{c=1}^d\beta_c = \min_{j\in\mathbf{d} } \sum_{c=1}^j\beta_c \bigg\}.
$$
One can thus write
\begin{align*}
&\widetilde{\sf P}_<^{\beta}({\sf h}_1,\dots {\sf h}_n)\{i\} = {\sf P}_<({\sf h}_1,\dots,{\sf h}_n)_i 
\\
&\qquad+ \sum_{d\in {\sf Cut}(\beta)}  \sum_{\bm{d} \in {\sf MultiCut}(\beta)[d]}(-1)^{n(\bm{d})+1} \prod_{c=1}^{n(\bm{d})-1}{\sf P}_<\big({\sf h}_{d_{c-1}+1},\dots,{\sf h}_{d_c}\big)\, {\sf P}_<\big({\sf h}_{d_{n(\bm{d})-1}+1},\dots,{\sf h}_{n}\big)\{i\}.
\end{align*}
For $d\in {\sf Cut}(\beta)$ and $1<j<d$ we have the equivalence
$$
\Big(\exists \, \bm{d} \in {\sf MultiCut}(\beta)[d],\enskip j\in \bm{d}\Big) \Leftrightarrow \Big(j \in {\sf Cut}\big( (\beta_1,\dots,\beta_d)\big)\Big).
$$
Likewise for $d<j<n$ we have 
$$
\Big(\exists \, \bm{d} \in {\sf MultiCut}(\beta)[d],\enskip j\in \bm{d} \Big) \Leftrightarrow \Big( j-d \in {\sf Cut}\big((\beta_{d+1},\dots,\beta_n)\big)\Big).
$$ 
This entails that we have
\begin{align*}
 &\sum_{\bm{d} \in {\sf MultiCut}(\beta)[d]}(-1)^{n(\bm{d})+1} \prod_{c=1}^{n(\bm{d})-1}{\sf P}_<\big({\sf h}_{d_{c-1}+1},\dots,{\sf h}_{d_c}\big) {\sf P}_<\big({\sf h}_{d_{n(\bm{d})-1}+1},\dots,{\sf h}_{n}\big)\{i\}
 \\
 &= - \widetilde{\sf P}_<^{\beta_{\leq d}}\big({\sf h}_1,\dots, {\sf h}_d\big) \, \widetilde{\sf P}_<^{\beta_{>d}}\big({\sf h}_{d+1},\dots, {\sf h}_n\big)\{i\},
\end{align*}
from which the statement of  the lemma follows.
\end{Dem}
    
    
On can rewrite Lemma \ref{rectilde1} in the context of the $\partial_\star$-derivatives.  For any multi-indice $k\in\bfN^{d_0}$ we have
\begin{equation} \label{eq_recstarder}
\begin{split}
\partial_{\star\alpha}^k &{\sf P}_<\big({\sf f}_1,\dots,{\sf f}_n\big)   \\ 
&=\partial^k {\sf P}_<\big({\sf f}_1,\dots,{\sf f}_n\big) - \sum_{c=1}^{n-1}\sum_{\substack{|\ell |<\sum_{j=1}^c\alpha_j \\ |k-\ell |>\sum_{j=c+1}^n\alpha_j }} \binom{k}{\ell} \, \partial_{\star\alpha_{\leq c}}^\ell{\sf P}_<\big({\sf f}_1,\dots,{\sf f}_c\big) \, \partial^{k-\ell}_{\star\alpha_{>c}}{\sf P}_<\big({\sf f}_{c+1},\dots,{\sf f}_n\big).
\end{split}
\end{equation}

\subsection{Local expansion properties of the ${\sf P}_<({\sf f}_1,\dots {\sf f}_n)$.}
\label{SectionExpansionSimplifiedIPP}

\hspace{0.1cm} Recall Hoshino's expansion result \eqref{EqLocalExpansionParaproduct} for ${\sf P}(f,g)$, for both $f$ and $g$ of positive regularity. We would like to give a similar expansion result for ${\sf P}_<(f_1,\dots,f_n)$ for an arbitrary $n\geq 2$.

Pick ${\sf f}_j\in{\sf C}^{\alpha_j}$ for each $1\leq j\leq n$. Let us make a first naive try at expanding ${\sf P}_<({\sf f}_1,\dots {\sf f}_n)(\cdot+h)$ as a function of $h\in\bbR^{d_0}$. For any $o>0$ we have 
\begin{equation} \label{EqDev1}  \begin{split}
     {\sf P}_<\big( &{\sf f}_1,\dots, {\sf f}_n \big)(\cdot+h) = {\sf P}_<\big( {\sf f}_1(\cdot + h),\dots, {\sf f}_n(\cdot+h)\big)
     \\
     &={\sf P}_<\bigg( \sum_{|k_1|<o} \frac{h^{k_1}}{k_1!} \, \partial^{k_1}{\sf f_1} + |h|^o R_h^o{\sf f_1}, \, {\sf f}_2(\cdot+h),\dots \bigg)
     \\
     &= \sum_{|k_1|<o} {\sf P}_<\bigg( \partial^{k_1}{\sf f_1}\,\frac{h^{k_1}}{k_1 !}, \; \sum_{|k_2|<o-|k_1|} \frac{h^{k_2}}{k_2 !} \, \partial^{k_2}{\sf f}_2 + |h|^{o-|k_1|} R^{o-|k_1|}_{h}{\sf f}_2 ,\dots \bigg)
     \\
     &\hspace{2cm}+{\sf P}_<\Big( |h|^o R_h^o{\sf f}_1, \, {\sf f}_2(\cdot+h),\dots \Big) = (\cdots) 
     \\
     &= \sum_{|k|<o}\sum_{\bbfk \in \mcP_n(k)}\frac{h^k}{k!} \binom{k}{\bbfk} \, {\sf P}_<\big(\partial^{k_1}{\sf f}_1,\dots,\partial^{k_n}{\sf f}_n\big)
     \\
     &\qquad+ \sum_{c=1}^{n} \sum_{\substack{|k|<o \\ \bbfk \in \mcP_{c-1}(k) } } \frac{h^k|h|^{o-|k|}}{\bbfk !} \, {\sf P}_<\Big( \partial^{k_1}{\sf f}_1, \dots,\partial^{k_{c-1}}{\sf f}_{c-1}, R^{o-|k|}_h {\sf f}_c, \, {\sf f}_{c+1}(\cdot+h), \, \dots \Big) 
\end{split} \end{equation}
\begin{equation*} \begin{split}
      &= T_h^o {\sf P}_< \big({{\sf f}_1},\dots,{{\sf f}_n}\big)
     \\
     &\qquad+ \sum_{c=1}^{n} \sum_{\substack{|k|<o \\ \bbfk \in \mcP_{c-1}(k) } } \frac{h^k |h|^{o-|k|}}{\bbfk !} \, {\sf P}_<\Big( \partial^{k_1}{\sf f}_1, \dots, \partial^{k_{c-1}}{\sf f}_{c-1}, R^{o-|k|}_h {\sf f}_c, \,  {\sf f}_{c+1}(\cdot+h),\,\dots \Big). 
\end{split} \end{equation*}
This formula does not give us the kind of expansion we are looking for as the last paraproducts in the right hand side of the equation contain some distributions with negative regularities so these paraproducts have no reason to define some functions. This would be the case if we had instead of some ${\sf P}_<$ terms some $\widetilde{\sf P}^\beta_<$ terms, for some appropriate tuples $\beta$ depending on the arguments. We will get our local expansion for ${\sf P}_<({\sf f}_1,\dots {\sf f}_n)(\cdot+h)$ by introducing the appropriate terms to force the appearance of these $\widetilde{\sf P}^\beta_<$. We proceed gradually and first introduce the quantity that will be the remainder term in this expansion. For $1\leq a\leq b\leq n-1, k\in\bbN^{d_0}$ and $\bbfk=(k_{a+1},\dots,k_{b-1}) \in \mcP_{b-a-1}(k)$ set
$$
\alpha_a( \bbfk,o) \defeq \Big(\alpha_{a+1}-|k_{a+1}|,\dots,\,\alpha_{b-1}-|k_{b-1}|,\, \alpha_b-o+|k|,\, \alpha_{b+1},\dots,\,\alpha_n \Big).
$$
and
\begin{align*}
\big(\triangle_{h,o}^\alpha {\sf P}_<\big)({\sf f}_{a+1},\dots,{\sf f}_n) \defeq \sum_{b=a+1}^{n} \hspace{-0.5cm}\sum_{\substack{|k|<o \\ \bbfk \in \mcP_{b-a-1}(k) } }\hspace{-0.5cm} \frac{h^k|h|^{o-|k|}}{\bbfk !} \, &\widetilde{\sf P}_<^{\alpha_a( \bbfk,o)}\Big( \partial^{k_{a+1}}{\sf f}_{a+1}\dots,\partial^{k_{b-1}}{\sf f}_{b-1},   \\
&\quad R^{o-|k|}_h {\sf f}_b, \, {\sf f}_{b+1}(\cdot+h),\dots,{\sf f}_n(\cdot+h) \Big),
\end{align*}
and for $i\geq-1$
\begin{align*}
\big(\triangle_{h,o}^\alpha {\sf P}_<\big)({\sf f}_{a+1},\dots,{\sf f}_n)\{i\} \defeq \sum_{b=a+1}^{n} \hspace{-0.5cm}\sum_{\substack{|k|<o \\ \bbfk \in \mcP_{b-a-1}(k) } }\hspace{-0.5cm} \frac{h^k|h|^{o-|k|}}{\bbfk !} \, &\widetilde{\sf P}_<^{\alpha_a( \bbfk,o)}\Big( \partial^{k_{a+1}}{\sf f}_{a+1}\dots,\partial^{k_{b-1}}{\sf f}_{b-1},   \\
&\quad R^{o-|k|}_h {\sf f}_b, \, {\sf f}_{b+1}(\cdot+h),\dots,{\sf f}_n(\cdot+h) \Big)\{i\}.
\end{align*}
We denote by $\delta_0$ the distance from $\bbZ$ to the set of all $\sum_{a\leq j\leq b} \alpha_j\notin \bbZ$ where $1\leq a\leq b\leq n$; it is positive from \textbf{\textit{Assumption (A)}}. Lemma \ref{LemResteTaylor} and Proposition \ref{EstiContinuityPTilde} give us uniform  continuity estimates on $(\triangle_{h,o}^\alpha {\sf P}_<)({\sf f}_{a+1},\dots,{\sf f}_n)\{i\}$. In particular if $o>\sum_{j=a+1}^n \alpha_j -\delta_0$, one has 
\begin{equation} \label{EqEstimateRemainderPMinus}
\big| \big(\triangle_{h,o}^\alpha {\sf P}_{<}\big) ({\sf f}_{a+1},\dots,{\sf f}_n)\{ i\} \big| \lesssim |h|^o 2^{-i(\sum_{j=a}^n\alpha_j - o)} \prod_{j=a+1}^n \norme{ {\sf f}_j}_{\alpha_j},
\end{equation}
and for $o<\sum_{j=a+1}^n \alpha_j$
$$
\big| \big(\triangle_{h,o}^\alpha {\sf P}_{<}\big) ({\sf f}_{a+1},\dots,{\sf f}_n)\big| \lesssim |h|^o \prod_{j=a+1}^n \norme{ {\sf f}_j}_{\alpha_j}.
$$

\ssk

\begin{prop} \label{PropAlgDevSimpl}
Pick ${\sf f}_1,\dots,{\sf f}_n$ in ${\sf C}^\infty$. Assume all the $\alpha_j$ are positive and $o > \sum_{j=1}^n \alpha_j - \delta_0$. Then we have 
\begin{equation} \label{devsimpparap} \begin{split}
\big(\triangle_{h,o}^\alpha {\sf P}_<\big)({\sf f}_1,\dots, {\sf f}_n)\{i\}  &=  {\sf P}_<\big({\sf f}_1,\dots,{\sf f}_n\big)_i(\cdot+h) - T_h^o {\sf P}_<\big({\sf f}_1,\dots,{\sf f}_n\big)_i
     \\
&\qquad - \sum_{a=1}^n\sum_{|k|<\sum_{i=1}^a\alpha_i}\partial_{\star\alpha_{\leq a}}^k{\sf P}_<({\sf f}_1,\dots,{\sf f}_a) \, \frac{h^k}{k!} \, \big(\triangle_{h,o-|k|}^\alpha {\sf P}_<\big) ({\sf f}_{a+1}, \dots, {\sf f}_n)\{i\}.
\end{split} \end{equation}
\end{prop}

\ssk

\begin{Dem}
We use in the proof the shorthand notation
$$
\alpha(\bbfk) \defeq \alpha_0(\bbfk,o) = \Big(\alpha_{1}-|k_1|, \dots,\, \alpha_{j-1}-|k_{j-1}|, \, \alpha_j-o+|k|, \, \alpha_{j+1},\dots, \alpha_n \Big).
$$
As all the $\alpha_j$ are positve we have ${\sf Cut}(\alpha(\bbfk))\subset \llbracket1,j-1\rrbracket$, so \eqref{EqInductiveTildePMinus} writes here

\begin{align*}
    \widetilde{\sf P}_<^{\alpha(\bbfk)} &\Big( \partial^{k_1}{\sf f}_1,\dots,\partial^{k_{j-1}}{\sf f}_{j-1}, \, R^{o-|k|}_h{\sf f}_j, \, {\sf f}_{j+1}(\cdot+h),\dots \Big) \{ i\}
    \\   
    &= {\sf P}_<\Big( \partial^{k_1} {\sf f}_1,\dots,\partial^{k_{j-1}}{\sf f}_{j-1}, R^{o-|k|}_h{\sf f}_j, \, {\sf f}_{j+1}(\cdot+h),\dots \Big)_i
    \\[0.2em]
    &\quad- \sum_{d\in {\sf Cut}(\alpha(\bbfk))} \widetilde{\sf P}_<^{\alpha(\bbfk)_{\leq d}}\Big( \partial^{k_1}{\sf f_1},\dots,\partial^{k_{d}}{\sf f}_d\Big) 
    \\[0.1em]
    &\hspace{3.2cm}\times\widetilde{\sf P}_<^{\alpha(\bbfk)_{>d}}\Big( \partial^{k_{d+1}}{\sf f}_{d+1},\dots,\partial^{k_{j-1}}{\sf f}_{j-1}, \; R^{o-|k|}_h{\sf f}_j, \; {\sf f}_{j+1}(\cdot+h),\dots \Big)  \{ i\}
\end{align*} 
\begin{align*}
     &= {\sf P}_<\Big( \partial^{k_1}{\sf f}_1,\dots,\partial^{k_{j-1}}{\sf f}_{j-1}, R^{o-|k|}_h{\sf f}_j, \, {\sf f}_{j+1}(\cdot+h), \, \dots \Big)_i
    \\[0.2em]
    &\quad-\hspace{-0.2cm}\sum_{d\in {\sf Cut}(\alpha(\bbfk) )}\hspace{-0.2cm} \partial^{\bbfk_{\leq d}}_\star{\sf P}_<\Big({\sf f}_1,\dots,{\sf f}_d\Big)
    \\[0.1em]
    &\hspace{3.2cm}\times
    \widetilde{\sf P}_<^{\alpha(\bbfk)_{>d}}\Big( \partial^{k_{d+1}}{\sf f}_{d+1}, \dots, \partial^{k_{j-1}}{\sf f}_{j-1}, R^{o-|k|}_h{\sf f}_j, \, {\sf f}_{j+1}(\cdot+h),\,\dots \Big)\{i\}.
\end{align*}
Note that as $o > \sum_{j=1}^n \alpha_j -\delta_0$ we have 
$$
{\sf Cut}(\alpha(\bbfk)) = \bigg\{ d \in \llbracket1,n\rrbracket \,;\, \sum_{j=1}^d \alpha(\bbfk)_j >0 \bigg\};
$$
we will use this fact to invert the sums over $m$ and $j$ below. Summing over $j,k$ and $\bbfk$ gives
\begin{align*}
 & (\triangle_{h,o}^\alpha {\sf P}_<)({\sf f}_1,\dots,{\sf f}_n)\{i\}-\sum_{j=1}^{n} \sum_{\substack{|k|<o \\ \bbfk \in \mcP_{j-1}(k) } }  \hspace{-0.3cm}  \frac{h^k|h|^{o-|k|}}{\bbfk !}{\sf P}_<\Big( \partial^{k_1}{\sf f}_1,\dots,\partial^{k_{j-1}}{\sf f}_{j-1},\, R^{o-|k|}_h{\sf f}_{j},\, {\sf f}_{j+1}(\cdot+h),\dots \Big)_i
    \\
    &=\sum_{j=1}^{n} \sum_{\substack{|k|<o \\ \bbfk \in \mcP_{j-1}(k) } }\hspace{-0.3cm}\frac{h^k|h|^{o-|k|}}{\bbfk !} \hspace{-0.3cm} \sum_{d\in {\sf Cut}(\alpha(\bbfk))} \partial^{\bbfk_{\leq d}}_\star{\sf P}_<\big({\sf f}_1,\dots,{\sf f}_d\big)
    \\
    &\hspace{2.7cm}\times\widetilde{\sf P}_<^{\alpha(\bbfk)_{>d}} \Big( \partial^{k_{d+1}}{\sf f}_{d+1},\dots,\partial^{k_{j-1}}{\sf f}_{j-1}, \,  R^{o-|k|}_h{\sf f}_{j}, \, {\sf f}_{j+1}(\cdot+h),\dots \Big)\{ i\}
    \\
    &= \sum_{d=1}^{n} \sum_{\substack{|k|<\sum_{i=1}^d\alpha_i \\ \bbfk\in \mcP_{d-1}(k) }} \frac{h^{k}}{\bbfk !} \partial^{\bbfk}_{\star\alpha_{\leq d}}{\sf P}_<\big({\sf f}_1,\dots,{\sf f}_d\big) \sum_{j=d+1}^n\sum_{\substack{ |\ell|<o-|k|  \\  \bbfl \in \mcP_{j-d-1}(\ell)    }} \frac{h^\ell |h|^{o-|k|-|\ell|}}{\bbfl!}  
    \\
    &\hspace{2.7cm}\times\widetilde{\sf P}_<^{\beta_d( \bbfl,o-|k|)}\Big( \partial^{\ell_{1}}{\sf f}_{d+1},\dots,\partial^{\ell_{j-d-1}}{\sf f}_{j-1}, \enskip R^{o-|k|-|\ell|}_h{\sf f}_{j}, \, {\sf f}_{j+1}(\cdot+h),\dots \Big)\{ i\}
    \\
    &= \sum_{d=1}^n \sum_{\substack{|k|<\sum_{i=1}^d\alpha_i} }\frac{h^k}{k!} \, \partial_{\star\alpha_{\leq d}}^k{\sf P}_<({\sf f}_1,\dots,{\sf f}_d)\,  \big(\triangle_{h,o-|k|}^\alpha {\sf P}_<\big)({\sf f}_{d+1},\dots,{\sf f}_n)\{i\}.
\end{align*}
The identity \eqref{devsimpparap} then follows from \eqref{EqDev1}.
\end{Dem}

\ssk

The terms $\triangle_{h,o-|k|}^\alpha{\sf P}_<({\sf f}_{d+1},\dots,{\sf f}_n)$ for which $o-\vert k\vert > \sum_{j=d+1}^n \alpha_j$, in \eqref{devsimpparap}, are still problematic as one cannot use Corollary \ref{PropContinuityPTilde} for them. 

\ssk

\begin{lem}\label{lem_grostheta1}
Assume all the $\alpha_j$ positive. For $1\leq a\leq n$ and $\sum_{j=a}^n \alpha_j - \delta_0 < o_1 < o_2$, we have
$$
\big(\triangle_{h,o_2}^\alpha {\sf P}_<\big)\big({\sf f}_a,\dots,{\sf f}_n\big)\{ i\} - \big(\triangle_{h,o_1}^\alpha {\sf P}_<\big)\big({\sf f}_a,\dots,{\sf f}_n\big )  \{i\} = \sum_{r_1<|k|<o_2} \frac{h^k}{k!} \, \partial_{\star\alpha_{\geq a}}^k {\sf P}_<\big({\sf f}_a,\dots,{\sf f}_n\big)\{i\}.
$$
\end{lem}

\ssk

\begin{Dem}
We prove it by induction over $n-a$ with the help of Proposition \ref{PropAlgDevSimpl} and the inductive relation \eqref{eq_recstarder} satisfied by the star derivatives. The result is true for $a=n$ as in this case $\triangle_{h,o}{\sf P}_<$ coincides with the Taylor remainder $|h|^r R_h^r$. To run the induction step we use Proposition \ref{PropAlgDevSimpl}  to see that $\big(\triangle_{h,o_2}^\alpha {\sf P}_<\big)\big({\sf f}_a,\dots,{\sf f}_n\big)\{ i\} - \big(\triangle_{h,o_1}^\alpha {\sf P}_<\big)\big({\sf f}_a,\dots,{\sf f}_n\big)\{ i\}$ is equal to

\begin{align*}
&= T_h^{o_2}{\sf P}_<\big({\sf f}_a,\dots,{\sf f}_n\big)_i - T_h^{o_1}{\sf P}_<\big({\sf f}_a,\dots,{\sf f}_n\big)_i
\\
&\quad -\sum_{j=a}^{n-1} \sum_{|p|<\sum_{s=a}^j\alpha_s} \partial^p_{\star\alpha_{\llbracket a,j\rrbracket}}{\sf P}_<({\sf f}_a,\dots,{\sf f}_j) \, \frac{h^p}{p!}
\\
&\hspace{4.5cm}\times\Big\{ \big(\triangle_{h,o_2-p}^\alpha {\sf P}_<\big)\big({\sf f}_{j+1},\dots,{\sf f}_n\big)\{ i\} -\big(\triangle_{h,o_1-p}^\alpha {\sf P}_<\big)\big({\sf f}_{j+1},\dots,{\sf f}_n\big)\{ i\} \Big\}.
\end{align*}
From the induction hypothesis the above quantity is equal to 
\begin{align*}
\sum_{o_1<|k|<o_2} \partial^k{\sf P}_< ({\sf f}_a,\dots,{\sf f}_n)_i \, \frac{h^k}{k!}  - \sum_{j=a}^{n-1} &\sum_{|p|<\sum_{s=a}^j\alpha_s} \partial^p_{\star\alpha_{\llbracket a,j\rrbracket}}{\sf P}_<({\sf f}_a,\dots,{\sf f}_j) \, \frac{h^p}{p!}
\\ 
&\times \sum_{r_1<|\ell|+|p|<o_2} \partial^\ell_{\star\alpha_{\llbracket j+1,n\rrbracket}}{\sf P}_<\big({\sf f}_{j+1},\dots,{\sf f}_n\big)\{ i\} \,  \frac{h^\ell}{\ell!}.
\end{align*}
We conclude using \eqref{eq_recstarder}.
\end{Dem}

\bigbreak

For $0\leq c\leq n-1$ we let
$$
\triangle_{yx}{\sf P}_<({\sf f}_{c+1},\dots,{\sf f}_n) \defeq \big(\triangle_{y-x,\sum_{j=c+1}^n \alpha_j}^\alpha {\sf P}_< \big) \big({\sf f}_{c+1},\dots,{\sf f}_n \big)(x),
$$
and for $i\geq-1$
$$
\triangle_{yx}{\sf P}_<({\sf f}_{c+1},\dots,{\sf f}_n)\{i\} \defeq \big(\triangle_{y-x,\sum_{j=c+1}^n \alpha_j}^\alpha {\sf P}_< \big) \big({\sf f}_{c+1},\dots,{\sf f}_n \big)\{i\}(x).
$$

From Lemma \ref{lem_grostheta1} we know that for any $o\in( \sum_{j=c+1}^n \alpha_j-\delta_0,\sum_{j=c+1}^n \alpha_j+\delta_0)$ one has the equality
$$
\triangle_{yx}{\sf P}_<({\sf f}_{c+1},\dots,{\sf f}_n)\{i\} = \big(\triangle_{y-x,o}^\alpha {\sf P}_< \big) \big({\sf f}_{c+1},\dots,{\sf f}_n \big)\{i\}(x).
$$
Then for any $o$ in a neighborhood of $\sum_{j=c+1}^n\alpha_j$, from \eqref{EqEstimateRemainderPMinus} the following estimate holds
\begin{equation}\label{eq_estiremainfinal1}
\big| \triangle_{yx}{\sf P}_<({\sf f}_{c+1},\dots,{\sf f}_n)\{i\}\big| \lesssim |y-x|^o \prod_{j=c+1}^n \norme{{\sf f}_j}_{\alpha_j} \, 2^{-i( \sum_{j=c+1}^n \alpha_j-o) }
\end{equation}

The following Lemma was already used in \cite{Hoshinocommutator} and enables us to get the optimal bound on $| \triangle_{yx}{\sf P}_<({\sf f}_{c+1},\dots,{\sf f}_n)|$, we reproduce its proof in Appendix \ref{SectionBasicsAnalysis}.

\begin{lem}\label{lem_seuilopti}
Assume we are given a family of absolutely convergent series $\big(X_{yx}=\sum_{i\geq -1}X_{yx}^i\big)$ indexed by $x,y\in\bfR^d$, for which there exists some positive constants $C>0$ and $\gamma>0$ such that the uniform bound 
$$
|X_{yx}^i|\leq C 2^{-i(\gamma-\theta)} \, |y-x|^\theta
$$
holds for any $\theta$ in a neighborhood of $\gamma$. Then we have
$$
|X_{yx}|\lesssim C|y-x|^\gamma,
$$
uniformly over $x,y\in\bfR^d$ such that $|y-x|\leq 1$.
\end{lem}

From \eqref{eq_estiremainfinal1} and Lemma \ref{lem_seuilopti} one has then for $|y-x|\leq1$
$$
\big| \triangle_{yx}{\sf P}_<({\sf f}_{c+1},\dots,{\sf f}_n)\big| \lesssim \prod_{j=c+1}^n \norme{{\sf f}_j}_{\alpha_j} |y-x|^{\sum_{j=c+1}^n \alpha_j}
$$

\ssk

It follows then that the following fact holds.

\ssk

\begin{prop} \label{PropLocalDevptPMinus}
Pick $\alpha=(\alpha_1,\dots,\alpha_n)\in(0,+\infty)^n$, for all ${\sf f}_j\in {\sf C}^{\alpha_j}_0$, where $1\leq j\leq n$, we have the local expansion

\begin{align*}
{\sf P}_<({\sf f}_1,\dots,{\sf f}_n)(y) 
       &= \sum_{|k|<\sum_{j=1}^n\alpha_j} \partial_{\star\alpha}^k{\sf P}_<\big({\sf f}_1,\dots,{\sf f}_n\big)(x)\frac{(y-x)^k}{k!}
        \\
        &\quad+ \sum_{c=1}^{n-1}\sum_{|k|<\sum_{j=1}^c\alpha_j} \partial_{\star\alpha_{\leq c}}^k{\sf P}_<\big({\sf f}_1,\dots,{\sf f}_c\big)(x) \; \frac{(y-x)^k}{k!} \triangle_{yx}{\sf P}_<({\sf f}_{c+1},\dots,{\sf f}_n)
        \\
        &\quad+ \big(\triangle_{yx}{\sf P}_<\big)({\sf f}_1,\dots, {\sf f}_n),
\end{align*}
where 
$$
\Big| \big(\triangle_{yx}{\sf P}_<\big)({\sf f}_{c+1},\dots, {\sf f}_n)\Big| \lesssim \bigg\{\prod_{j=c+1}^n \norme{{\sf f}_j}_{\alpha_j}\bigg\} \, |y-x|^{\sum_{j=c+1}^n\alpha_j}.
$$
\end{prop}

\ssk

\begin{Dem}
    From Propositions \ref{PropAlgDevSimpl} and \ref{lem_grostheta1}, we have
    \begin{align*}
        (\triangle_{yx} {\sf P}_<)({\sf f}_1,\cdots, {\sf f}_n)
     &={\sf P}_<\big({\sf f}_1,\cdots,{\sf f}_n\big)(y) - \sum_{|k|<\theta}\frac{(y-x)^k}{k!}\partial^k{\sf P}_<\big({\sf f}_1,\cdots,{\sf f}_n\big)(x)
     \\
     &\quad-\sum_{c=1}^n\sum_{|p|<\sum_{j=1}^c\alpha_i}\partial^p_\star{\sf P}_<\big({\sf f}_1,\cdots,{\sf f}_c\big)(x) \frac{(y-x)^p}{p!}  \big( \triangle_{yx}{\sf P}_<\big)({\sf f}_{c+1},\cdots,{\sf f}_n)
     \\
     &\quad-\sum_{c=1}^n\sum_{\substack{|p|<\sum_{j=1}^c\alpha_j \\ |\ell|>\sum_{j=c+1}^n \alpha_j
 }} \frac{(y-x)^p}{p!} \frac{(y-x)^\ell}{\ell!}   \\
     &\hspace{4cm} \times \partial^p_\star{\sf P}_<\big({\sf f}_1,\cdots,{\sf f}_c\big)(x) \partial^{\ell}_\star {\sf P}_<\big({\sf f}_{c+1},\cdots,{\sf f}_n\big)(x).
    \end{align*}
Using Equation \ref{eq_recstarder} gives the statement of the proposition.
\end{Dem}

\ssk

We obtain the fact that one can work with ${\sf f}_j\in {\sf C}^{\alpha_j}_\circ$ rather than with ${\sf f}_j\in {\sf C}^\infty$ by an elementary continuity reasoning. In that setting we would write the expansion with the lighter notation $\partial_\star^k$ in place of $\partial_{\star\alpha_{\leq c}}^k$ as the use of the subscripts $\alpha_{\leq c}$ would be redundant.

\medskip

\section{Local expansion properties of the $\partial_\star^k{\sf P}_<(f_1,...,f_n)$}
\label{SectionLocalExpansionDStarK}

The quantities $\partial_{\star\alpha_{\leq c}}^k{\sf P}_<({\sf f}_1,\dots,{\sf f}_c)$ appear as coefficients in the local expansion of the simplified paraproduct ${\sf P}_<({\sf f}_1,\dots,{\sf f}_n)$. These coefficients also have a local expansion property that we describe in this section. It will be convenient for that purpose to introduce in Section \ref{Subsec_twisted2} some operators $\widetilde{\sf P}_<^{\beta^1\hspace{-0.05cm},\beta^2}$ indexed by two tuples of integers, as an intermediate tool. The local expansion formula for $\partial_\star^k{\sf P}_<(f_1,...,f_c)$ follows from a similar expansion for the $\widetilde{\sf P}_<^\alpha$ operators. The latter is given in Proposition \ref{PropExpansionTildePAlphaMinus} and takes a form similar to the expansion formula for ${\sf P}_<$. Its proof has the same architecture as the proof of Proposition \ref{PropLocalDevptPMinus}. It is described in Section \ref{subsection_devtwisted}.

\subsection{The operators $\widetilde{\sf P}_<^{\beta^1\hspace{-0.05cm},\beta^2}$.}  
\label{Subsec_twisted2}

\hspace{0.1cm} Their definition requires the following notation.

\ssk

\begin{defn*}
For $\beta^1, \beta^2$ in $\bbR^n$ such that $\beta^1_i\geq \beta^2_i$ for all $1\leq i\leq n$ we set 
\begin{equation*} \begin{split}
{\sf Multi}&{\sf Cut}\big(\beta^1, \beta^2\big)   \\
&\defeq \Big\{ \bm{d}=\big(0=d_0<d_1<\dots<d_{n(\bm{d})}=n\big) \,;\, \forall e\in\llbracket 1,{n(\bm{d})}-1\rrbracket, \, d_e\hspace{-0.05cm}\in \hspace{-0.05cm} {\sf Cut}(\beta^1)\cup {\sf Cut}( \beta^2) \Big\}.
\end{split} \end{equation*}
For $({\sf h}_i)_{1\leq i \leq n}\subset{\sf C}^{0+}$ we set
 $$
\widetilde{\sf P}^{\beta^1\hspace{-0.05cm},\beta^2}_<\big({\sf h}_{1},\dots, {\sf h}_n \big) \defeq \sum_{\bm{d} \in {\sf MultiCut}(\beta^1\hspace{-0.05cm},\beta^2)}(-1)^{n(\bm{d})+1} \prod_{e=1}^{n(\bm{d})}{\sf P}_<\big({\sf h}_{d_{e-1}+1},\dots,{\sf h}_{d_e}\big).
$$
\end{defn*}

\ssk

The following statement is proved in Appendix \ref{SectionAlgebraicProofs}.

\ssk

\begin{lem} \label{rectilde2}
Pick $\beta^1, \beta^2$ in $\bbR^d$ satisfying \textbf{\textit{Assumption (A)}}. If $\beta_i^2\leq \beta_i^1 $ for all $1\leq i\leq n$ then for any ${\sf h}_1,\dots,{\sf h}_n$ in ${\sf C}^\infty$ we have 
\begin{align*}
\widetilde{\sf P}^{\beta^1\hspace{-0.05cm},\beta^2}_< \big({\sf h}_1,\dots {\sf h}_n \big) \hspace{-0.04cm}= \widetilde{\sf P}^{\beta^1}_<\big({\sf h}_1,\dots,{\sf h}_n\big) - \hspace{-0.5cm} \sum_{d\in {\sf Cut}(\beta^2)\backslash {\sf Cut}(\beta^1)} \hspace{-0.5cm} \widetilde{\sf P}^{\beta^1_{\leq m},\beta^2_{\leq m}}_< \big({\sf h}_1,\dots, {\sf h}_m \big) \, \widetilde{\sf P}^{\beta^1_{>m},\beta^2_{> m}}_<\big({\sf h}_{m+1},\dots, {\sf h}_n\big).
\end{align*}
\end{lem}

\ssk

We will use in the end the $\widetilde{\sf P}^{\beta^1\hspace{-0.05cm},\beta^2}_<$ operators in settings where $\sum_{i=1}^n \beta_i^2>0$. In that case we have ${\sf Cut}(\beta^1)\subset {\sf Cut}(\beta^2)$, so $\widetilde{\sf P}_<^{\beta^1\hspace{-0.05cm},\beta^2}({\sf h}_{1},\dots, {\sf h}_n) = \widetilde{\sf P}_<^{\beta^2}({\sf h}_{1},\dots, {\sf h}_n)$, and we will be able to use the continuity property of Proposition \ref{PropContinuityPTilde}. General $\widetilde{\sf P}^{\beta^1\hspace{-0.05cm},\beta^2}_<$ operators will be useful in the algebraic steps.

\medskip

\subsection{Local expansion properties of the $\widetilde{\sf P}_<^{\mathbf{\beta}}\big({\sf f}_1,\dots,{\sf f}_n\big)$.} 
\label{subsection_devtwisted}

\hspace{0.1cm} Pick $\beta\in\bbR^n$ such that $\sum_{i=1}^n \beta_i >0$.  Proceeding as in \eqref{EqDev1} we see that
\begin{equation*} \begin{split}
     \widetilde{\sf P}_<^{\beta}\big({\sf f}_1 &, \dots, {\sf f}_n\big)(\cdot+h) = \widetilde{\sf P}_<^{\beta} \big( {\sf f}_1(\cdot + h),\dots, {\sf f}_n(\cdot+h)\big)   \\
      &=T_h^o \widetilde{\sf P}_<^{\beta}\big({\sf f}_1,\dots,{\sf f}_n\big)   \\
      &\qquad+ \sum_{m=1}^{n} \sum_{\substack{|k| < o \\ \bbfk \in \mcP_{m-1}(k) } } \frac{h^k|h|^{o-|k|}}{\bbfk !}\widetilde{\sf P}_<^{\beta} \Big( \partial^{k_1}{\sf f}_1,\dots,\partial^{k_{m-1}}{\sf f}_{m-1}, R^{o-|k|}_h{\sf f}_{m}, \, {\sf f}_{m+1}(\cdot+h), \,\dots \Big) 
 \end{split}  \end{equation*}
With the same motivations as in Section \ref{SectionExpansionSimplifiedIPP} we set here
$$
\beta_a( \bbfk,o)\defeq \Big(\beta_{a+1}-|k_{a+1}|,\dots,\,\beta_{b-1}-|k_{b-1}|,\, \beta_b-o+|k|,\, \beta_{b+1},\dots,\, \beta_n \Big).
$$
and
\begin{align*}
\big(\triangle_{h,o} \widetilde{\sf P}_<^{\beta}\big)({\sf f}_{a+1},\dots,{\sf f}_n) \defeq \sum_{b=a+1}^{n} \hspace{-0.3cm}\sum_{\substack{|k|<o \\ \bbfk \in \mcP_{b-a-1}(k) } }\hspace{-0.3cm} \frac{h^k|h|^{o-|k|}}{\bbfk !} \, &\widetilde{\sf P}_<^{\beta{>a}, \beta_a( \bbfk,r)}\Big( \partial^{k_{a+1}}{\sf f}_{a+1}\dots,\partial^{k_{b-1}}{\sf f}_{b-1},   \\
&\qquad R^{o-|k|}_h {\sf f}_b, \, {\sf f}_{b+1}(\cdot+h),\dots,{\sf f}_n(\cdot+h) \Big),
\end{align*}
and for $i\geq-1$
\begin{align*}
\big(\triangle_{h,o} \widetilde{\sf P}_<^{\beta}\big)({\sf f}_{a+1},\dots,{\sf f}_n)\{i\} \defeq \sum_{b=a+1}^{n} \hspace{-0.15cm}\sum_{\substack{|k|<o \\ \bbfk \in \mcP_{b-a-1}(k) } }\hspace{-0.15cm} \frac{h^k|h|^{o-|k|}}{\bbfk !} &\, \widetilde{\sf P}_<^{\beta{>a}, \beta_a( \bbfk,r)}\Big( \partial^{k_{a+1}}{\sf f}_{a+1}\dots,\partial^{k_{b-1}}{\sf f}_{b-1},   \\
& R^{o-|k|}_h {\sf f}_b, \, {\sf f}_{b+1}(\cdot+h),\dots,{\sf f}_n(\cdot+h) \Big)\{i\}.
\end{align*}
For any $k\in\bbN^{d_0}$ and $\bbfk=(k_a,\dots,k_b)\in\mcP_{b-a+1}(k)$ we also set
$$
\partial^{\bbfk}_{\star} \widetilde{\sf P}^{\beta_{\llbracket a,b\rrbracket}}_<({\sf f}_a,\dots,{\sf f}_b) \defeq \widetilde{\sf P}^{\beta_{\llbracket a,b\rrbracket}, \beta_{\llbracket a,b\rrbracket}-\vert\bbfk\vert}_<\big(\partial^{k_a}{\sf f}_a, \dots, \partial^{k_b}{\sf f}_b
\big)
$$
and
$$
\partial_\star^k \widetilde{\sf P}^{\beta_{\llbracket a,b\rrbracket}}_<({\sf f}_a,\dots,{\sf f}_b) \defeq \sum_{\bbfk\in\mcP_{b-a+1}(k)} \binom{k}{\bbfk} \partial^{\bbfk}_\star \widetilde{\sf P}^{\beta_{\llbracket a,b\rrbracket}}_<({\sf f}_a,\dots,{\sf f}_b).
$$
If $\sum_{i=a}^b \beta_i>0$ and $|k|<\sum_{i=a}^b \beta_i$, then $\widetilde{\sf P}^{\beta_{\llbracket a,b \rrbracket}, \beta_{\llbracket a,b \rrbracket}-\vert\bbfk\vert}_< = \widetilde{\sf P}^{\beta_{   \llbracket a,b \rrbracket } - \vert\bbfk\vert}_<$ and $\partial_\star^k \widetilde{\sf P}^{\beta_{\llbracket a,b\rrbracket}}_<({\sf f}_a,\dots,{\sf f}_b) = \partial_{\star \beta_{\llbracket a,b\rrbracket}}^k {\sf P}_<({\sf f}_a,\dots,{\sf f}_b)$ in that case.

\ssk

We define
$$
{\sf I}(\beta) \defeq \Big\{ c \in \llbracket 1,n-1\rrbracket \,;\, \sum_{j=1}^c \beta_j > 0 \; \emph{ and } \; \sum_{j=c+1}^n \beta_j > 0 \Big\}.
$$

\ssk

\begin{prop}\label{PropAlgDevSimpl2}
For $o > \sum_{j=1}^n \beta_j - \delta_0$, we have 
\begin{align*}
\big(\triangle_{h,r}^\beta \widetilde{\sf P}_<^\beta\big) ({\sf f}_1,\dots,{\sf f}_n)\{i\}
&=\widetilde{\sf P}_<^\beta \big({\sf f}_1,\dots,{\sf f}_n\big)\{i\}(\cdot+h) - T_h^o \widetilde{\sf P}^\beta_<\big({\sf f}_1, \dots, {\sf f}_n\big)\{i\}   \\
&\hspace{-0.5cm}- \sum_{c\in {\sf I}(\beta)} \sum_{|k|<\sum_{j=1}^c \beta_j}\partial^k_{\star\beta_{\leq c}}{\sf P}_<  \big({\sf f}_1,\dots,{\sf f}_c\big) \, \frac{h^k}{k!} \, \big(\triangle_{h,o-|k|}^\beta \widetilde{\sf P}_<^{\beta_{>c}}\big)\big({\sf f}_{c+1},\dots,{\sf f}_n\big)\{i\},
\end{align*}
\end{prop}

\ssk

\begin{Dem}
Recall that for $k\in\bbN^{d_0}$ and $\bbfk\in\mcP_{j-1}(k)$
$$
\beta(\bbfk,o) =  \Big( \beta_1 - |k_1|, \dots, \beta_{j-1} - |k_{j-1}|, \beta_j-o+|k| , \beta_{j+1}, \dots, \beta_n \Big).
$$
We are going to apply Lemma \ref{rectilde2} with the tuples $\alpha$ and $\beta(\bbfk,o)$, which verify indeed $\beta(\bbfk,o)_a\leq \alpha_a$ for any $1\leq a \leq n$. Moreover for $a>j$ we have $\beta(\bbfk,o)_a=\alpha_a$, and as consequence ${\sf Cut}(\beta(\bbfk,o))\backslash {\sf Cut}(\alpha)\subset \llbracket 1,j-1\rrbracket$. Then Lemma \ref{rectilde2} gives
\begin{align*}
    &\widetilde{\sf P}_<^{\beta(\bbfk,o)}\Big( \partial^{k_1}{\sf f}_1,\dots,\partial^{k_{j-1}}{\sf f}_{j-1}, \, R^{o-|k|}_h{\sf f}_{j}, \, {\sf f}_{j+1}(\cdot+h),\dots \Big) \{i\}
    \\[0.2em]
    &= \widetilde{\sf P}_<^{\beta}\Big( \partial^{k_1}{\sf f}_1,\dots,\partial^{k_{j-1}}{\sf f}_{j-1}, R^{o-|k|}_h{\sf f}_{j}, \, {\sf f}_{j+1}(\cdot+h),\dots \Big)\{i\}
    \\[0.2em]
    &\quad-\sum_{c\in {\sf Cut}(\beta(\bbfk,o)) \backslash {\sf Cut}(\alpha)} \widetilde{\sf P}_<^{ \alpha_{\leq c},\beta(\bbfk,o)_{\leq c}} \Big( \partial^{k_1}{\sf f}_1,\dots,\partial^{k_c}{\sf f}_c\Big) 
       \\[0.2em]
    &\hspace{3.5cm}\times \widetilde{\sf P}_<^{\beta_{>c},\beta(\bbfk,o)_{>c} }\Big( \partial^{k_{c+1}}{\sf f}_{c+1},\dots,  
\partial^{k_{j-1}}{\sf f}_{j-1}, R^{o-|k|}_h{\sf f}_{j}, \, {\sf f}_{j+1}(\cdot+h),\dots \Big) \{i\}
    \\[0.2em]  
     &= \widetilde{\sf P}^{\beta}_< \Big( \partial^{k_1}{\sf f}_1,\dots,\partial^{k_{j-1}}{\sf f}_{j-1}, \enskip R^{o-|k|}_h{\sf f}_{j}, \, {\sf f}_{j+1}(\cdot+h),\dots \Big)\{i\}
    \\[0.2em] 
    &\quad-\sum_{c\in {\sf Cut}(\beta(\bbfk,o))\backslash {\sf Cut}(\beta)} \partial^{\bbfk_{\leq c}}_\star{\sf P}_<\Big({\sf f}_1,\dots,{\sf f}_c\Big)
    \\[0.2em]
    &\hspace{3.5cm}\times\widetilde{\sf P}_<^{\alpha_{>c}, \beta(\bbfk,o)_{>c} }\Big( \partial^{k_{c+1}}{\sf f}_{c+1},\dots,\partial^{k_{j-1}}{\sf f}_{j-1}, R^{o-|k|}_h{\sf f}_j, \, {\sf f}_{j+1}(\cdot+h),\dots \Big)\{i\},
\end{align*}
where we used that 
$$
\widetilde{\sf P}_<^{\beta(\bbfk,o)_{\leq c},\beta_{\leq c}} \big( \partial^{k_1}{\sf f}_1,\dots, \partial^{k_c}{\sf f}_c\big) =\widetilde{\sf P}_<^{\beta(\bbfk,o)_{\leq c}}\big( \partial^{k_1}{\sf f}_1,\dots,\partial^{k_c}{\sf f}_c\big) 
$$ 
since $\sum_{i=1}^c \beta(\bbfk,o)_i>0$ for $c\in {\sf Cut}(\beta(\bbfk,o))$.

We now sum over $j,k$ and $\bbfk$ and invert the sums over $c$ and $j$. In order to implement this sum inversion we use the inclusion $ {\sf Cut}\big(\beta(\bbfk,o)\big)\backslash {\sf Cut}\big(\beta \big)\subset {\sf I}(\beta)$. This gives
\begin{align*}
\big(\triangle_{h,o} \widetilde{ \sf P}^{\beta}_<\big) &\big({\sf f}_1,\dots,{\sf f}_n\big)\{i\}-\sum_{j=1}^{n} \sum_{\substack{|k|<o \\ \bbfk \in \mcP_{j-1}(k) } }\frac{h^k|h|^{o-|k|}}{\bbfk !} \, \widetilde{\sf  P}_<^{\beta}  \Big( \partial^{k_1}{\sf f}_1,\dots,\partial^{k_{j-1}}{\sf f}_{j-1}, 
 \\
 &\hspace{8.5cm}R^{o-|k|}_h{\sf f}_j, \, {\sf f}_{j+1}(\cdot+h), \, \dots \Big)\{i\} 
 \\
    &=\sum_{j=1}^{n} \sum_{\substack{|k|<o \\ \bbfk \in \mcP_{j-1}(k) } }\hspace{-0.3cm}\frac{h^k|h|^{o-|k|}}{\bbfk !} \hspace{-0.3cm} \sum_{c\in {\sf Cut}(\beta(\bbfk,o)) \backslash {\sf Cut}(\beta) } \partial^{\bbfk}_\star{\sf P}_<\big({\sf f}_1,\dots,{\sf f}_c\big)
    \\
    &\hspace{1cm}\times\widetilde{\sf P}^{\beta_{>c},\beta(\bbfk)_{>c}}_{<}\Big( \partial^{k_{c+1}}{\sf f}_{c+1},\dots,\partial^{k_{j-1}}{\sf f}_{j-1}, R^{o-|k|}_h{\sf f}_j, \, {\sf f}_{j+1}(\cdot+h), \, \dots \Big)\{i\}
    \\
    &= \sum_{c\in {\sf I}(\beta)} \sum_{\substack{|k|<\sum_{i=1}^c\beta_i \\ \bbfk\in \mcP_{c-1}(k) } } \, \frac{h^k}{\bbfk!} \, \partial^{ \bbfk }_\star {\sf P}_<\big({\sf f}_1,\dots,{\sf f}_c\big)
\end{align*}
\begin{align*}
    &\hspace{1cm}\times\sum_{j=c+1}^n\sum_{\substack{ |p|<o \\  \bbfp \in \mcP_{j-c-1}(p)    }} \, \frac{h^p|h|^{o-|k|-|p|}}{\bbfp !} \, \widetilde{\sf P}_<^{\beta_{>c}, \beta(\bbfk,o)_{>c} }\Big( \partial^{p_{1}}{\sf f}_{c+1},\dots,\partial^{p_{j-c-1}}{\sf f}_{j-1},
    \\
    &\hspace{8.5cm}R^{o-|p|}_h{\sf f}_j, \, {\sf f}_{j+1}(\cdot+h), \, \dots \Big)\{i\} 
    \\
    &= \sum_{c\in {\sf I}(\beta)} \sum_{\substack{|k|<\sum_{i=1}^c\beta_i}} \frac{h^k}{ k!} \, \partial^k_\star\widetilde{\sf P}^{\beta_{\leq c}}_<({\sf f}_1,\dots,{\sf f}_c) \, \big(\triangle_{h,o-|k|} \widetilde{\sf P}_<^{\beta_{>c}}\big) ({\sf f}_{c+1}, \dots, {\sf f}_n)\{i\}.
\end{align*}
The result follows from this identity.
\end{Dem}

\ssk 

\begin{lem} \label{lem_grostheta2}
For $0\leq c\leq n-1$, for $o_2 > o_1> \sum_{j=c+1}^n\alpha_j -\delta_0$, we have
$$
\big(\triangle_{h, o_2} \widetilde{\sf P}_<^{\beta_{>c}}\big) ({\sf f}_{c+1},\dots,{\sf f}_n)\{i\}   - \big(\triangle_{h,o_1} \widetilde{\sf P}_<^{\beta_{>c}}\big)({\sf f}_{c+1},\dots,{\sf f}_n)\{i\} = \sum_{o_1<|k|<o_2} \, \frac{h^k}{k!} \, \partial^k_\star \widetilde{\sf P}^{\beta_{>c}}_<({\sf f}_{c+1},\dots,{\sf f}_n)\{i\}.
$$
\end{lem}

\ssk

\begin{Dem}
The proof follows the same induction as for Lemma \ref{lem_grostheta1}. The result is true for $c=n-1$ as $\triangle_{h,o}\widetilde{\sf P}_<^{\beta}$ coincides still with the Taylor remainder $|h|^o R_h^o$. Suppose it to be true for $(n-c-1)$ functions. Proposition \ref{PropAlgDevSimpl2} gives then
\begin{align*}
&\big(\triangle_{h, o_2} \widetilde{\sf P}^{\beta_{>c}}_<\big)\big({\sf f}_{c+1},\dots,{\sf f}_n\big)\{i\}-\big(\triangle_{h, o_1} \widetilde{\sf P}^{\beta_{>c}}_<\big)\big({\sf f}_{c+1}, \dots, {\sf f}_n\big)\{i\} 
\\
&= T_h^{o_2}\widetilde{\sf P}^{\beta_{>c}}_<\big({\sf f}_{c+1},\dots,{\sf f}_n\big)_i - T_h^{o_1}\widetilde{\sf P}^{\beta_{>c}}_<\big({\sf f}_{c+1}, \dots, {\sf f}_n\big)_i
\\
&\quad - \sum_{j\in {\sf I}(\beta_{>c})} \sum_{|p|<\sum_{a=c+1}^j\beta_a} \partial^p_\star {\sf P}_<({\sf f}_{c+1},\dots,{\sf f}_j) \, \frac{h^p}{p!} \,
\\
&\hspace{3.4cm}\times\Big\{ \big(\triangle_{h, o_2-\vert p\vert} \widetilde{\sf P}_<^{\beta_{>j}}\big)\big({\sf f}_{j+1},\dots,{\sf f}_n\big)\{i\} - \big(\triangle_{h, o_1-\vert p\vert} \widetilde{\sf P}^{\beta_{>j}}_<\big)\big({\sf f}_{j+1},\dots,{\sf f}_n\big)\{i\} \Big\}.
    \end{align*}
From the induction hypothesis this quantity is equal to 
\begin{align*}
&\sum_{o_1<|k|<o_2} \frac{h^k}{k!} \, \partial^k \widetilde{\sf P}_<^{\beta_{>c}} ({\sf f}_{c+1},\dots,{\sf f}_n)\{i\} - \sum_{j\in {\sf I}(\beta_{>c})} \sum_{|p|<\sum_{a=c+1}^j\beta_a} \, \frac{h^p}{p!} \, \partial^p_\star{\sf P}_<({\sf f}_{c+1}, \dots, {\sf f}_j) 
\\ 
&\hspace{7cm}\times \sum_{o_1<|\ell|+|p|<o_2} \frac{h^\ell}{\ell!} \,\partial^\ell_\star \widetilde{\sf P}^{\beta_{>c}}_<({\sf f}_{j+1}, \dots, {\sf f}_n)\{i\}   \\
&=\sum_{o_1<|k|<o_2} \, \frac{h^k}{k!} \, \sum_{\bbfk\in\mcP_{n-c}(k)} \binom{k}{\bbfk} \Lambda_{\bbfk,i},
\end{align*}
where
\begin{align*}
\Lambda_{\bbfk,i} &\defeq \widetilde{\sf P}_<^{\beta_{>c}} \big(\partial^{k_1}{\sf f}_{c+1},\dots,\partial^{k_{n-c}}{\sf f}_n\big)\{i\} -  
\\
&\hspace{-0.5cm}\sum_{j\in {\sf Cut}(\beta_{>c} - \vert\bbfk\vert) \backslash {\sf Cut}(\beta_{>c})} \hspace{-0.5cm} \widetilde{\sf P}_<^{\beta_{\llbracket c+1,j\rrbracket}-\vert\bbfk\vert} \big(\partial^{k_1}{\sf f}_{c+1},\dots,\partial^{k_{n-j+1}}{\sf f}_n\big) \, \widetilde{\sf P}_<^{\beta_{>j},\beta_{>j} - \vert\bbfk\vert} \big(\partial^{k_{n-j}}{\sf f}_{j+1},\dots,\partial^{k_{n-c}}{\sf f}_n\big)\{i\}.
\end{align*}
 Lemma \ref{rectilde2} gives 
 $$
 \Lambda_{\bbfk,i} = \widetilde{\sf P}_<^{\beta_{>c}, \beta_{>c} -\vert\bbfk\vert} \big(\partial^{k_1}{\sf f}_{c+1},\dots,\partial^{k_{n-c}}{\sf f}_n\big)\{i\},
 $$
and the result follows.
\end{Dem}

\ssk

For $0\leq a\leq n-1$ we define
$$
\big(\triangle_{yx}\widetilde{\sf P}_<^{\beta}\big)({\sf f}_{a+1},\dots,{\sf f}_n) \defeq \big(\triangle_{y-x,\sum_{j=a+1}^n \beta_j} \widetilde {\sf P}_<^{\beta} \big) ({\sf f}_{a+1},\dots,{\sf f}_n)(x);
$$
From the same arguments of Section \ref{SectionSimplifiedIPP}, for $o$ in a neighborhood of $\sum_{j=1}^n\beta_j$ one has the estimate
$$
\big| 
\big(\triangle_{yx}\widetilde{\sf P}_<^{\beta}\big)({\sf f}_{a+1},\dots,{\sf f}_n)\{i\}
  \big| \lesssim \prod_{j=1}^n \norme{{\sf f}_j}_{\beta_j} |y-x|^o 2^{-i(\sum_{j=1}^n\beta_j-o)}
$$
where 
$$
\big(\triangle_{yx}\widetilde{\sf P}_<^{\beta}\big)({\sf f}_{a+1},\dots,{\sf f}_n)\{i\}\defeq \big(\triangle_{y-x,\sum_{j=a+1}^n \beta_j} \widetilde {\sf P}_<^{\beta} \big) ({\sf f}_{a+1},\dots,{\sf f}_n)\{i\}(x).
$$
Then Lemma \ref{lem_seuilopti} gives the estimate
$$
\big|\big(\triangle_{yx}\widetilde{\sf P}_<^{\beta}\big)({\sf f}_{a+1},\dots,{\sf f}_n)\big|\lesssim \bigg\{\prod_{j=a+1}^n \norme{{\sf f}_j}_{\beta_j}\bigg\} \, |y-x|^{\sum_{j=a+1}^n \beta_j}.
$$

\ssk

\begin{prop} \label{PropExpansionTildePAlphaMinus}
Pick $\beta=(\beta_1,\dots,\beta_n)\in\bbR^n$ with $\sum_{j=1}^n \beta_j > 0$ and ${\sf f}_j\in {\sf C}^{\beta_j}_\circ$ for $1\leq j\leq n$. Then we have the local expansion
    \begin{equation} \label{EqRemainderTildePMinusAlpha} \begin{split}
        \widetilde{\sf P}^{\beta}_<({\sf f}_1,\dots,{\sf f}_n)(y) &= \sum_{|k|<\sum_{j=1}^n\beta_j} \partial_\star^k\widetilde{\sf P}^\beta_<\big({\sf f}_1,\dots,{\sf f}_n\big)(x) \, \frac{(y-x)^k}{k!}
        \\
        &\quad+ \sum_{c=1}^{n-1}\sum_{|k|<\sum_{j=1}^c\beta_j} \partial_\star^k\widetilde{\sf P}^\beta_<\big({\sf f}_1,\dots,{\sf f}_c\big)(x) \, \frac{(y-x)^k}{k!} \, \triangle_{yx}\widetilde{\sf P}_<^{\beta}({\sf f}_{c+1},\dots,{\sf f}_n)
        \\
        &\quad+ \triangle_{yx}\widetilde{\sf P}_<^{\beta}({\sf f}_{1},\dots,{\sf f}_n).
\end{split} \end{equation}
where 
$$
\Big| \big(\triangle_{yx}\widetilde{\sf P}^\beta_<\big)({\sf f}_{c+1},\dots, {\sf f}_n)\Big| \lesssim \bigg\{\prod_{j=c+1}^n \norme{{\sf f}_j}_{\beta_j}\bigg\} \, |y-x|^{\sum_{j=c+1}^n\beta_j}.
$$
\end{prop}

\ssk

\begin{Dem}
We proceed as in the proof of Proposition \ref{PropLocalDevptPMinus}. Proposition \ref{PropAlgDevSimpl2} and Proposition \ref{lem_grostheta2} give    
    \begin{align*}
        (\triangle_{yx} \widetilde{\sf P}^{\beta}_<)({\sf f}_1,\dots, {\sf f}_n)
     &=\widetilde{\sf P}_<^{\beta}\big({\sf f}_1,\dots,{\sf f}_n\big)(\cdot+h) - \sum_{|k|<o} \frac{h^k}{k!} \, \partial^k\widetilde{\sf P}^{\beta}_<\big({\sf f}_1,\dots,{\sf f}_n\big)
     \\
     &\quad-\sum_{c\in {\sf I}(\beta)} \sum_{|p|<\sum_{j=1}^c\beta_i}\partial^p_\star \widetilde{\sf P}^{\beta}_<({\sf f}_1,\dots,{\sf f}_c) \, \frac{h^p}{p!}  \,\big(\triangle_{yx} \widetilde{\sf P}^{\beta}_<\big)({\sf f}_{c+1},\dots,{\sf f}_n)
     \\
     &\quad- \sum_{c\in {\sf I}(\beta)} \sum_{\substack{|p|<\sum_{j=1}^c\beta_j \\ |\ell|>\sum_{j=c+1}^n \beta_j
 }} \frac{h^p}{p!} \, \frac{h^\ell}{\ell!} \, \partial^p_\star\widetilde{\sf P}^\beta_<({\sf f}_1,\dots,{\sf f}_c) \, \partial^{\ell}_\star \widetilde{\sf P}^{\beta}_<({\sf f}_{c+1},\dots,{\sf f}_n).
    \end{align*}
From Lemma \ref{rectilde2} we have for any $|k|<\sum_{j=1}^n\beta_j$ that $\partial_\star^k{\sf P}_<\big({\sf f}_1,\dots,{\sf f}_n\big)$ is equal to
\begin{align*}
\partial^k\widetilde{\sf P}^{\beta}_<\big({\sf f}_1,\dots,{\sf f}_n\big) -\sum_{c\in {\sf I}(\beta)}\sum_{\substack{|p|<\sum_{j=1}^c\beta_j  \\ |k-p|> \sum_{j=c+1}^n \beta_j }} \binom{k}{p} \, \partial^p_\star\widetilde{\sf P}^{\beta}_<({\sf f}_1,\dots,{\sf f}_c) \, \partial^{k-p}_\star \widetilde{\sf P}^{\beta}_<({\sf f}_{c+1},\dots,{\sf f}_n).
\end{align*}
This identity concludes the proof.
\end{Dem}

\section{The regularity structure of iterated paraproducts}
\label{RSIPP}

We fix $\alpha\in\bbR^{n}$ in this section. We introduced in Section \ref{IntroExpansionPPSystems} the spaces $T$ and $T^+$ of symbols of the regularity structure that we will associate to some iterated paraproducts. The vector space $T$ is spanned by 
$$
\mcB \defeq \Big\{  \llbracket a,b\rrbracket_{\bbfl} \, X^p \Big\}_{1\leq a<b\leq n ,\, \bbfl\in\mcP_{b-a}(\ell), \,\ell\in\bbN^{d_0},\, p\in \bbN^{d_0}} \cup \big\{X^p\big\}_{p\in\bbN^{d_0}}
$$
and the algebra is generated by
$$
\mcB^+ \defeq \Big\{ \llbracket a,b\rrbracket_{\bbfl}^\bbfk \Big\}_{\textrm{\textbf{\textsf{condition}}}(a,b,\bbfk,\bbfl)} \cup \big\{ X^{\epsilon_i} \big\}_{1\leq i\leq d},
$$
where one says that $(a,b,\bbfk,\bbfl)$ satisfies $\textrm{\textbf{\textsf{condition}}}(a,b,\bbfk,\bbfl)$ if $1\leq a<b\leq n, \bbfk=(k_a,\dots,k_b)\in\mcP_{b-a+1}(k)$ for some $k\in\bbN^{d_0}$, and $\bbfl\in\mcP_{b-a}(\ell)$ for some $\ell\in\bbN^{d_0}$, and we have
$$
\max(\vert k\vert, \vert\ell\vert) < \sum_{1\leq j\leq n} \vert \alpha_j\vert
$$
and
\begin{equation*}
\big\vert \llbracket a,b\rrbracket_{\bbfl}^\bbfk \big\vert_\alpha > 0.
\end{equation*}
We introduce in this section some splitting maps $\Delta : T\rightarrow T\otimes T^+$ and $\Delta^+ : T^+\rightarrow T\otimes T^+$ and prove in Proposition \ref{PropRS} that $((T,\Delta),(T^+,\Delta^+))$ is indeed a concrete regularity structure. We refer the reader to Appendix \ref{SectionBasicsRS} for some basics on the subject. 

\ssk

We use below the notation 
$$
M\big((\sigma_1\otimes \sigma_2), (X^{m_1}\otimes X^{m_2})\big) \defeq (\sigma_1 X^{m_1})\otimes (\sigma_2X^{m_2}).
$$
For $\tau = \llbracket a,b\rrbracket_{\bbfl} \in \mcB$, with $\bbfl=(\ell_a,\dots,\ell_b)$, we set 
$$
\oplus(\tau) = \Big\{ \llbracket a,c\rrbracket_{\bbfl_{<c}}^\bbfp \in \mcB^+ \,;\, a\leq c\leq b, \ell_c= 0, \, \bbfp\in\mcP_{c-a+1}(p), \, p\in\bbN^{d_0} \Big\} \cup \{ {\bf 1}^+ \}.
$$
For $\tau = \llbracket a,b\rrbracket_{\bbfl} \in \mcB$ and $\sigma = \llbracket a,c\rrbracket_{\bbfl_{<c}} X^\bbfp \in \oplus(\tau)$ we define $(\tau\backslash {\bf 1}^+)=\tau$ and if $c\leq b-1$ we set
$$
(\tau\backslash\sigma) \defeq \sum_{p=p_1+p_2} \sum_{{\bbfp}_1 \in \mcP_{b-c}(p_1)} \frac{ p ! }{ \bbfp_1! \, p_2! \, \bbfp! } \, \llbracket c+1,b\rrbracket_{\bbfl_{>c-a+1} + {\bbfp}_1} X^{p_2},
$$
and for $c=b$ we set $(\tau\backslash\sigma) \defeq \frac{1}{\bbfp !} \, X^p$. For $p\in\bbN^{d_0}$ we set 
$$
\Delta(X^p) = \Delta^+(X^p) \defeq \sum_{p_1+p_2=p}\binom{p}{p_1} X^{p_1}\otimes X^{p_2}.
$$
We define the map $\Delta$ on $T$ by setting
$$
\Delta\big(\llbracket a,b\rrbracket_{\bbfl} X^p\big) = M\big( \Delta(\llbracket a,b\rrbracket_{\bbfl})\,,\, \Delta(X^p)\big)
$$
and for $\tau = \llbracket a,b\rrbracket_{\bbfl} \in \mcB$
$$
\Delta(\tau) = \sum_{\sigma\in \oplus(\tau)} (\tau\backslash\sigma) \otimes \sigma.
$$
For $\mu = \llbracket a,b\rrbracket_{\bbfl}^{\bbfk} \in\mcB^+$ we set 
$$
\oplus(\mu) = \Big\{ \llbracket a,c\rrbracket_{\bbfl_{<c}}^{\bbfk_{\leq c-a+1} + \bbfp} \in\mcB^+ \,;\, a\leq c\leq b, \ell_c= 0, \, \bbfp\in\mcP_{c-a+1}(p), \, p\in\bbN^{d_0} \Big\} \cup \{ {\bf 1}^+ \}.
$$
For $\mu = \llbracket a,b\rrbracket_{\bbfl}^{\bbfk} \in \mcB^+$ and $\nu = \llbracket a,c\rrbracket_{\bbfl_{<c}}^{\bbfk_{\leq c-a+1} + \bbfp} \in \oplus(\mu)$ with $\bbfp\in\mcP_{c-a+1}(p)$, we define for $c\leq b-1$ 
\begin{equation} \label{EqMuOverNu}
(\mu\backslash\nu) \defeq \sum_{p=p_1+p_2} \sum_{{\bbfp}_1 \in \mcP_{b-c}(p_1)} \frac{ p ! }{ \bbfp_1! \, p_2! \, \bbfp! } \, 
\llbracket c+1 , b\rrbracket_{\bbfl_{>c-a+1} + {\bbfp}_1}^{\bbfk_{>c-a+1}} X^{p_2},
\end{equation}
and for $c=b$ set $(\tau\backslash\sigma) \defeq \frac{1}{\bbfp !} \, X^{p}$. All the terms in this sum have the same homogeneity $\vert\llbracket c+1,b\rrbracket_{\bbfl_{>c}}^{\bbfk_{>c-a+1}}\vert_\alpha + \vert p\vert$. We define the map $\Delta^+$ on $T^+$ by setting
$$
\Delta^+\big(\llbracket a,b\rrbracket_{\bbfl}^{\bbfk} X^p\big) = M\big( \Delta^+(\llbracket a,b\rrbracket_{\bbfl}^{\bbfk}) \,,\, \Delta^+(X^p)\big)
$$
and for $\mu = \llbracket a,b\rrbracket_{\bbfl}^{\bbfk} \in \mcB^+$
$$
\Delta^+(\mu) = \sum_{\substack{\nu\in \oplus(\mu)  \\  \vert(\mu\backslash\nu)\vert_\alpha > 0}} (\mu\backslash\nu) \otimes \nu.
$$
With the notation of \eqref{EqMuOverNu}, the condition $\vert(\mu\backslash\nu)\vert_\alpha > 0$ means that we only consider here those $\nu\in\oplus(\mu)$ such that $\vert\llbracket c+1,b\rrbracket_{\bbfl_{>c}}^{\bbfk_{> c-a+1}}\vert_\alpha >0$.

\ssk

\begin{prop} \label{PropRS}
The space $\big((T,\Delta),(T^+,\Delta^+)\big)$ is a concrete regularity structure.
\end{prop}

\ssk

In the proof of this proposition we use the following generalisation of the Vandermonde identity, which states that for any integer $i\geq 1$, for any $p,q,r$ in $\bbN^{d_0}$ such that $p+q=r$, and any $\bbfr \in\mcP_i(r)$, one has
\begin{equation} \label{eq_vandermonde}
\sum_{\substack{\bbfp\in\mcP_i(p), \bbfq \in\mcP_i(q) \\ \bbfp+\bbfq=\bbfr}} \frac{\bbfr!}{\bbfp ! \, \bbfq !} = \frac{r!}{p! \, q!}.
\end{equation}

\ssk

\begin{Dem}
We prove here that we have the comodule identity
$$
(\Delta\otimes \textrm{Id})\Delta = (\textrm{Id} \otimes \Delta^+)\Delta.
$$
The proof of the coassociativity identity
$$
(\Delta^+\otimes \textrm{Id})\Delta^+ = (\textrm{Id} \otimes \Delta^+)\Delta^+
$$
is almost identical and left to the reader. We also let the reader check the other conditions involved in the definition of a concrete regularity structure spelled out in Definition \ref{DefnRS} in Appendix \ref{SectionBasicsRS}.

It suffices to prove the comodule identity for $\tau = \llbracket a,b\rrbracket_{\bbfl} \in \mcB$ with $\bbfl=(\ell_a,\dots,\ell_b)$. To lighten the computations we use the convention $\llbracket c+1,c\rrbracket=\mathbf{1}^{(+)}$ for any $a\leq c\leq b$. We have
$$
\Delta\big(\llbracket a,b\rrbracket_{\bbfl}\big) =  \sum_{c, \bbfk_1, k_2,\bbfk'} \frac{k!}{\bbfk'! \, \bbfk_1! \, k_2! } \, \Big( \llbracket c+1,b\rrbracket_{\bbfl_{>c}+\bbfk_1} X^{k_2} \Big) \otimes\big( \llbracket a,c\rrbracket_{\bbfl_{<c}}^{\bbfk'} \big)
$$
where the sum runs over the $a\leq c \leq b$ such that $\ell_c = 0$ and the multi-indices $k=k_1+k_2$ such that $\big\vert \llbracket a,c\rrbracket_{\bbfl_{<c}}^{\bbfk} \big|_\alpha > 0$, over $\bbfk_1\in\mcP_{b-c}(k_1)$ and $\bbfk'\in\mcP_{c-a+1}(k)$. Then $(\Delta\otimes\textrm{Id})\Delta(\tau)$ is equal to

\begin{align*}
\sum_{\substack{c,\bbfk_1, k_2, \bbfk' \\ d,\bbfp_1, p_2, \bbfp'}} \frac{k! \,  p!}{\bbfk'! \, \bbfk_1! \, k_{21}! \, k_{22}! \, \bbfp'! \, \bbfp_1! \, p_2!} \, &\llbracket d+1,b\rrbracket_{\bbfl_{>d} + (\bbfk_1)_{>b-d} + \bbfp_1} X^{k_{21}+p_2}   \\
&\otimes
\llbracket c+1,d\rrbracket_{\bbfl_{\llbracket c+1,d-1\rrbracket}+ (\bbfk_1)_{\leq b-d}}^{\bbfp'} X^{k_{22}} \otimes \llbracket a,c\rrbracket_{\bbfl_{<c}}^{\bbfk'}.
\end{align*}
where the sum runs over $1<c\leq d\leq b$ such that $\ell_c$ and $\ell_d + (k_1)_{d-c}$ are null, and the multi-indices $k=k_1+k_2, \, k_2=k_{21}+k_{22},\, p=p_1+p_2$ such that 
$$
\big\vert \llbracket a,c\rrbracket_{\bbfl_{<c}}^{\bbfk} \big|_\alpha>0
$$
and 
$$
\big| \llbracket c+1,d\rrbracket_{\bbfl_{ \llbracket c+1,d-1 \rrbracket } + (\bbfk_1)_{\leq d-c}}^{\bbfp'} \big|_\alpha > 0
$$ 
and $\bbfp_1\in\mcP_{b-d}(p_1), \,\bbfp'\in\mcP_{d-c}(p)$. On the other hand $(\textrm{Id}\otimes \Delta^+)\Delta(\tau)$ is equal to
\begin{align*}
\sum_{\substack{c,\bbfk_1, k_2, \bbfk' \\ d,\bbfp_1, p_2, \bbfp'}} \frac{k! \, p! }{\bbfk'! \, \bbfk_1! \, k_2! \, \bbfp'! \, \bbfp_1! \, p_2!} \, \llbracket c+1,b\rrbracket_{\bbfl_{>c} + \bbfk_1} X^{k_2} 
&\otimes
\llbracket e+1,c\rrbracket_{\bbfl_{\llbracket e+1,c-1\rrbracket} + \bbfp_1}^{\bbfk_{> e-a+1}} X^{p_2}  \otimes 
\llbracket a,e\rrbracket_{\bbfl_{<e}}^{\bbfk_{\leq e-a+1} + \bbfp'},
\end{align*}
where the sum runs over $a<e\leq c \leq b$ such that $\ell_c=\ell_e = 0$ and multi-indices $p=p_1+p_2$ such that $\bbfp'\in\mcP_{e-a+1}(p)$ and
$$
\big| \llbracket a,e\rrbracket_{\bbfl_{<e}}^{\bbfk_{\leq e-a+1} + \bbfp'} \big|_\alpha >0
$$
and 
$$
\big| \llbracket e+1,c\rrbracket_{\bbfl_{\llbracket e+1,c-1\rrbracket} + \bbfp_1}^{\bbfk_{> t}} \big|_\alpha >0.
$$ 

Both sums take the form
$$
\sum_{\substack{c,\bbfk_1, k_2, \bbfq \\ d,\bbfp_1, p_2, \bbfq'}} C_{\bbfp_1,p_2,\bbfq'}^{\bbfk_1,k_2,\bbfq} \llbracket d+1,b\rrbracket_{\bbfl_{>d}+\bbfk_1} X^{k_2} 
\otimes 
\llbracket c+1;d\rrbracket_{\bbfl_{\llbracket c+1,d-1\rrbracket} + \bbfp_1}^{\bbfq} X^{p_2} \otimes 
\llbracket a,c\rrbracket_{\bbfl_{<c}}^{\bbfq'},
$$
where the sum runs over $a\leq c\leq d\leq b$ such that $\ell_c, \ell_d \neq 0$, over multi-indices and tuples of multi-indices $\bbfk_1, k_2, \bbfq,\bbfp_1, p_2, \bbfq'$ such that the first two terms in each tensor products are in $T^+$ and $q+q'=k_1+k_2+p_1+p_2$ and $q\leq k_1+k_2$.

We check that the constants $C_{\bbfp_1,p_2,q'}^{\bbfk_1,k_2,q}$ coincide in both expressions using the Vandermonde identity \eqref{eq_vandermonde}. Both are equal to
$$
C_{\bbfp_1,p_2,q'}^{\bbfk_1,k_2,q} = \frac{1}{\bbfk_1! \, k_2! \, \bbfp_1! \, p_2! \, \bbfq! \, \bbfq'!} \, \frac{k! \, q'!}{q!}.
$$
This concludes the proof of the statement.
\end{Dem}

\ssk

We note that Hoshino was the first to investigate in \cite{ComutBesov} the algebraic structure behind the iterated paraproducts, in a restricted setting compared to the present general setting.

\ssk

For $\tau\in T$, one can re-index the sum defining $\Delta(\tau)$ by its different components $\tau_1$ on the canonical basis $\mcB$ of the $T$ factor in $T\otimes T^+$ and write

$$
\Delta(\tau) \eqdef \sum_{\tau_1 \leq \tau} \tau_1 \otimes (\tau/\tau_1).
$$
(This identity defines $(\tau/\tau_1)$.) Below we write $\tau_1<\tau$ to mean that $\tau_1$ appears in this decomposition and $\tau_1 \neq \tau$. Similarly we can rewrite 
$$
\Delta^+(\mu) \eqdef \sum_{\mu_1 \leq^+ \mu} \mu_1 \otimes (\mu/^+\mu_1).
$$
(This identity defines $(\mu/^+\mu_1)$.) Below we write $\mu_1<^+\mu$ to mean that $\mu_1$ appears in this decomposition and $\mu_1 \neq \mu$.

\bigskip

\section{Local expansion properties of iterated paraproducts}
\label{SectionProofThmsMain}

We prove Theorem \ref{ThmMain1} in this section. The core of the proof rests on a representation formula
$$
{\sf P}_{\bbfl}(f_a,\dots,f_b) = \sum_{c\geq 0} \sum_{\tau_1,\dots,\tau_c} {\sf P}_<\Big([\tau/\tau_1]^{\bm{f}} , \dots , [\tau_c]^{\bm{f}} \Big)
$$
of the ${\sf P}_{\bbfl}$ operators in terms of the ${\sf P}_<$ operators and some functions $[\sigma]^{\bm{f}}$ that we build from $\bm{f}=(f_a,\dots,f_b)$. One can then infer the local expansion properties of ${\sf P}_{\bbfl}(f_a,\dots,f_b)$ from the local expansion properties of the ${\sf P}_<$ operators obtained in Section \ref{SectionSimplifiedIPP} and Section \ref{SectionLocalExpansionDStarK}. We describe in Section \ref{SectionBuildingBlocks} the generic construction of the bracket maps $[\sigma]$ if the initial data is a pair of maps $\sf (\Pi,g)$ of a particular type. They will be specified in Section \ref{SectionRepresentationFormula} in terms of a fixed tuple $\bm{f}=(f_1,\dots,f_n)$ of distributions. The actual proof of Theorem \ref{ThmMain1} occupies all of Section \ref{SectionProofTheorems}. The inductive mechanics of this proof is detailed at the begining of this section.


\subsection{Building blocks for a representation of $\sf P$ in terms of ${\sf P}_<$.}
\label{SectionBuildingBlocks}

\hspace{0.1cm} Recall from Appendix \ref{SectionBasicsRS} the basic notions and notations on regularity structures. For any integer $n_0$ define $\mathcal{Z}_{n_0}$ as the set of $C^{n_0}$ functions $\phi$ supported in the unit ball of ${\bbR^{d_0}} $ and such that $\norme{\phi}_{C^{n_0}}\leq 1$. Given a pair of maps $({\sf \Pi},{\sf g})$ with $\sf \Pi$ a linear map from $T$ to $\mcD'(\bbR^{d_0})$ and  $\sf g$ a map from $\bbR^{d_0}$ into the set of characters on the algebra $T^+$, we define 
$$
{\sf \Pi}_x \defeq ({\sf \Pi}\otimes {\sf g}_x^{-1})\Delta
$$ 
for any $x\in \bbR^{d_0}$. For a real-valued function $\varphi$ on $\bbR^{d_0},\,x\in\bbR^{d_0}$ and $\epsilon>0$ we define 
$$
\varphi_x^\epsilon(y) \defeq \epsilon^{-d_0} \varphi\big(\epsilon^{-1} (y-x)\big).
$$
We define the size $\llparenthesis\sf \Pi , g\rrparenthesis$ of $\sf (\Pi,g)$ by setting first for $\tau\in T_{|\tau|}$ 
$$
\Vert\tau\Vert_{\sf (\Pi,g)} \defeq  \sup_{x\in{\bbR^{d_0}}}\sup_{\phi \in \mathcal{Z}_{n_0}}\sup_{\epsilon\in(0,1)} \epsilon^{-\vert\tau\vert} \big|\langle {\sf \Pi}_x\tau, \phi_x^\epsilon \rangle\big|,
$$
and for $\nu\in T_{|\nu|}^+$
$$
\Vert\nu\Vert_{\sf (\Pi,g)} \defeq  \sup_{x,y\in{\bbR^{d_0}}}\frac{{\sf g}_{yx}(\nu)}{|y-x|^{|\nu|}},
$$
and recursively for $\tau\in \mcB^{(+)}$
$$
\Vert\tau\Vert_{\sf (\Pi,g)}^* \defeq \max \Big( \Vert\tau\Vert_{\sf (\Pi,g)},\, \max_{\sigma<^{(+)}\tau} \Vert\tau/\sigma\Vert_{\sf (\Pi,g)} \Vert\sigma\Vert_{\sf (\Pi,g)}^*\Big),
$$
and then by defining
$$
\llparenthesis\sf \Pi,g\rrparenthesis \defeq \max\Big(\Vert\tau\Vert_{\sf (\Pi,g)}^*  \,,\, \Vert\mu\Vert_{\sf (\Pi,g)}  \Big),
$$
for a maximum over $\tau\in\mcB $  and $\mu\in\mcB^+$.We have
$$
{\sf \Pi}(\tau) = \sum_{\tau_1 \leq \tau}  {\sf g}_x(\tau/\tau_1) \, {\sf \Pi}_x(\tau_1),
$$
that is
$$
{\sf \Pi}_x(\tau) = {\sf \Pi}(\tau)- \sum_{\tau_1 < \tau} {\sf g}_x(\tau/\tau_1) \, {\sf \Pi}_x(\tau_1).
$$
Iterating we obtain the formula
\begin{equation} \label{eq_devPi}
    {\sf \Pi}_x(\tau)= {\sf \Pi}(\tau) - \sum_{e\geq 1}(-1)^{e-1} \sum_{\tau_e<\cdots<\tau_1<\tau}  {\sf g}_x(\tau/\tau_1)\cdots {\sf g}_x(\tau_{q-1}/\tau_e) \, {\sf \Pi}(\tau_e),
\end{equation}
where the sum over $e$ is finite as the sets $A$ and $A^+$ that index the homogeneities of the regularity structure $((T,\Delta),(T^+,\Delta^+))$ are locally finite and bounded from below. Likewise for $\tau/\rho\in T^+$
\begin{equation} \label{eq_devg}
\begin{split}
    {\sf g}_{yx}\big(\tau/\rho\big) =& \; {\sf g}_y\big(\tau/\rho\big) - {\sf g}_x\big(\tau/\rho\big) 
    \\
    &-\sum_{e\geq 1}(-1)^{e-1} \sum_{\rho<\tau_e<\cdots<\tau_1<\tau}  {\sf g}_x\big(\tau/\tau_1\big)\cdots {\sf g}_x\big(\tau_{e-1}/\tau_e\big) \Big( {\sf g}_y\big(\tau_e/\tau_1\big) - {\sf g}_x\big(\tau_e/\rho\big)\Big).
\end{split}
\end{equation}
For $\tau\in T$ we define an element $[\tau] = ( [\tau]_i )_{i\geq -1}$ of ${\sf C}^{-\infty}$ setting
$$
[\tau]_i \defeq \Delta_i \big({\sf \Pi}(\tau)\big)- \sum_{e\geq 1} (-1)^{e-1} \sum_{\tau_e<\cdots<\tau_1<\tau} \Delta_{<i-1}\big({\sf g}\big(\tau/\tau_1\big)\big)\cdots \Delta_{<i-1}\big({\sf g}\big(\tau_{q-1}/\tau_e\big)\big)  \Delta_i\big({\sf \Pi}\big(\tau_e\big)\big).
$$
Likewise for $\tau/\rho\in T^+$ we set
\begin{align*}
[\tau/\sigma]_i &\defeq \Delta_i \big({\sf g}(\tau/\rho)\big)
\\
&\qquad -\sum_{e\geq 1} (-1)^{e-1} \sum_{\rho<\tau_e<\cdots<\tau_1<\tau} \Delta_{<i-1}\big({\sf g}\big(\tau/\tau_1\big)\big)\cdots \Delta_{<i-1}\big({\sf g}\big(\tau_{q-1}/\tau_e\big)\big)  \Delta_i\big({\sf g}\big(\tau_e/\rho\big)\big).
\end{align*}

\ssk

\begin{prop} \label{lemregcrochet}
For any $\tau\in T_{|\tau|}$ and $\tau/\rho\in T_{|\tau/ \rho|}^+$ we have
$$
\Vert [\tau]\Vert_{|\tau|} + \Vert [\tau/\rho] \Vert_{|\tau/\rho|} \lesssim \llparenthesis\sf \Pi , g\rrparenthesis,
$$
so $[\tau]\in {\sf C}^{|\tau|}$ and $[\tau/\rho]\in {\sf C}^{|\tau/\rho|}$ if $\llparenthesis\sf \Pi , g\rrparenthesis<\infty$.
\end{prop}

\ssk

The proof of this statement uses the following result stated in Proposition 8 of Bailleul \& Hoshino's work \cite{PCandRS1}. We denote below by $K_j(x-y)$ the translation-invariant kernel of the Littlewood-Paley projector $\Delta_j$ and set 
$$
K_{<i-1} \defeq \sum_{-1\leq j<i-1} K_j.
$$

\ssk

\begin{lem} \label{regu2parap}
Let $F=(F_x)_{x\in\bbR^{d_0}}$ be a family of distributions on $\bbR^{d_0}$ indexed by $\bbR^{d_0}$. Set
$$
(Q_iF)(z) \defeq \int K_{<i-1}(z-x) F_x\big(K_i(z-\cdot)\big) dx
$$
and assume that
$$
\Vert Q_iF\Vert_\infty \leq C_F 2^{-i r_1}
$$
for some positive constant $C_F$ and $r_1\in\bbR$. Let $G$ be a function on $(\bbR^{d_0})^2$ such that we have
$$
|F(x,y)|\leq C_G |y-x|^{r_2}
$$
for all $x,y$, for some exponent $r_2>0$ and some positive constant $C_G$. Set
\begin{equation*} \begin{split}
\big(Q^+_iF\big)(z) \defeq \iint K_{<i-1}(z-x) K_{<i-1}(z-y) F(x,y) \, dxdy.
\end{split} \end{equation*}
Then $QF=(Q_iF)_{i\geq -1} \in {\sf C}^{r_1}$ and $Q^+G=(Q_i^+G)_{i\geq -1}\in {\sf C}^{r_2}$ with
$$
\norme{QF}_{r_1} + \norme{Q^+G}_{r_2} \lesssim C_F+C_G.
$$
\end{lem}

\ssk

\begin{Dem}[of Proposition \ref{lemregcrochet}]
We proceed by induction. For $\tau\in T_{|\tau|}$ and $\tau/\rho\in T_{|\tau/ \rho|}^+$ we set $F_{\tau x} = {\sf \Pi}_x\tau$ and $G_{\tau/\rho}(x,y)={\sf g}_{yx}(\tau/\rho)$, for all $x,y\in\bbR^{d_0}$. Writing 
$$
{\sf \Pi}_x\tau = {\sf \Pi}_z\tau + \sum_{\sigma<\tau} {\sf g}_{zx}(\tau/\sigma) \, {\sf \Pi}_z\sigma
$$
we see that
$$
\big(Q_iF_\tau\big)(z) = \int K_{<i-1}(z-x) ({\sf \Pi}_z\tau)\big(K_i(z-\cdot)\big) dx + \sum_{\sigma<\tau} \int K_{<i-1}(z-x) \, {\sf g}_{zx}(\tau/\sigma) \, ({\sf \Pi}_z\sigma)\big(K_i(z-\cdot)\big) dx
$$
with
$$
\big\vert ({\sf \Pi}_z\tau)\big(K_i(z-\cdot)\big) \big\vert \lesssim  2^{-i\vert\tau\vert} \Vert\tau\Vert_{\sf (\Pi,g)}
$$
uniformly in $z$, with a similar estimate with $\sigma$ in place of $\tau$, and 
$$
\int \big\vert K_{<i-1}(z-x) \, {\sf g}_{zx}(\tau/\sigma) \big\vert dx \lesssim 2^{i\vert\tau/\sigma\vert} \Vert\tau/\sigma\Vert_{\sf (\Pi,g)}.
$$
It follows that 
$$
\Vert Q_iF_\tau\Vert_\infty \lesssim 2^{-i \vert\tau\vert} \max\big\{\Vert\tau/\sigma\Vert_{\sf (\Pi,g)} \,;\, \sigma\leq \tau\big\},
$$
go we get from Lemma \ref{regu2parap} that $QF_\tau \in {\sf C}^{|\tau|}$ with $\norme{QF_\tau}_{|\tau|}\lesssim \max\big\{\Vert\tau/\sigma\Vert_{\sf (\Pi,g)} \,;\, \sigma\leq \tau\big\}$. Note that 
$$
Q_iF_\tau = \Delta_i({\sf \Pi}\tau) - \sum_{e\geq 1} (-1)^{e-1} \sum_{\sigma_e<\cdots<\sigma_1<\tau} \Delta_{i-1}\Big( {\sf g}(\tau/\sigma_1) \cdots {\sf g}(\sigma_{e-1}/\sigma_e)\Big) \, \Delta_i({\sf \Pi}\sigma_e).
$$
On the other hand one has directly from Lemma \ref{regu2parap} that $Q^+G_{\tau/\rho}\in {\sf C}^{|\tau/\rho|}$ with norm bounded above by a constant multiple of $\Vert\tau/\rho\Vert_{\sf (\Pi,g)}$. We actually have from \eqref{eq_devg} the following formula for
\begin{align*}
\big( Q^+G_{\tau/\rho} \big)_i &= \Delta_{<i-1}\big({\sf g}(\tau/\sigma)\big)  
\\
&\qquad- \sum_{e\geq 1} (-1)^{e-1}\sum_{\sigma<\sigma_e<\cdots<\sigma_1<\tau} \Big\{\Delta_{<i-1}\Big({\sf g}\big(\tau/\sigma_1\big)\cdots {\sf g}\big(\sigma_{e-1}/\sigma_e\big)\Big) \, \Delta_{<i-1} \big( {\sf g}\big(\sigma_e/\sigma\big) \big)
\\ 
&\hspace{6cm}- \Delta_{<i-1}\Big({\sf g}(\tau/\sigma_1)\cdots{\sf g}(\sigma_{e-1}/\sigma_e) \, {\sf g}(\sigma_e/\sigma)\Big)\Big\}
\end{align*}
It follows from induction that $\big((Q^+G_{\tau/\rho})_i [\rho]_i\big)_{i\geq -1}$ defines an element of ${\sf C}^{|\tau|}$ with norm bounded by a constant multiple of $\llparenthesis\sf \Pi,g\rrparenthesis$. The conclusion of Proposition \ref{lemregcrochet} will the follow after we check that
\begin{equation}\label{eq_PQtocrochet}
[\tau]_i = Q_iF_\tau +  \sum_{\sigma<\tau} \big(Q_i^+G_{\tau/\sigma}\big) [\sigma]_i.
\end{equation}
To see that one has this identity we notice that for any $\sigma<\tau$ one has
\begin{align*}
\big( Q_i^+G_{\tau/\sigma}\big) [\sigma]_i &= \underset{\sigma_{e_2}<\cdots<\sigma<\nu_{e_1}<\cdots<\tau}{\sum_{e_1,e_2\geq 0}}  (-1)^{e_1+e_2} \Delta_{<i-1}\Big({\sf g}(\tau/\nu_1)\cdots {\sf g}(\nu_{e_1-1}/\nu_{e_1}) \big)\Delta_{<i-1}\big({\sf g}(\nu_{e_1}/\sigma)\Big) 
\\
&\hspace{3.5cm}\times\Delta_{<i-1}\big({\sf g}(\sigma/\sigma_1)\big) \cdots \Delta_{<i-1}\big({\sf g}(\sigma_{e_2-1}/\sigma_{e_2}\big)\Delta_i\big({\sf \Pi}(\sigma_{e_2})\big)
\\
&\qquad-  \underset{\sigma_{e_2}<\cdots<\sigma<\nu_{e_1}<\cdots<\tau}{\sum_{e_1,e_2\geq 0}} (-1)^{e_1+e_2}\Delta_{<i-1}\big({\sf g}(\tau/\nu_1)\cdots {\sf g}(\nu_{e_1-1}/\nu_{e_1}){\sf g}(\nu_{e_1}/\sigma)\big) 
\\
&\hspace{3.5cm}\times\Delta_{<i-1}\big({\sf g}(\sigma/\sigma_1)\big) \cdots \Delta_{<i-1}\big({\sf g}(\sigma_{e_2-1}/\sigma_{e_2})\big)\Delta_i\big({\sf \Pi}(\sigma_{e_2})\big),
\end{align*}
so summing over $\sigma<\tau$ one recognizes a telescopic sum which simplifies indeed to \eqref{eq_PQtocrochet}.
\end{Dem}

\medskip

\subsection{A representation formula.}
\label{SectionRepresentationFormula}

\hspace{0.1cm} We fix here a tuple $\bm{f}=(f_1,\dots,f_n)$ of smooth functions and prove in Proposition \ref{PropoRepresentationProposition} below that one can represent any ${\sf P}_{\bbfl}(f_a,\dots,f_b)$ as a sum of ${\sf P}_<$ terms involving the bracket functions $[\tau]$ constructed in Section \ref{SectionBuildingBlocks}. We prove a similar representation formula for $\widetilde{\sf P}^{\llbracket a,b\rrbracket-\vert\bbfk\vert}_{\bbfl}(\partial^{k_a}f_a,\dots,\partial^{k_b}f_b)$. We need some preparatory results before stating and proving Proposition \ref{PropoRepresentationProposition}.

\ssk

As in \eqref{EqPiMap1}, \eqref{EqPiMap2} and \eqref{EqGMap} in Section \ref{SectionIntro} we associate to $\bm{f}$ the pair of maps 
$$
{\sf M}^{\bm{f},\alpha} = {\sf M}^{\bm{f}} = ({\sf \Pi}^{\bm{f}} , {\sf g}^{\bm{f}})
$$ 
where
$$
{\sf \Pi}^{\bm{f}}\big(\llbracket a,b\rrbracket_{\bbfl} X^p\big)(y) \defeq y^p\,{\sf P}_{\bbfl}\big(f_a,\dots,f_b\big)(y)
$$
and
$$
{\sf g}^{\bm{f}}\big(\llbracket a,b\rrbracket^{\bbfk}_{\bbfl} X^q\big)(y) \defeq y^q \, \widetilde{\sf P}^{\alpha_{\llbracket a,b\rrbracket}-\vert\bbfk\vert}_{\bbfl}\big(\partial^{k_a}f_a,\dots, \partial^{k_b}f_b\big)(y).
$$
We denote by $[\cdot]^{\bm{f}}$ the bracket maps associated to $({\sf \Pi}^{\bm{f}}, {\sf g}^{\bm{f}})$. For each $\llbracket a,b\rrbracket_{\bbfl} X^p \in \mcB$ we define an element $\big({\sf \Pi}^{\bm{f}}(\llbracket a,b\rrbracket_{\bbfl} X^p)_i\big)_{i\geq -1}$ of ${\sf C}^{0+}$ by setting

$$
{\sf \Pi}^{\bm{f}}(\llbracket a,b\rrbracket_{\bbfl} X^p)_i \defeq \Delta_i^p\big( {\sf P}_{\bbfl}(f_a,\cdots,f_b) \big).
$$
Likewise, for $\llbracket a,b\rrbracket_{\bbfl}^{\bbfk} X^q \in \mcB^+$ with $\bbfk=(k_a,\dots,k_b)$, we define an element $\big({\sf g}^{\bm{f}}(\llbracket a,b\rrbracket_{\bbfl}^{\bbfk} X^q )_i\big)_{i\geq -1}$ of ${\sf C}^{0+}$ by setting
$$
{\sf g}^{\bm{f}}(\llbracket a,b\rrbracket_{\bbfl}^{\bbfk} X^q )_i \defeq \Delta_i^q\Big(\widetilde{\sf P}_{\bbfl}^{\alpha_{\llbracket a,b\rrbracket} - \vert\bbfk\vert} \big(\partial^{k_a}f_a,\cdots,\partial^{k_b}f_b\big) \Big).
$$
For $j\geq 0$ we set
$$
{\sf g}^{\bm{f}}(\llbracket a,b\rrbracket_{\bbfl}^{\bbfk} X^q )_{<j} \defeq \sum_{i=-1}^{j-1} {\sf g}^{\bm{f}}(\llbracket a,b\rrbracket_{\bbfl}^{\bbfk} X^q )_i.
$$
The statement of Proposition \ref{PropoRepresentationProposition} below, and the next two preparatory results, require a notation that we now introduce. For $\tau = \llbracket a,b\rrbracket_{\bbfl} X^p \in \mcB$ we write 
$$
\sigma\prec\tau \textrm{ if } \sigma<\tau \textrm{ and } \sigma = \llbracket c,b\rrbracket_{\bbfl'}X^{p'} \textrm{ with } c>a.
$$ 
We also write 
$$
\sigma\leqq\tau \textrm{ if } \sigma<\tau \textrm{ but not } \sigma\prec\tau.
$$ 
For a descending sequence $\sigma_e\leqq \cdots \leqq\sigma_1 \leqq \tau$ we have $\sigma_j =  \llbracket a,b\rrbracket_{\bbfl} X^{p_j}$ with $0\leq p_e<\cdots<p_1<p$, and $\sigma_j/\sigma_{j+1} = \binom{p_j}{p_{j+1}}X^{p_j-p_{j+1}}$. For $\mu = \llbracket a,b\rrbracket_{\bbfl}^{\bbfk} X^q \in \mcB^+$ we write 
$$
\nu\prec\mu \textrm{ if } \nu<\mu \textrm{ and } \nu = \llbracket c,b\rrbracket_{\bbfl'}^{\bbfk'} X^{q'}  \textrm{ with } c>a.
$$ 

\ssk

\begin{lem}\label{lem_gtogtilde}
We have for $\tau\in T$ and $i\geq 1$
\begin{equation} \label{EqPiFTauI}
{\sf \Pi}^{\bm{f}}(\tau)_i = \Delta_i( {\sf \Pi}^{\bm{f}}(\tau)) - \sum_{e\geq 1}(-1)^e \hspace{-0.2cm} \sum_{\sigma_e\leqq \cdots \leqq\sigma_1 \leqq \tau} \hspace{-0.3cm} \Delta_{<i-1}\big({\sf g}^{\bm{f}}\big(\tau/\sigma_1\big)\big)\cdots \Delta_{<i-1}\big({\sf g}^{\bm{f}}\big(\sigma_{e-1}/\sigma_e\big)\big) \, \Delta_i\big({\sf \Pi}^{\bm{f}}(\sigma_e)\big)
\end{equation}
and for $\tau/\sigma\in T^+$ with $\sigma\prec\tau$
\begin{equation} \label{EqGFTauI} \begin{split}
& {\sf g}^{\bm{f}}(\tau/\sigma)_{<i-1} = 
\\
& \Delta_{<i-1}\big({\sf g}^{\bm{f}}(\tau/\sigma)\big) - \sum_{e\geq 1}(-1)^e \hspace{-0.5cm} \underset{\sigma \prec \sigma_e}{\sum_{\sigma_e\leqq \cdots \leqq \sigma_1 \leqq \tau}} \hspace{-0.4cm} \Delta_{<i-1}\big({\sf g}^{\bm{f}}\big(\tau/\sigma_1\big)\big)\cdots \Delta_{<i-1}\big({\sf g}^{\bm{f}}\big(\sigma_{e-1}/\sigma_e\big)\big) \, \Delta_{<i-1}\big({\sf g}^{\bm{f}}(\sigma_e/\sigma)\big).   
\end{split} \end{equation}
\end{lem}

\ssk

\begin{Dem}
We consider first the identity \eqref{EqPiFTauI}. Denote by $(\star)_i(\cdot)$ its right hand side. It suffices to treat the case of $\tau =  \llbracket 1,n\rrbracket_{\bbfl} X^p$. One has $\tau/\sigma_1 = \binom{p}{p_1}X^{p-p_1} $ and $\sigma_j/\sigma_{j+1} = \binom{p_j}{p_{j+1}}X^{p_j-p_{j+1}}$ for $1\leq j\leq e-1$; moreover for $k\in\bbN^{d_0}$ we have $\Delta_{<i-1}\big( {\sf g}^{\bm{f}}(X^k)  \big)(x) = x^k$ for all $i\geq 1$. It follows that $(\star)_i(x)$ is equal to

\begin{align*}
&\Delta_i( {\sf \Pi}^{\bm{f}}(\llbracket 1,n\rrbracket_{\bbfl} X^p))(x) - \sum_{e\geq 1} (-1)^e \sum_{0\leq p_e<\cdots<p_1<p} \prod_{j=1}^{e-1} \binom{p_j}{p_{j+1}}x^{p_j-p_{j+1}} 
\Delta_i\Big( {\sf \Pi}^{\bm{f}}\big(\llbracket 1,n\rrbracket_{\bbfl}  X^{p_q} \big) \Big)(x)
\\
&= \Delta_i\Big( {\sf \Pi}^{\bm{f}}\big(\llbracket 1,n\rrbracket_{\bbfl} X^p \big) \Big)(x) 
- \sum_{r<p} C_{pr} \, x^{p-r} \Delta_i\Big( {\sf \Pi}^{\bm{f}}\big(\llbracket 1,n\rrbracket_{\bbfl}  X^{r} \big)  \Big)(x)
\end{align*}
where
$$
C_{pr}\defeq  \sum_{e\geq 1} (-1)^e \sum_{ r<p_{e-1}<\cdots<p_1<p } \prod_{j=0}^{e-1} \binom{p_j}{p_{j+1}}.
$$
We note that the constants $C_{pr}$ satisfy the inductive relation

$$
C_{pr} = -\binom{p}{r}  + \sum_{r<s<p}\binom{p}{s} C_{sr}, 
$$
so
$$
C_{pr} = (-1)^{p-r+1} \binom{p}{r}.
$$
One then has
\begin{align*}
(\star)_i(x) &= \Delta_i\Big( {\sf \Pi}^{\bm{f}}\big(\llbracket 1,n\rrbracket_{\bbfl} X^p \big) \Big)(x) 
+ 
\sum_{r<p} (-1)^{p-r} \binom{p}{r} x^{p-r} \Delta_i\Big( {\sf \Pi}^{\bm{f}}\big(\llbracket 1,n\rrbracket_{\bbfl}  X^{r} \big) \Big)(x)
\\
&=\int_{{\bbR^{d_0}}} K_i(y-x) \sum_{r=0}^p (-1)^{p-r}\binom{p}{r}  x^{p-r} y^r \, {\sf P}_{\bbfl}(f_1,\cdots,f_n)(y) \, dy
\\
&= {\sf \Pi}^{\bm{f}}(\tau)_i(x).
\end{align*}

One uses a similar reasoning to prove identity \eqref{EqGFTauI}. It suffices to treat the case $\tau = \llbracket 1;N \rrbracket_{\bbfl} X^p$ and  $\sigma = \llbracket n+1,N\rrbracket_{\bbfl_{>n} + {\color{blue}\bm{s}}} X^q$. One has in that case
$$
\tau / \sigma =  \sum_{ s_1,r } \frac{1}{{{\color{blue}\bm{s}}}! (s_1 - |{{\color{blue}\bm{s}}}|)!}  \binom{p}{r} \llbracket 1;n \rrbracket_{\bbfl_{< n} }^{s_1} X^r
$$
where the sum runs over the multi-indices $s_1,r$ such that $ p=q+\norme{{\color{blue}\bm{s}}}+r-s_1$ and such that $\llbracket 1;n \rrbracket_{\bbfl_{<n} }^{s_1}\in T^+$ and $r\geq 0$. We write $D_{p,q}$ for the set of such $s_1,r$. Writing $(\star\star)_i(\cdot)$ for the right hand side of \eqref{EqGFTauI}, we have this time 
\begin{align*}
    (\star\star)_i(x) &=  \sum_{(s_1,r)\in D_{p,q}}  \frac{1}{{{\color{blue}\bm{s}}}! (s_1 - |{{\color{blue}\bm{s}}}|)!} \binom{p}{r}\Delta_{<i-1}\Big( {\sf \Pi}^{\bm{f}}\big(\llbracket 1,n\rrbracket_{\bbfl}^{s_1} X^r \big) \Big)(x) 
\\
&\quad- 
\sum_{p'<p} \, \sum_{(s_1,r')\in D_{p',q}}  C_{pp'}  x^{p-p'} \binom{p'}{r'}\frac{1}{{{\color{blue}\bm{s}}}! (s_1 - |{{\color{blue}\bm{s}}}|)!} \Delta_{<i-1}\Big( {\sf \Pi}^{\bm{f}}\big(\llbracket 1,n\rrbracket_{\bbfl}^{s_1}  X^{r'} \big) \Big)(x)
\\
&=\sum_{(s_1,r)\in D_{p,q}}  \frac{1}{{{\color{blue}\bm{s}}}! (s_1 - |{{\color{blue}\bm{s}}}|)!} \binom{p}{r}\Delta_{<i-1}\Big( {\sf \Pi}^{\bm{f}}\big(\llbracket 1,n\rrbracket_{\bbfl}^{s_1} X^r \big) \Big)(x) 
\\
& \quad - \sum_{(s_1,r)\in D_{p,q}} \,  \sum_{r' < r} \frac{1}{{{\color{blue}\bm{s}}}! (s_1 - |{{\color{blue}\bm{s}}}|)!} \binom{p}{r} (-1)^{r-r'} x^{r-r'}  \binom{r}{r'} \Delta_{<i-1}\Big( {\sf \Pi}^{\bm{f}}\big(\llbracket 1,n\rrbracket_{\bbfl}^{s_1}  X^{r'} \big) \Big)(x).
\end{align*}
where we used that, for any fixed $s_1$, if $(s_1,r')\in D_{p',q}$ and $(s_1,r)\in D_{p',q}$ then $p-p'=r-r'$. This gives indeed ${\sf g}^{\bm{f}}(\tau/\sigma)_{<i-1}$.
\end{Dem}

\ssk

\begin{cor}\label{lem_sumgtilde}
For any $\tau\in\mcB$ and $i\geq 1$ we have the relation
$$
[\tau]^{\bm{f}}_i = {\sf \Pi}^{\bm{f}}(\tau)_i - \sum_{e\geq 1} (-1)^{e-1} \sum_{\sigma_e\prec\cdots\prec\sigma_1\prec\tau} {\sf g}^{\bm{f}}\big(\tau/\sigma_1\big)_{<i-1} \cdots {\sf g}^{\bm{f}}\big(\sigma_e/\sigma_{e-1}\big)_{<i-1} \, {\sf \Pi}^{\bm{f}}\big(\sigma_e\big)_i.
$$
Likewise for $\tau/\sigma\in \mcB^+$ we have
\begin{align*}
[\tau/\sigma]^{\bm{f}}_i &= {\sf g}^{\bm{f}}(\tau/\sigma)_i - \sum_{e\geq 1} (-1)^{e-1} \sum_{\sigma\prec\sigma_e\prec\cdots\prec\sigma_1\prec\tau} {\sf g}^{\bm{f}}\big(\tau/\sigma_1\big)_{<i-1}\cdots {\sf g}^{\bm{f}}\big(\sigma_e/\sigma_{e-1}\big)_{<i-1} \, {\sf g}^{\bm{f}}\big(\sigma_e/\sigma\big)_i.
\end{align*}
\end{cor}

\ssk

\begin{Dem} 
Plugging the identity of Lemma \ref{lem_gtogtilde} giving ${\sf \Pi}^{\bm f}$ and ${\sf g}^{\bm f}$ into the right hand of the identity to prove, developing the products, one recovers the definition of $[\tau]_i^{\bm f}$ and $[\tau/\sigma]_i^{\bm f}$ by  noting that any descending sequence $\tau_e<\cdots<\tau_1<\tau$ takes the form
$$
\cdots\leqq\tau_{3,0}\prec \tau_{2,e_2} \leqq \cdots  \leqq  \tau_{2,0} \prec \tau_{1,e_1} \leqq \cdots \leqq \tau_{1,1} \leqq \tau.
$$
The conclusion follows.
\end{Dem}

\ssk

\begin{prop} \label{PropoRepresentationProposition}
For any $\tau = \llbracket a,b\rrbracket_{\bbfl} X^p \in T$ we have
\begin{equation} \label{eq_ecrituresompart}
\Delta_i^p\big( {\sf P}_{\bbfl}(f_a,\cdots,f_b) \big) = \sum_{e\geq 0}\hspace{3pt}\sum_{\sigma_e\prec\cdots\prec\sigma_1\prec\tau} {\sf P}_<\big([\tau/\sigma_1]^{\bm{f}}, \dots,[\sigma_e]^{\bm{f}} \big)_i
\end{equation}
and for $\sigma\leq \tau$ with $\tau/\sigma =  \llbracket c,d\rrbracket_{\bbfl'}^{\bbfk} X^q \in T^+$ we have
\begin{equation}\label{eq_ecrituresompart2}
\Delta_i^q\Big(\widetilde{\sf P}_{\bbfl'}^{\alpha_{\llbracket c,d\rrbracket} - \vert\bbfk\vert} \big(\partial^{k_c}f_c,\cdots,\partial^{k_d}f_d\big) \Big) = \sum_{e\geq 0}\hspace{3pt}\sum_{\sigma\prec\sigma_e\prec\cdots\prec\sigma_1\prec\tau} {\sf P}_<\big([\tau/\sigma_1]^{\bm{f}}, \dots, [\sigma_e/\sigma]^{\bm{f}} \big)_i
\end{equation}
\end{prop}

\ssk

\begin{Dem}
We prove \eqref{eq_ecrituresompart} and let the reader prove \eqref{eq_ecrituresompart2} as its proof is almost identical. We proceed by developing the sum 
$$
\sum_{e\geq 0} \; \sum_{\sigma_e\prec\cdots\prec\sigma_1\prec\tau} {\sf P}_<\big([\tau/\sigma_1]^{\bm{f}}, \cdots,[\sigma_e]^{\bm{f}} \big)_i
$$ 
and use the identities of Corollary \ref{lem_sumgtilde} to see that a number of cancellations give in the end ${\sf \Pi}^{\bm{f}}(\tau)_i$.

A non-increasing map $\pi : \llbracket 0,e\rrbracket \rightarrow \bbN$ is said to be \emph{admissible} if it is such that $\pi(e)=0,\pi(e-1)=1$ and $\pi(j)-\pi(j+1)\in \{0;1\}$ for every $0\leq j\leq e-1$. For any such $\pi$ and any integer $0\leq m \leq \pi(0)$ we define $j_\pi(m)$ as the smallest integer $j$ such that $\pi(j)=m$.

We associate to any $i\geq -1$, to any descending chain $\boldsymbol{\nu}: \nu_e\prec\cdots\prec\nu_0=\tau$, and to any admissible $\pi$ the element of ${\sf C}^{-\infty}$
$$
{\sf Q}_\pi \big( \tau/\nu_1, \dots , \nu_{e-1}/\nu_e , \nu_e \big)_{i_e} \defeq  \sum_{(i_j)_{0\leq j\leq e-1} \in  D_{\pi,i_e}} \; \prod_{j=0}^{e-1} {\sf g}^{\bm{f}}\big(\nu_j/\nu_{j+1}\big)_{i_{j}}  \, {\sf \Pi}^{\bm{f}}(\nu_e)_{i_e},
$$
where
$$
D_{\pi,i_e} \defeq \Big\{ (i_j)_{0\leq j\leq e-1}\in \llbracket -1,+\infty \llbracket^{e} ,\quad \forall j\in \llbracket 0,e-1 \rrbracket, \, i_j< i_{j_\pi( \pi(j)-1)}-1 \Big\}.
$$

 For every descending chain $\boldsymbol{\sigma}: \sigma_{e'}\prec\cdots\prec\tau$, from the identity of Lemma \ref{lem_gtogtilde} giving $[\sigma_j/\sigma_{j+1}]$ and $[\sigma_{e'}]$ in terms of ${\sf g}^{\bm{f}}$ and ${\sf \Pi}^{\bm{f}}$, developing the products gives the identity
$$
{\sf P}_<\big([\tau/\sigma_1],\dots,[\sigma_{e'}]\big)_i = \sum_{\boldsymbol{\nu},\pi}\lambda_{\boldsymbol{\sigma}}^{\boldsymbol{\nu},\pi} {\sf Q}_\pi \big( \tau/\nu_1,\dots,\nu_e \big)_i.
$$
where the sum runs over the set of descending sequences $\nu_e\prec\cdots\prec\nu_0=\tau$ and the set of admissible maps $\pi$, and where $ \lambda_{\boldsymbol{\sigma}}^{\boldsymbol{\nu},\pi}=0$ except if $\boldsymbol{\sigma}$ is a subsequence of $\boldsymbol{\nu}$ of size $e'$ such that $\pi(0)-e'\in \{  0,1  \} $  and $\sigma_{e'-m} = \nu_{j_\pi (m)}$ for every $0\leq m \leq e'$, in which case we have $ \lambda_{\boldsymbol{\sigma}}^{\boldsymbol{\nu},\pi}=(-1)^{e-e'}$. Then
$$
\sum_{e'\geq 0}\hspace{3pt}\sum_{\sigma_{q'}\prec\cdots\prec\sigma_1\prec\tau} {\sf P}_<\big([\tau/\sigma_1],\cdots,[\sigma_{e'}] \big)_i = \sum_{e\geq 0} \sum_{\nu_e\prec\cdots\prec\nu_1\prec\tau}\sum_\pi\lambda^{\boldsymbol{\nu},\pi} {\sf Q}_\pi \big( \tau/\nu_1,\cdots,\nu_e  \big)
$$
where $\lambda^{\bm{\nu},\pi} =  \sum_{\bm{\sigma}} \lambda_{\bm{\sigma}}^{\bm{\nu},\pi}$, for a sum over the set of finite descending sequences $\bm{\sigma}: \sigma_{e'}\prec\cdots\prec\tau$. We actually have $\lambda^{\bm{\nu},\pi}=0$ for every non-empty sequence $\bm{\nu}$. Indeed for any given $\bm{\nu}\neq \emptyset$ of size $e$ and any admissible $\pi$ there are only two descending sequences such that $\lambda_{\bm{\sigma}}^{\bm{\nu},\pi}\ne 0$. These sequences $\bm{\sigma}^1$ and $\bm{\sigma}^2$ are of size $\pi(0)$ and $\pi(0)-1$, respectively, and 
$$
\sigma^1_m = \nu_{j_\pi(\pi(0)-m)}
$$
and 
$$
\sigma^2_m = \nu_{j_\pi(\pi(0)-1-m)}.
$$
The two coefficient $\lambda_{\sigma}^{\boldsymbol{\nu},\pi}$ for these two $\boldsymbol{\sigma}$ are of opposite sign, which implies indeed that $\lambda^{\boldsymbol{\nu},\pi}=0$.
 \end{Dem}

\ssk

Last, before turning to the next section, we recall a variation on Lemma 6.6 of Gubinelli, Imkeller \& Perkowski's work \cite{GIP}.

\ssk

\begin{lem} \label{lem_critmodel}
Let $\sf \Pi$ be a linear map from $T$ to $\mcD'(\bbR^{d_0})$ and $\sf g$ be a map from $\bbR^{d_0}$ into the set of characters of the algebra $T^+$. The pair $({\sf \Pi},{\sf g)}$ is a model iff one has both
    \begin{equation} \label{eq_critmodelPi}
    \big|\langle{\sf \Pi}_x\tau , K_{<i,x} \rangle \big| \lesssim 2^{-i|\tau|} \qquad (\forall \, \tau\in \mcB)
    \end{equation}
uniformly over $i\geq-1$ and $x\in{\bbR^{d_0}}$, and 
    \begin{equation} \label{eq_critmodelg}
    |{\sf g}_{yx}(\mu)| \lesssim |y-x|^{|\mu|}  \qquad (\forall \, \mu\in \mcB^+)
    \end{equation}
uniformly on $(x,y)$ in any compact subset of $\bbR^{d_0}$.
\end{lem}

\medskip

\subsection{Proof of Theorem \ref{ThmMain1}. \hspace{0.15cm}}
\label{SectionProofTheorems}

For $\tau, \sigma\in T$ and a descending sequence $\bm{\sigma}(e) = (\sigma_e\prec\cdots\prec\sigma_1)$ we write
$$
\sigma \prec \bm{\sigma}(e) \prec \tau \textrm{ if }\big(\sigma \prec \sigma_e \textrm{ and } \sigma_1\prec\tau\big)
$$
and set
$$
|\tau/\bm{\sigma}(e)|_\alpha \defeq \big( |\tau/\sigma_1|_\alpha,\cdots,|\sigma_e/\sigma|_\alpha \big) \in\bbR^{e+1}.
$$
We prove by induction on $n$ the following three facts at a time.
{\it \begin{enumerate}
    \item[$\textit{(a)}_n$] For any tuple $\beta=(\beta_j)_{1\leq j\leq n}\in\bbR^n$ such that $\sum_{j=1}^n\beta_j>0$ the map 
    $$
    \big(g_1, \dots, g_n\big) \mapsto \widetilde{\sf P}^\beta_{\bbfl}(g_1,\dots,g_n)
    $$
    has a continuous extension from $\prod_{j=1}^nC^{\beta_j}_\circ$ into $L^\infty$.   \vspace{0.15cm}
\end{enumerate}
\hspace{0.35cm}For any tuples $\alpha=(\alpha_1, \dots, \alpha_n)$ and $\bm{f}=(f_1, \dots, f_n)$ of smooth functions:   \vspace{0.1cm}
\begin{enumerate}
    \item[$\textit{(b)}_n$] for any homogeneous $\tau=\llbracket a,b\rrbracket_{\bbfl} X^p\in T$ we have 
    $$
\big|\langle {\sf \Pi}^{\bm{f}}_x\tau, K_{<i,x} \rangle \big|\lesssim  \norme{f_a}_{\beta_a}\cdots\norme{f_b}_{\beta_b} \, 2^{-i|\tau|_\alpha},
    $$
    uniformly over $x\in{\bbR^{d_0}}$ and $i\geq 0$;   \vspace{0.15cm}
    
    \item[$\textit{(c)}_n$] for any homogeneous $\tau=\llbracket a,b\rrbracket_{\bbfl}^{\bbfk} X^p\in T^+$ we have
$$
\big| {\sf g}^{\bm{f}}_{y,x}(\tau) \big| \lesssim \norme{f_a}_{\beta_a}\cdots\norme{f_b}_{\beta_b} \, |y-x|^{|\tau|_\alpha},
$$
uniformly over $x,y\in{\bbR^{d_0}}$.
\end{enumerate}   }

\ssk

Theorem \ref{ThmMain1} follows as a consequence. The result holds true for $n=1$. We will use in the induction the following two algebraic identities proved in Appendix \ref{SectionProofLemmaRelStarDeriv}. 

\ssk

\begin{lem} \label{LemRelStarDerive}
We fix a tuple $\bm{f}=(f_1,\dots,f_n)$ of smooth functions.
\begin{enumerate}    
    \item[\textbf{(i)}] Pick $k\in\bbN^{d_0}$ with $\bbfk\in \mcP_n(k)$. Set $\partial^k\bm{f}=(\partial^{k_1}f_1, \dots, \partial^{k_n}f_n)$. We work here in the regularity structure $\mathscr{T}_{\alpha-\vert\bbfk\vert}$ with the pair of maps $({\sf \Pi}^{\partial^k\bm{f}},{\sf g}^{\partial^k\bm{f}})$ and its associated bracket maps $[\cdot]^{\partial^k\bm{f}}$. For $\tau = \llbracket 1,n\rrbracket_{\bbfl}X^m\in T$ with $|k|< |\llbracket 1,n\rrbracket_{\bbfl}|_{\alpha}$ we have
    $$
    \sum_{e\geq 0} \hspace{2pt}\sum_{\sigma_e\prec\cdots\prec\sigma_1\prec \tau} \widetilde{\sf P}_<^{|\tau/\bm{\sigma}(e)|_{\alpha-\vert\bbfk\vert}} \Big( [\tau/\sigma_1]^{\partial^k\bm{f}}, \dots, [\sigma_e]^{\partial^k\bm{f}}  \Big)
    =
    \mathbf{1}_{m=0} \, \widetilde{\sf P}_{\bbfl}^{\alpha-\vert\bbfk\vert}\big(\partial^{k_1}f_1,\cdots,\partial^{k_n}f_n\big).
    $$
    
        \item[\textbf{(ii)}] We work here in the regularity structure $\mathscr{T}_\alpha$ with the pair of maps $({\sf \Pi}^{\bm{f}}, {\sf g}^{\bm{f}})$ and its associated bracket maps $[\cdot]^{\bm{f}}$. For $\tau/\sigma = \llbracket 1,n\rrbracket^{\bbfk}_{\bbfl}X^m \in T^+$ and $|p|<|\tau/\sigma|_{\alpha}$ we have the identity
    \begin{align*}
    &\sum_{e\geq 0} \hspace{2pt}\sum_{\sigma\prec\bm{\sigma}(e)\prec\tau} \, \sum_{\bbfp\in\mcP_{e+1}(p)}\binom{p}{\bbfp} \, \widetilde{\sf P}_<^{|\tau/\bm{\sigma}(e)|_{\alpha} -\vert\bbfp\vert}\Big( \partial^{p_1} [\tau/\sigma_1]^{\bm{f}},\cdots,\partial^{p_{q+1}}[\sigma_e/\sigma]^{\bm{f}}  \Big)
    \\
    &\hspace{3cm}=\mathbf{1}_{m=0} \; \bigg\{ \sum_{\bbfk\in\mcP_n(k)}\sum_{\bbfp\in\mcP_n(p)}\binom{k}{\bbfk}\binom{p}{\bbfp} \, \widetilde{\sf P}_{\bbfl }^{\alpha-\vert \bbfk+\bbfp\vert}\Big(\partial^{k_1+p_1}f_1,\cdots,\partial^{k_n+p_n}f_n\Big) \bigg\}.
    \end{align*}    
\end{enumerate} 
\end{lem}

\ssk

We proceed with the induction step 
$$
\big((a)_{n-1},(b)_{n-1},(c)_{n-1}\big) \Longrightarrow \big((a)_n,(b)_n,(c)_n\big).
$$

\ssk

-- \textbf{\textit{We begin by proving $\textrm{(a)}_n$.}} Pick $\beta=(\beta_1,\dots,\beta_n)\in\bbR^n$ with $\sum_{i=1}^n\beta_i>0$. We work with the regularity structure $\mathscr{T}_\beta$. For $\bbfl=(\ell_1,\dots,\ell_{n-1})\in(\bbR^{d_0})^{n-1}$ set $\bbfl^- \defeq (\ell_1,\dots,\ell_{n-2})$ and
$$
\tau_n(\bbfl) \defeq \llbracket 1,n-1\rrbracket_{\bbfl^-} {X}^{\ell_{n-1}}.
$$
Write $\llbracket n\rrbracket$ for $\llbracket n,n\rrbracket\in T$ and $[n]^{\bm{g}}$ for its associated bracket map. From the continuity result of Proposition \ref{PropContinuityPTilde} for the $\widetilde{\sf P}_<^\gamma$ it suffices to prove that
\begin{equation} \label{eq_ptildetoptilde}
    \widetilde{\sf P}_{\bbfl }^{\beta}\big(g_1,\dots, g_n\big) = \sum_{e\geq 0}\hspace{2pt}\sum_{\sigma_e\prec\cdots\prec \tau_n(\bbfl )} \widetilde{\sf P}_{<}^{|\tau_n(\bbfl )/\boldsymbol{\sigma}^+|_{\beta}}\big([\tau_n(\bbfl )/\sigma_1]^{\bm{g}}, \dots, [\sigma_e]^{\bm{g}}, [n]^{\bm{g}} \big)
\end{equation}
where
$$
\big| \tau_n(\bbfl )/\boldsymbol{\sigma}^+ \big|_{\beta} \defeq \Big( |\tau_n(\bbfl )/\sigma_1|_{\beta}, \, |\sigma_1/\sigma_2|_{\beta}, \dots, \, |\sigma_e|_{\beta},\, \beta_n \Big),
$$
the symbol $+$ meaning that we added $\beta_n$ at the end of the uplet $|\tau_n(\bbfl )/\boldsymbol{\sigma}|_\beta$. Indeed, if one has \eqref{eq_ptildetoptilde}, the induction hypothesis and Proposition \ref{lemregcrochet} then ensure that any term $[\nu]^{\bm{g}}$ appearing in the right hand side of \eqref{eq_ptildetoptilde} is an element of ${\sf C}^{\vert\nu\vert}$ that depends continuously on $\bm{g}\in \prod_{j=1}^n C^{\beta_j}_\circ$. Since 
$$
|\tau_n(\bbfl )/\sigma_1|_\beta + |\sigma_1/\sigma_2|_\beta + \cdots + |\sigma_e|_\beta + \beta_n = \sum_{j=1}^n \beta_j > 0
$$
we can use Proposition \ref{PropContinuityPTilde} to conclude that $(a)_n$ holds true.

The remaining of this paragraph is dedicated to proving \eqref{eq_ptildetoptilde} by induction. The recursive definition of the $\widetilde{\sf P}_<^\alpha$ given in Lemma \ref{rectilde1} writes here
\begin{equation} \label{EqInductionStepThm1A} \begin{split}
\sum_{e\geq 0} \hspace{2pt} \sum_{\sigma_e\prec\cdots\prec \tau_n(\bbfl )} &\widetilde{\sf P}_{<}^{|\tau_n(\bbfl )/\boldsymbol{\sigma}^+|_{\beta}} \big([\tau_n(\bbfl )/\sigma_1]^{\bm{g}}, \dots, [\sigma_e]^{\bm{g}}, [n]^{\bm{g}} \big)
\\
&= \sum_{e\geq 0} \hspace{2pt} \sum_{\sigma_e\prec\cdots\prec \tau_n(\bbfl )} {\sf P}_{<} \big([\tau_n(\bbfl )/\sigma_1]^{\bm{g}}, \dots, [\sigma_e]^{\bm{g}}, [n]^{\bm{g}} \big)
\\
&\quad-\sum_{\substack{\sigma\prec\tau_n(\bbfl )\\ |\sigma|_{\beta}+\beta_n <0}} \sum_{\sigma\prec\sigma_{e_1}\prec\cdots\prec\tau_n(\bbfl )} \widetilde{\sf P}_{<}^{|\tau_n(\bbfl )/\boldsymbol{\sigma}|_{\beta}}\big([\tau_n(\bbfl )/\sigma_1]^{\bm{g}}, \dots, [\sigma_{e_1}/\sigma]^{\bm{g}} \big)
\\
&\hspace{3cm}\times \sum_{\nu_{e_2}\prec\cdots\prec\sigma} \widetilde{\sf P}_{<}^{|\sigma/\boldsymbol{\nu}^+|_{\beta}} \big([\sigma/\nu_1]^{\bm{g}}, \dots, [\nu_{e_2}]^{\bm{g}}, [n]^{\bm{g}} \big).
\end{split} \end{equation}
From Proposition \ref{PropoRepresentationProposition} one has 
$$
\sum_{e\geq 0} \hspace{2pt} \sum_{\sigma_e\prec\cdots\prec \tau_n(\bbfl )} {\sf P}_{<} \big([\tau_n(\bbfl )/\sigma_1]^{\bm{g}}, \dots, [\sigma_e]^{\bm{g}}, [n]^{\bm{g}} \big) = {\sf P}_{\bbfl}\big( g_1,\dots,g_n \big).
$$
Since any $\sigma\prec\tau_n(\bbfl )$ has the form $\sigma = \llbracket m+1 , n-1\rrbracket_{\bbfl +p_1}X^{p_2+s_2}$, one has $\tau_n(\bbfl )/\sigma = \llbracket 1,m\rrbracket_{\bbfl }^pX^{s_1}$, with $\ell_{n-1}=s_1+s_2$ and $p=p_1+p_2$. If $s_1=0$ item {\it (ii)} of Lemma \ref{LemRelStarDerive} ensures that the sum over the descending sequences $\sigma\prec\sigma_{e_1}\prec\cdots\prec\tau_n(\bbfl )$ is null. The terms $\sigma\prec\tau_n(\bbfl )$ that may give some non-trivial contributions to the sum \eqref{EqInductionStepThm1A} are thus of the form $\sigma = \llbracket m+1 , n-1\rrbracket_{\bbfl +p_1}X^{p_2+l_{n-1}}$, for which $\tau_n(\bbfl )/\sigma = \llbracket 1 , m\rrbracket_{\bbfl }^p$ with $p_1+p_2=p$. For such $\sigma$, item {\it (ii)} of Lemma \ref{LemRelStarDerive} gives 
$$
\sum_{\sigma\prec\sigma_{e_1}\prec\cdots\prec\tau_n(\bbfl )} \hspace{-0.2cm} \widetilde{\sf P}_{<}^{|\tau_n(\bbfl )/\boldsymbol{\sigma}|_{\beta}} \big([\tau_n(\bbfl )/\sigma_1]^{\bm{g}}, \dots, [\sigma_{e_1}/\sigma]^{\bm{g}} \big) = \sum_{\bbfp \in\mcP_m(p)}\binom{p}{\bbfp } \widetilde{\sf P}_{\bbfl _{<m}}^{\beta_{\leq m}-\vert\bbfk\vert}\big(\partial^{p_1}g_1,\dots,\partial^{p_m}g_m \big)
$$
and we have from the induction hypothesis
$$
\sum_{\nu_{e_2}\prec\cdots\prec\sigma} \widetilde{\sf P}_{<}^{|\sigma/\boldsymbol{\nu}^+|_{\beta_{>m}}} \big([\sigma/\nu_1]^{\bm{g}}, \dots, [\nu_{e_2}]^{\bm{g}}, [n]^{\bm{g}} \big) = \widetilde{\sf P}_{ J_{p_1,p_2}(\bbfl _{>m}) }^{\beta}\big(g_{m+1}, \dots, g_n \big)
$$
where 
$$
J_{p_1,p_2}(\bbfl _{>m}) = \sum_{\bbfa\in\mcP_{n-m-2}(p_1)}\binom{p_1}{\bbfa}\binom{k}{p_1}\big( \ell_{m+1}+a_1, \dots, \ell_{n-2}+a_{n-m-2}, \ell_{n-1}+p_2 \big).
$$
We recognize then in \eqref{EqInductionStepThm1A} the recursive relation satisfied by the $\widetilde{\sf P}_{\bbfl }^{\beta}$, which proves \eqref{eq_ptildetoptilde}.

\ssk

-- \textbf{\textit{We now turn to $\textrm{(b)}_n$.}} We would like to implement the same strategy as in point $(a)_n$: Write an iterated paraproduct as a sum of simplified iterated paraproducts and use their local expansion properties.. The problem with this strategy is that Proposition \ref{PropContinuityPTilde} requires some positivity assumption on some regularity exponents to hold -- which does not necessarily hold true here. To circumvent this issue,  for any $r\geq -1$, we look at the expansion properties of the iterated paraproduct ${\sf P}_{\bbfl }\big(f_1,\dots,f_{n-1},\Delta_r(f_n)\big)$ and treat $\Delta_r(f_n)$ as a function of high enough regularity in the estimates. We verify a posteriori that the remainders are summable over $r\geq -1$ and provide the right expression. 

We use the same notations as in the proof of point $(a)_n$. Pick $\alpha_n'> \alpha_n$ big enough such that $\sum_{s=j}^{n-1}\alpha_s+ \alpha_n' > 0$ for all $1\leq j < n$. Set 

$$
\bm{\alpha'}\defeq (\alpha_1,\dots,\alpha_{n-1},\alpha_n') 
$$ 
and, for any $r\geq -1$, let 
$$
\bm{f}^r \defeq \big(f_1,\dots, f_{n-1} , \Delta_r(f_n)\big)
$$
and
$$
{{\sf M}^r}' = \big( {{\sf \Pi}^r}', {{\sf g}^r}' \big) \defeq {\sf M}^{\bm{f}^r,\bm{\alpha'}}
$$ 
(The last notation was introduced at the beginning of Section \ref{SectionRepresentationFormula}.) One has
\begin{equation} \label{eq_prvthmb} \begin{split}
{\sf P}_{\bbfl }\big(f_1,\dots,\Delta_r f_n\big) &= \sum_{i\geq -1} \Delta_{<i-1}^{\ell_{n-1}}\big({\sf P}_{\bbfl^-}\big(f_1,\dots,f_{n-1}\big)\big) \Delta_i\big( \Delta_r(f_n) \big)
\\
&= \sum_{i\geq -1} \sum_{e\geq 0}  \hspace{2pt}\sum_{\sigma_e\prec\cdots\prec \tau_n(\bbfl )} {\sf P}_{<}\Big([\tau_n(\bbfl )/\sigma_1]^{ \sf M},\dots,[\sigma_e]^{\sf M} \Big)_{<i-1}\Delta_i \big( \Delta_r (f_n) \big)
\\
&= \sum_{e\geq 0}\hspace{2pt}\sum_{\sigma_e\prec\cdots\prec \tau_n(\bbfl )} {\sf P}_{<}\Big([\tau_n(\bbfl )/\sigma_1]^{{{\sf M}^r}'},\dots,[\sigma_e]^{{{\sf M}^r}'},[n]^{{{\sf M}^r}'} \Big).
\end{split}\end{equation}
From Proposition \ref{lemregcrochet} and the induction hypothesis, every term $[\sigma]$ appearing in the paraproduct ${\sf P}_<$ is an element of ${\sf C}^{|\sigma|}$ that depends continuously on $\bm{f}\in\big\{\prod_{j=1}^{n-1}C^{\alpha_j}_\circ\big\}\times C_\circ^{\alpha'_n}$. The assumption on $\alpha_n'$ ensures that the homogeneities of the element in the iterated paraproducts ${\sf P}_{<}\big([\tau_n(\bbfl )/\sigma_1]^{{{\sf M}^r}'},\dots,[\sigma_e]^{{{\sf M}^r}'},[n]^{{{\sf M}^r}'} \big)$ add up to a positive quantity, however $\vert\sigma_e\vert_{\bm{\alpha}}$ may be non-positive. This is cured by noticing that the assumption on $\alpha_n'$ ensures that
\begin{equation} \label{EqTrickPtoTildeP}
{\sf P}_{<}\Big([\tau_n(\bbfl )/\sigma_1]^{{{\sf M}^r}'},\dots,[\sigma_e]^{{{\sf M}^r}'},[n]^{{{\sf M}^r}'} \Big) = \widetilde{\sf P}_<^{|\tau_n(\bbfl )/\boldsymbol{\sigma}^+ |_{\boldsymbol{\alpha'}}} \Big([\tau_n(\bbfl )/\sigma_1]^{{{\sf M}^r}'}, \dots, [\sigma_e]^{{{\sf M}^r}'}, [n]^{{{\sf M}^r}'} \Big),
\end{equation}
where 
$$
|\tau_n(\bbfl )/\boldsymbol{\sigma}^+|_{\boldsymbol{\alpha'}} \defeq \big( \, |\tau_n(\bbfl )/\sigma_1|_{\boldsymbol{\alpha'}},\, |\sigma_1/\sigma_2|_{\boldsymbol{\alpha'}},\dots,|\sigma_e|_{\boldsymbol{\alpha'}},\, \alpha_n'  \, \big),
$$ 
so one can use Proposition \ref{PropExpansionTildePAlphaMinus} on the local expansion of terms of the type $\widetilde{\sf P}_<^\gamma$. The remainder term in \eqref{EqRemainderTildePMinusAlpha} is $(\triangle_{y-x,\theta}\widetilde{\sf P}_<^\gamma)(\dots)(x)$ with $\theta=\sum_{j=1}^n \gamma_j$. We infer from this generic expansion property, \eqref{eq_prvthmb} and \eqref{EqTrickPtoTildeP}, that ${\sf P}_{\bbfl }(f_1,\dots,\Delta_r(f_n))$ has a corresponding expansion with remainder
$$
\sum_{e\geq 0}\hspace{2pt}\sum_{\sigma_e\prec\cdots\prec \tau_n(\bbfl )} \big(\triangle_{y-x,\theta}\widetilde{\sf P}_{<}^{|\tau_n(\bbfl )/\boldsymbol{\sigma}^+|_{\boldsymbol{\alpha'}}}\big) \Big([\tau_n(\bbfl )/\sigma_1]^{{{\sf M}^r}'}, \dots, [\sigma_e]^{{{\sf M}^r}'}, [n]^{{{\sf M}^r}'} \Big)(x),
$$
with $\theta = \sum_{j=1}^n \alpha_j'-\delta$, here. The following result is proved in Appendix \ref{SectionBasicsAnalysis} by induction on $n$.

\ssk

\begin{lem} \label{LemProofPointB}
For every point $x\in{\bbR^{d_0}}$ and $i\geq-1$ one has the identity
\begin{equation} \label{eq_triangletopi} \begin{split}
\int_{{\bbR^{d_0}}} K_{<i}(h) \sum_{e\geq 0}\hspace{2pt}\sum_{\sigma_e\prec\cdots\prec \tau_n(\bbfl )} 
&\big(\triangle_{h,\theta} \widetilde{\sf P}_{<}^{|\tau_n(\bbfl )/\boldsymbol{\sigma}^+|_{\boldsymbol{\alpha'}}}\big) \Big([\tau_n(\bbfl )/\sigma_1]^{{{\sf M}^r}'}, \dots, [\sigma_e]^{{{\sf M}^r}'}, [n]^{{{\sf M}^r}'} \Big)(x) \, dh
\\
&=\int_{{\bbR^{d_0}}} K_{<i}(h) \big({{\sf \Pi}_x^r}'\llbracket 1,n\rrbracket_{\bbfl}\big)(x+h) \, dh.
\end{split} \end{equation}
\end{lem}

\ssk

We then have from \eqref{eq_triangletopi} and Lemma 6.3 in \cite{GIP}
\begin{equation} \label{EqAlmostPointBn} \begin{split} 
\left|\int_{{\bbR^{d_0}}} K_{<i}(x-y) \big({\sf \Pi}'^{(r)}_x \llbracket 1,n\rrbracket_{\bbfl }\big)(y) \, dy \right| 
&
\lesssim \bigg\{\prod_{j=1}^{n-1} \norme{f_j}_{\alpha_j}\bigg\} \norme{\Delta_r f_n}_{\alpha_n'}  2^{-i\theta}
\\
&\lesssim 2^{-i\theta} \, 2^{r(\alpha_n'-\alpha_n)} \bigg\{ \prod_{j=1}^{n} \norme{f_j}_{\alpha_j} \bigg\}.
\end{split} \end{equation}
There is an integer $i(n)$ depending only on $n$ such that we have for $j\leq n$ and $i\geq -1$
$$
\Delta_i \big({\sf P}_\bbfl(f_1,\dots,f_j)\big) = \sum_{r\leq i+i(n)} \Delta_i \big( {\sf P}_{\bbfl}\big(f_1, \dots, f_{j-1}, \Delta_r(f_j)\big)\big).
$$
Using the identity \eqref{eq_devPi} on ${\sf \Pi}_x(\tau)$, we see that we have
$$
\Big\langle {\sf \Pi}_x(\llbracket 1,n\rrbracket_{\bbfl}) , K_{<i,x} \Big\rangle = \sum_{r\leq i+i(n)} \Big\langle {{\sf \Pi}_x^r}'\big(\llbracket 1,n\rrbracket_{\bbfl }\big),K_{<i,x} \Big\rangle,
$$
so the expected bound
\begin{align*}
\Big| \Big\langle {\sf \Pi}_x(\llbracket 1,n\rrbracket_{\bbfl}) , K_{<i,x} \Big\rangle \Big| &\lesssim \bigg\{\prod_{j=1}^n \norme{f_j}_{\alpha_j} \bigg\} \sum_{r=-1}^{i+i(n)} 2^{r(\alpha_n'-\alpha_n)} \, 2^{-i\theta} \lesssim \bigg\{\prod_{j=1}^n \norme{f_j}_{\alpha_j}\bigg\} \, 2^{-i{\sum_{j=1}^n\alpha_j}}
\end{align*}
follows from \eqref{EqAlmostPointBn}.

\ssk

-- \textbf{\textit{We finally prove $\textrm{(c)}_n$.}} Pick $\alpha=(\alpha_1,\dots,\alpha_n)\in\bbR^n$, a multi-indice $k\in\bbR^{d_0}$ such that $|k|<\sum_{j=1}^n \alpha_j$ and $\bbfk \in\mcP_n(k)$. We work in the regularity structure $\mathscr{T}_{\alpha-\vert\bbfk\vert}$. From item {\it (ii)} of Lemma \ref{LemRelStarDerive} we have for any smooth functions $f_1,\dots f_n$ the equality
$$
\widetilde{\sf P}_{\bbfl }^{\alpha-\vert\bbfk\vert} \big(\partial^{k_1}f_1,\dots,\partial^{k_n}f_n\big) = \sum_{e\geq 0} \hspace{2pt}\sum_{\sigma_e\prec\cdots\prec\sigma_1\prec \llbracket 1,n\rrbracket_{\bbfl }} \widetilde{\sf P}_<^{|\llbracket 1,n\rrbracket_{\bbfl }/\boldsymbol{\sigma}|_{\alpha-\vert\bbfk\vert}}\Big( [\llbracket 1,n\rrbracket_{\bbfl }/\sigma_1]^{{\sf M}_{\bbfk}},\dots,[\sigma_e]^{{\sf M}_{\bbfk}}  \Big),
$$
where ${\sf M}_{\bbfk}={\sf M}_{\partial^{\bbfk}\bm{f} , \alpha-\vert\bbfk\vert}$. Proposition \ref{lemregcrochet} and point $(b)_n$ ensure by induction that all the terms $[\nu]^{{\sf M}_{\bbfk}}$ are some elements of ${\sf C}^{|\nu|_{\alpha-\bbfk}}$ that depend continuously on all the $f_j\in C^{\alpha_j}_\circ$. As above it follows from Proposition \ref{PropExpansionTildePAlphaMinus} that $\widetilde{\sf P}_{\bbfl }^{\alpha-\vert\bbfk\vert} \big(\partial^{k_1}f_1,\dots,\partial^{k_n}f_n\big)$ has a local expansion with remainder 
$$
R_{\bm{f},\alpha}\big(\llbracket 1,n\rrbracket_{\bbfl }^{\bbfk}\big)(x,h) \defeq \sum_{e\geq 0} \sum_{\sigma_e\prec\cdots\prec \llbracket 1,n\rrbracket_{\bbfl }}
\big(\triangle_{h,\theta} \widetilde{\sf P}_<^{|\llbracket 1,n\rrbracket_{\bbfl }/\boldsymbol{\sigma}|_{\alpha-\vert\bbfk\vert}}\big)  \Big( \llbracket 1,n\rrbracket_{\bbfl }/\sigma_1]^{{\sf M}_{\bbfk}},\dots, [\sigma_e]^{{\sf M}_{\bbfk}} \Big)(x)
$$
where $\theta = |\llbracket 1,n\rrbracket_{\bbfl }|_{{\alpha-\vert\bbfk\vert}}$. From Proposition \ref{PropExpansionTildePAlphaMinus} this remainder has $|h|^\theta \prod_{j=1}^n \norme{f_j}_{\alpha_j}$ as an $x$-uniform upper bound. Point $(c)_n$ will thus be proved after we show that for any $\tau = \llbracket 1,n\rrbracket_{\bbfl }X^s$ one has 
\begin{equation} \label{EqLastBit} \begin{split}
\sum_{e\geq 0 } \sum_{\sigma_e\prec\cdots\prec \tau}
\big( \triangle_{h,\theta} \widetilde{\sf P}_<^{|\tau/\boldsymbol{\sigma}|_{\alpha-\vert\bbfk\vert}} \big) \Big( [ \tau/\sigma_1]^{{\sf M}_{\bbfk}},\dots, [\sigma_e]^{{\sf M}_{\bbfk}} \Big)(x) = \mathbf{1}_{s=0}\enskip {\sf g}^{\bm{f}}_{x+h,x}\big(\llbracket 1,n\rrbracket_{\bbfl }^{\bbfk}\big).
\end{split}
\end{equation}
The remainder of this paragraph is dedicated to proving this identity by induction on $n$. Recall that we write
\begin{align*}
\partial^p_\star{\sf P}_<\big([\tau/\sigma_1]^{{\sf M}_{\bbfk}}, \dots, &[\sigma_{m-1}/\sigma_{m}]^{{\sf M}_{\bbfk}} \big)   \\
&= \sum_{\bbfp \in\mcP_m(p)} \binom{p}{\bbfp} \, \widetilde{\sf P}_<^{|\tau/\boldsymbol{\sigma}_{\leq m}|_{\alpha-\vert\bbfk\vert} - \vert\bbfp\vert} \big([\tau/\sigma_1]^{{\sf M}_{\bbfk}},\dots,[\sigma_{m-1}/\sigma_{m}]^{{\sf M}_{\bbfk}} \big).
\end{align*}
From the definition of $\triangle_{h,\theta}\widetilde{\sf P}_<$ the left hand side of \eqref{EqLastBit} is equal to 
\begin{equation*} \begin{split}
 &\sum_{e\geq 0}\hspace{2pt}\sum_{\sigma_e\prec\cdots\prec \tau} 
   \widetilde{\sf P}_{<}^{|\tau/\boldsymbol{\sigma}|_{\alpha-\vert\bbfk\vert}}\big([\tau/\sigma_1]^{{\sf M}_{\bbfk}},\dots,[\sigma_e]^{{\sf M}_{\bbfk}}\big)(x+h)
\end{split} \end{equation*}

\begin{equation} \label{eqprvthmc} \begin{split}
  &- \sum_{e\geq 0}\hspace{2pt}\sum_{\sigma_e\prec\cdots\prec \tau} \sum_{|p|<|\tau|_{\alpha-\vert\bbfk\vert}} \partial^p_\star\widetilde{\sf P}_<\big(\, 
 [\tau/\sigma_1]^{{\sf M}_{\bbfk}},\dots,[\sigma_e]^{{\sf M}_{\bbfk}} \big)(x) \, h^p 
\\
  & - \underset{\sigma_e\prec\cdots\prec \tau}{\sum_{e\geq 0}} \sum_{m=1}^e \sum_{|p|<|\tau/\sigma_m|_{\alpha-\vert\bbfk\vert }} \partial^p_\star{\sf P}_<\Big([\tau/\sigma_1]^{{\sf M}_{\bbfk}},\dots,[\sigma_{m-1}/\sigma_{m}]^{{\sf M}_{\bbfk}} \Big) (x)
  \\
  &\hspace{4,5cm} \times \frac{h^p}{p!} \, \big(\triangle_{h,|\sigma_m|_{\alpha-\vert\bbfk\vert}} \widetilde{\sf P}_<^{|\sigma_m/\boldsymbol{\sigma}_{>m}|_{\alpha-\vert\bbfk\vert }} \big)  \Big([\sigma_m/\sigma_{m-1}]^{{\sf M}_{\bbfk}},\dots ,[\sigma_e]^{{\sf M}_{\bbfk}}\big)(x),
\end{split} \end{equation}
From item {\it (i)} of Lemma \ref{LemRelStarDerive}, the first double sum in \eqref{eqprvthmc} is equal to 
$$
\mathbf{1}_{s=0}\, \widetilde{\sf P}_{\bbfl }^{\alpha-\vert\bbfk\vert }\big(\partial^{k_1}f_1,\dots,\partial^{k_n}f_n\big)(x+h) = \mathbf{1}_{s=0} \; {\sf g}_{x+h}^{\bm{f}}\big(\llbracket 1,n\rrbracket_{\bbfl }^{\bbfk}\big).
$$ 
Lemma \ref{LemRelStarDerive} also gives that the second line of \eqref{eqprvthmc} is equal to 
$$
\mathbf{1}_{s=0}\,\sum_{|p|<|\tau|_{\alpha-\vert\bbfk\vert }}{\sf g}_{x}^{\bm{f}}\big(\llbracket 1,n\rrbracket_{\bbfl }^{\bbfk+p}\big) \, h^p.
$$ 
The $\sigma \in T$ such that $\sigma\prec\tau = \llbracket 1,n\rrbracket_{\bbfl }X^s$ have a form $\sigma = \llbracket m+1,n\rrbracket_{\bbfl +v_1}X^{v_2+s_2}$, in which case $\tau/\sigma=\llbracket 1,m\rrbracket^{v}_{\bbfl }X^{s_1}$ with $s=s_1+s_2$ and $v=v_1+v_2$. For such $\sigma\in T$ Lemma \ref{LemRelStarDerive} gives 
\begin{align*}
\sum_{e_1\geq 0}\sum_{\sigma_{e_1}\prec\cdots\prec \tau} \partial^p_\star{\sf P}_<\Big([\tau/\sigma_1]^{{\sf M}_{\bbfk}}, \dots, [\sigma_{e_1}/\sigma]^{{\sf M}_{\bbfk}}\Big)(x)
&= \mathbf{1}_{s_1=0} \; {\sf g}^{\bm{f}}_x\big(\llbracket 1,m\rrbracket_{\bbfl }^{\bbfk+v+p}\big).
\end{align*}
Also, we have by induction that
$$
\sum_{e_2\geq 0}\sum_{\nu_{e_2}\prec\cdots\prec \sigma} \triangle_{h,|\sigma|_{\alpha-\vert\bbfk\vert}} {\sf P}_<\Big([\sigma/\nu_1]^{{\sf M}_{\bbfk}},\dots,[\nu_{e_2}]^{{\sf M}_{\bbfk}}\Big)= \mathbf{1}_{v_2+s_2=0} \; {\sf g}_{x,x+h}^{\bm{f}} \big(\llbracket m+1,n\rrbracket_{\bbfl +v_1}^{\bbfk}   \big).
$$
If $s\ne 0$, then either $s_1\ne 0$ or $s_2+p_2\ne0$, and all the terms of \eqref{eqprvthmc} add up to 0. Suppose now that $s=0$. The $\sigma\in T$ we have to consider are of the form $\sigma = \llbracket m+1,n\rrbracket_{\bbfl +v}$, for which $\tau/\sigma=\llbracket 1,m\rrbracket^v_{\bbfl}$ and the right hand side of \eqref{eqprvthmc} writes as
\begin{align*}
    R_{\bm{f},\alpha}\big(\llbracket 1,n\rrbracket_{\bbfl }^{\bbfk}\big) &= 
    {\sf g}_{x+h}^{\bm{f}}\big(\llbracket 1,n\rrbracket_{\bbfl }^{\bbfk}\big) - \sum_{|p|<|\llbracket 1,n\rrbracket_{\bbfl }^{\bbfk}|_{\alpha}} {\sf g}_{x}^{\bm{f}}\big(\llbracket 1,n\rrbracket_{\bbfl }^{\bbfk+p}\big) \, h^p
    \\
    &\quad- \sum_{m,p,v} {\sf g}_x^{\bm{f}}\big(\llbracket 1,m\rrbracket_{\bbfl }^{\bbfk+p+v}\big) {\sf g}^{\bm{f}}_{x+h,x}\big(\llbracket m+1,n\rrbracket_{\bbfl +v}^{\bbfk}X^p \big),
    \end{align*}
where the sum over $m,p,v$ runs over $1\leq m<n$ and multi-indices $p,v$ such that $\big|\llbracket 1,m\rrbracket_{\bbfl }^{v+p}\big|_{\alpha-\vert\bbfk\vert }>0$ and $\big|\llbracket m+1,n\rrbracket_{\bbfl +v}\big|_{\alpha-\vert\bbfk\vert }>0$. This sum corresponds to a sum over $\sigma\in T^+$ such that $\sigma\prec\llbracket 1,n\rrbracket^{\bbfk}_{\bbfl }$ in the regularity structure $\mathscr{T}_\alpha$. It follows that we finally have
\begin{align*}
R_{\bm{f},\alpha}\big(\llbracket 1,n\rrbracket_{\bbfl }^{\bbfk}\big)&={\sf g}_{x+h}^{\bm{f}} \big( \llbracket 1,n\rrbracket_{\bbfl }^{\bbfk} \big)  -  \sum_{\sigma<\llbracket 1,n\rrbracket_{\bbfl }^{\bbfk}} {\sf g}_{x+h,x}^{\bm{f}}(\sigma) {\sf g}_x^{\bm{f}}(\tau/\sigma) = {\sf g}_{x+h,x}^{\bm{f}} \big( \llbracket 1,n\rrbracket_{\bbfl }^{\bbfk}\big),
\end{align*}
which concludes the proof of \eqref{EqLastBit}, and closes the induction step in the proof of point $\textit{(c)}_n$.

\medskip

\section{Back to paracontrolled systems}\label{section_parapsystem}
\label{SectionPCSystems}

This section is dedicated to proving Theorem \ref{ThmMain2}. We set ourselves in the setting of Section \ref{IntroUniversalStructure}, with its finite alphabet $\mcL=(l_1,\dots,l_{\vert\mcL\vert})$ and its associated set $\mcW$ of finite words $w=l_{i_1}\dots l_{i_w}$. An a priori notion of size $\vert\cdot\vert_\mcL$ is given on $\mcL$, extended to $\mcW$ setting $\vert l_{i_1}\dots l_{i_w}\vert_\mcL = \vert l_{i_1}\vert_\mcL + \cdots +\vert l_{i_w}\vert_\mcL$.


\subsection{The regularity structure $\mathscr{T}_\mcL$.}
\label{SectionRSPCSystems}

\hspace{0.1cm} The following construction is identical to the construction of Section \ref{RSIPP}. We define a set of symbols
$$
\mcB \defeq \Big\{  [w]_{\bbfl} X^p \,;\, w=l_{i_1}\dots l_{i_j}\in\mcW ,\, p,\ell\in\bbN^{d_0}, \, \bbfl \in \mcP_{j-1}(\ell) \Big\} \cup \{X^k\}_{k\in\bbN^{d_0}},
$$
and 
\begin{align*}
\mcB^+ &\defeq \Big\{ [w]_{\bbfl }^{\bbfk} \,;\,  w=l_{i_1}\dots l_{i_j}\in\mcW, \, k,\ell\in\bbN^{d_0},\, \bbfl \in \mcP_{j-1}(\ell),\, \bbfk\in\mcP_j(k),  \; |w|_{\mcL}-|k|+|\ell|>0  \Big\}
\\
&\hspace{10cm}\cup \big\{X^{e_i}\big\}_{1\leq i\leq d_0}.
\end{align*}
We let $T$ be the vector space freely generated by $\mcB$, and $T_+$ be the algebra freely generated by $\mcB^+$, with unit ${\bf 1}^+$. We also set
$$
\big| [w]_{\bbfl } X^p \big|_{\mcL} \defeq |w|_{\mcL} + |\ell| + |p|
$$
and define $\vert\cdot\vert_\mcL$ on $T^+$ as a multiplicative function such that $|X^{e_i}|_{\mcL}=1$ and
$$
\big| [w]_{\bbfl }^{\bbfk} \big|_{\mcL} \defeq |w|_{\mcL} + |\ell| - |k|.
$$
Proceeding as in Section \ref{RSIPP}, for $\tau=[l_{i_1}\cdots l_{i_n}]_{\bbfl} \in T$ we set
$$
\oplus(\tau)\defeq\Big\{ [l_{i_1}\cdots l_{i_j}]_{\bbfl _{<j}}^{\bbfp }\in \mcB^+ \,;\, 1\leq j\leq n,\; \ell_j=0, \bbfp\in\mcP_j(p), p\in\bbN^{d_0} \Big\}\cup\{\mathbf{1}^+\},
$$
and for $\mu=[l_{i_1}\cdots l_{i_n}]_{\bbfl}^{\bbfk} \in T^+$ we set
$$
\oplus(\mu)\defeq\Big\{   [l_{i_1}\cdots l_{i_j}]_{\bbfl _{<j}}^{\bbfk+\bbfp }\in \mcB^+  \,;\, 1\leq j\leq n,\; \ell_j=0, \bbfp\in\mcP_j(p), p\in\bbN^{d_0} \Big\}\cup\{\mathbf{1}^+\}.
$$
Set for $\tau=[l_{i_1}\cdots l_{i_n}]_{\bbfl } \in T$ and $\sigma=[l_{i_1}\cdots l_{i_j}]_{\bbfl }^{\bbfp }\in \oplus(\tau)$ with $j<n$
$$
(\tau \backslash \sigma) \defeq \sum_{p=p_1+p_2}  \binom{p}{p_1}  [l_{i_{j+1}}\cdots l_{i_n}]_{\bbfl +p_1}X^{p_2},
$$
and for $j=n$ set $(\tau \backslash \sigma) \defeq \frac{1}{\bbfp !} X^p$.
For $\mu=[l_{i_1}\cdots l_{i_n}]_{\bbfl }^{\bbfk}$ and $\nu=[l_{i_1}\cdots l_{i_j}]_{\bbfl }^{\bbfk +\bbfp }\in \oplus(\mu)$ with $j<n$
$$
(\mu\backslash\nu) \defeq \sum_{p=p_1+p_2}  \binom{p}{p_1}  [l_{i_{j+1}}\cdots l_{i_n}]^{\bbfk}_{\bbfl +p_1}X^{p_2},
$$
and for $j=n$ set $(\mu \backslash \nu) \defeq \frac{1}{\bbfp !} X^p$. Finally set
$$
\Delta(\tau) \defeq \sum_{\sigma\in\oplus(\tau)} (\tau \backslash \sigma) \otimes \sigma,
$$
and
$$
\Delta^+(\mu) \defeq \sum_{\substack{\nu\in\oplus(\mu) \\ {|(\mu\backslash\nu)|_\mcL > 0}}} (\mu\backslash\nu) \otimes \nu.
$$
Proceeding as in Section \ref{RSIPP} shows that 
$$
\mathscr{T}_\mcL = \big((T,\Delta),(T^+,\Delta^+)\big)
$$ 
is a concrete regularity structure. Given $\alpha=(\alpha_1,\dots,\alpha_{\vert\mcL\vert})$ with $\sum_{j=1}^{\vert\mcL\vert}\alpha_j > 0$, and $[l]\in C^{\alpha_l}_\circ$ for all $1\leq j\leq n$, we define from Theorem \ref{ThmMain1} a model on $\mathscr{T}_\mcL$ setting
$$
{\sf \Pi}( [l_{i_1}\dots l_{i_n}]_{\bbfl }) \defeq {\sf P}_{\bbfl }\big( [l_{i_1}], \dots, [l_{i_n}] \big),
$$
$$
{\sf g}\big( [l_{i_1}\dots l_{i_n}]_{\bbfl }^{\bbfk} \big) \defeq \widetilde{\sf P}^{(\vert l_{i_1}\vert_\mcL,\dots, \vert l_{i_n}\vert_\mcL) - \vert\bbfk\vert}_{\bbfl }\big( \partial^{k_1}[l_{i_1}], \dots, \partial^{k_n} [l_{i_n}]\big),
$$
and ${\sf \Pi}( [l_{i_1}\dots l_{i_n}]_{\bbfl } X^p)(\cdot) = \cdot^p \, {\sf \Pi}( [l_{i_1}\dots l_{i_n}]_{\bbfl })(\cdot)$, with the notation \eqref{EqPiMap1}.

\medskip

\subsection{Paracontrolled systems and modelled distributions.}
\label{SectionPCSystemsModelledDistributions}

\hspace{0.1cm} We prove Theorem \ref{ThmMain2} in the refined form of Theorem \ref{ThmMain2Refined}. Let $r>0$ and $(u_w)_{w\in\mcU_{<r}}$ be a system $r$-paracontrolled by the $[l]$, as in \eqref{EqPCSystem}. For each $\tau=[w]_{\bbfl}X^{p}\in \mcB$, with $w=l_{j_1}\dots l_{j_w}$ and $\bbfl\in\mcP_{j_w-1}(\ell)$ with $l,p\in\bfN^{d_0}$ such that $|\tau|_\mcL < r$, set
\begin{align*}
u_\tau \defeq \sum_{\substack{ w'=l_{i_2}\cdots l_{i_n}\in\mcW \\ ww'\in\mcW_{<r} }} \, \sum_{\bbfk\in\mcP_n(\ell+p)}\binom{k}{\bbfk} \, \widetilde{\sf P}^{(\gamma-|ww'|, |l_{i_2}|, \dots, |l_{i_n}|) - \vert\bbfk\vert} \Big( \partial^{k_1}f_{s's}^\#,\,\partial^{k_2}[l_{i_2}], \dots, \partial^{k_n}[l_{i_n}] \Big),
\end{align*}
From Theorem \ref{ThmMain1}, each $u_\tau$ defines a bounded function as $r - |\tau|_\mcL > 0$. Define the $T$-valued function
$$
\bm{u}(x) \defeq \sum_{\tau\in\mcB} u_\tau(x) \, \tau.
$$

\begin{thm} \label{ThmMain2Refined}
One has $\bm{u} \in \mcD^r(T,{\sf g})$.
\end{thm}

\ssk

\begin{Dem}
We use Theorem \ref{ThmMain1} to prove that statement, but in a regularity structure that takes into account the $u_w^\sharp$ on the same footing as the $[l]$. We introduce for that purpose a new alphabet 
$$
\mcA \defeq \mcL \sqcup \mcW
$$ 
and set $|\lambda|_{\mcA} \defeq |a|_\mcA$ for $\lambda\in\mcL$ and $|\lambda|_{\mcA} \defeq r-|\lambda|_\mcA$ for $\lambda\in\mcW$. We write $\mcW_\mcA$ for the set of words written with the alphabet $\mcA$. To avoid any confusion when writing words with the alphabet $\mcA$ we will write $(w)$ the letter of $\mcA$ associated with $w\in\mcW$. We extend our collection $([l])_{l\in\mcL}$ into $([\lambda])_{\lambda\in\mcA}$ setting $[w] \defeq u_w^\sharp \in C^{r-\vert w\vert}$ for $\lambda=w\in\mcW$. As above, Theorem \ref{ThmMain1} provides a regularity structure associated with $\mcA$ and a model $\overline{\sf M}=(\overline{\sf \Pi}, \overline{\sf g})$ on it associated with $([\lambda])_{\lambda\in\mcA}$. There is a canonical injection $\iota :\mathscr{T}_\mcL \hookrightarrow \mathscr{T}_\mcA$ that commutes with the coproducts, and $\overline{\sf M}$ is an extension ${\sf M}$. We can thus freely pass from $\overline{\sf g}$ to ${\sf g}$ in some computations below.

\ssk

Within $\mathscr{T}_\mcA$, working with $\overline{\sf M}$, for $\tau=[w]_{\bbfl }X^p$ one can rewrite
\begin{equation*}
\overline{u}_\tau(x) = \sum_{w' \in\mcW} \overline{\sf g}_x\big([(w'w)w'\, ]^{l+p} \big).
\end{equation*}
For $w\in\mcW$ we let 
$$
\rho(w) \defeq [(w)w]
$$ 
where $(w)w\in\mcA$ is the word beginning with the letter $(w)\in\mcA$ followed by $w\in\mcA$ -- for $w=l_{i_1}\cdots l_{i_n}$ it represents the function ${\sf P}\big(u_w^\sharp, [l_{i_1}], \dots, [l_{i_n}]\big)$. Then $\bm{u}$ can be re-written in $\mathscr{T}_\mcA$ under the form
$$
\bm{u}(x) = \sum_{w\in\mcU_{<r}} \sum_{\sigma < \rho(w)} \overline{\sf g}_x(\rho(w)/\sigma) \, \sigma,
$$
Note that any $\sigma<\rho(w)$ has form $[w']_{\bbfl}X^p$ where $w'$ is a subword of $w$. We now prove that $\bm{u} \in \mcD^r(T,{\sf g})$ by proving that for any $w\in\mcU_{<r}$ the map $\bm{h}_w(x) = \sum_{\sigma<\rho(w)} \overline{\sf g}_x(\rho(s)/\sigma) \, \sigma$ is an element of $\mcD^{\vert\rho(w)\vert}(T_\mcA,\overline{g})$. For any $x,y\in{\bbR^{d_0}}$, and for $\tau=[w]_{\bbfl }X^p$, one has
\begin{align*}
\widehat{\overline{\sf g}_{yx}}(\bm{h}_w(x)) &= \sum_{\nu\leq\sigma<\rho(w)}  \overline{\sf g}_x\big(\rho(w)/\sigma \big) \, {\sf g}_{yx}(\sigma/\nu) \, \nu = \sum_{\nu<\tau} \Big( \overline{\sf g}_y(\rho(w)/\nu) - \overline{\sf g}_{yx}(\rho(w)/\nu) \Big) \, \nu
\\
&= \bm{h}_w(y) - \sum_{\nu<\tau} \overline{\sf g}_{yx}\big(\rho(w)/\nu\big) \, \nu.
\end{align*}
Theorem \ref{ThmMain1} ensures that $\big| \overline{\sf g}_{yx}(\rho(w)/\nu) \big| \lesssim |y-x|^{|\rho(w)/\nu|}$, with $|\rho(w)/\nu| = r - |\nu|_\mcA$.
\end{Dem}

\ssk

Paracontrolled systems in the generality of Section \ref{IntroUniversalStructure} were first introduced in Bailleul \& Mouzard's work \cite{BailleulMouzard}.


\appendix
\section{Appendix}

\subsection{Basics in regularity structures.}
\label{SectionBasicsRS}

\hspace{0.15cm} We recall here some basic facts about regularity structure. We refer the reader to \cite{RSGuide} for a thorough introduction to the subject, and to \cite{Hai14} for the original work of M. Hairer on the subject.

\ssk

\begin{defn} \label{DefnRS}
    A \textbf{concrete regularity structure} is a pair $\mathscr{T} = (T,T^+)$ of graded vector spaces 
    $$
    T=\bigoplus_{r\in A} T_r, \qquad T^+=\bigoplus_{s\in A^+} T_s^+,
    $$
    such that the following holds.
    
    \begin{itemize}
        \item[--] The spaces $T_r$ and $T_s^+$ are finite dimensional for any $r\in A$ and $s\in A^+$. One has $A^+\subset [0,+\infty)$ and both $A$ and $A^+$ are bounded from below and have no accumulation points.   \vspace{0.1cm}
        \item[--] The vector space $T^+$ is a connected graded bialgebra with coproduct $\Delta^+$ and grading $A^+\subset [0,+\infty[$.   \vspace{0.1cm}
        \item[--] The vector space is endowed with a linear splitting map $\Delta : T\rightarrow T\otimes T^+$ such that
        $$
        (\Delta\otimes \emph{Id})\Delta = (\emph{Id}\otimes \Delta^+)\Delta.
        $$
        \item[--] We have
        $$
        \Delta T_{r_1} \subset \bigoplus_{r_2\in A} T_{r_2} \otimes T^+_{r_1-r_2}, \qquad \Delta^+ T_{s_1}^+ \subset \bigoplus_{s_2\in A^+} T_{s_2}^+ \otimes T^+_{s_1-s_2}.
        $$
    \end{itemize}
\end{defn}

\ssk

We suppose here that the vector spaces $T$ and $T^+$ come with some bases $\mcB$ and $\mcB^+$. Then for any $\tau\in T$ we have a decomposition 
$$
\Delta \tau= \sum_{\sigma\in \mcB} (\tau/\sigma) \otimes \sigma
$$
for some elements $\tau/\sigma\in T$. Likewise we define $\tau/\sigma \in T^+$ for $\tau\in T^+$ and $\sigma\in\mcB^+$ from the identity
$$
\Delta^+ \tau= \sum_{\sigma\in \mcB^+} (\tau/\sigma) \otimes \sigma.
$$
For $\sigma,\tau\in\mcB$ we write $\sigma\leq\tau$ if $\tau/\sigma \ne 0$ and $\sigma<\tau$ if $\sigma$ and $\tau$ are distinct and $\sigma\leq\tau$. For $\tau,\sigma,\nu\in \mcB$ we have
$$
\Delta^+(\tau/\sigma) = \sum_{\sigma\leq \nu \leq \tau}  \tau/\nu \otimes \nu/\sigma. 
$$
We denote by $G^+$ the set of real-valued characters of the algebra $T^+$. We endow $G^+$ with a group structure by defining the convolution product of ${\sf g}_1$ and ${\sf g}_2$ as
$$
({\sf g}_1*{\sf g}_2)(\tau) = ({\sf g}_1\otimes {\sf g}_2)\Delta^+\tau,
$$
for all $\tau\in T$. We write ${\sf g}^{-1}$ for the inverse of a character ${\sf g}\in G^+$ in this group structure. For any map $x\in{\bbR^{d_0}} \mapsto {\sf g}_x\in G^+$ we define for any $x,y\in{\bbR^{d_0}}$ the character
$$
{\sf g}_{yx} \defeq {\sf g}_y  * {\sf g}_x^{-1} .
$$
Similarly we define for any map $ { \sf \Pi } : T \rightarrow {\mathcal D}'({\bbR^{d_0}})$, any point $x\in{\bbR^{d_0}}$, a new map ${\sf \Pi}_x : T \rightarrow \mcD'({\bbR^{d_0}})$ by setting
$$
{\sf\Pi}_x = \big( {\sf \Pi}\otimes {\sf g_x}^{-1}\big)\Delta.
$$
For any function $\phi$, point $x\in{\bbR^{d_0}} $ and $\epsilon>0$ we set
$$
\phi_x^\lambda(\cdot) \defeq \epsilon^{-d} \, \phi \Big( \frac{\cdot-x}{\epsilon}\Big).
$$
Finally for any integer $n_0$ also define $\mathcal{Z}_{n_0}$ as the set of $C^{n_0}$ functions $\phi$ supported in the unit ball of ${\bbR^{d_0}} $ and such that $\norme{\phi}_{C^{n_0}}\leq 1$.

\bigbreak

\begin{defn} \label{def_model}
    Pick $n\geq \vert\beta_0\vert$. A \textbf{model} $\sf M=(\Pi,g)$ over a regularity structure $\mcT$ is a pair of maps
    $$
    {\sf \Pi}  : T \to C^{\beta_0}({\bbR^{d_0}}),   \qquad   {\sf g} : {\bbR^{d_0}} \to G^+
    $$
    with the following properties.    
    \begin{itemize}
    \item[--] For any $x\in {\bbR^{d_0}}$ and $\tau\in T_{\vert\tau\vert}$ we have
    $$
    \big| {\sf \Pi}_x(\tau)(\varphi_x^\epsilon) \big| \lesssim \epsilon^{|\tau|}
    $$
    uniformly in $x$ in  compact subsets of $\bbR^{d_0}$, in $\epsilon\in (0,1)$ and in $\varphi\in \mathcal{Z}_{n_0}$.

        \item[--] For any $x,y \in {\bbR^{d_0}}$ and $\mu\in T^+_{\vert\mu\vert}$ we have
    $$
    |{\sf g}_{yx}(\mu)| \lesssim |y-x|^{|\mu|}
    $$
    uniformly for $x,y$ in compact subsets of $\bbR^{d_0}$.
     \end{itemize}
\end{defn}

\ssk

\begin{defn*}
    Let $\mathscr{T}$ be a regularity structure and ${\sf M}=({\sf \Pi},{\sf g})$ be a model on it. For any $r\in\bfR$, a modelled distributions $\bm{f}\in\mcD^r(T, {\sf g})$ is a map $\bm{f} : {\bbR^{d_0}} \rightarrow \bigoplus_{r'<r}T_{r'}$ such that
    \begin{align*}
    &\max_{r'<r}\sup_{x\in{\bbR^{d_0}} } \norme{\bm{f}(x)}_{r'} <+\infty,
    \\
    &\max_{r'<r}\sup_{x,y\in{\bbR^{d_0}}} \frac{\big\Vert \bm{f}(y) - \widehat{\sf g}_{yx}(\bm{f}(x)) \big\Vert_{r'}}{|y-x|^{r-r'}} < +\infty.
    \end{align*}
\end{defn*}





\medskip

\subsection{Basics on analysis and proofs of three lemmas.}
\label{SectionBasicsAnalysis}

\hspace{0.1cm} For any function $f$ and any multi-index $\ell\in\bbN^{d_0}$ we define the modified Littlewood-Paley projector
$$
(\Delta_i^\ell f)(x) \defeq \int_{{\bbR^{d_0}} }K_i(x-y)(y-x)^\ell f(y) \, dy. 
$$

\ssk 

\begin{lem}
For $f\in C^r$ with $r>0$ one has $\big(\Delta_i^\ell f\big)_{i\geq -1} \in  {\sf C}^{r+|\ell|}$ and 
$$
\big\Vert (\Delta_i^\ell f\big)_{i\geq -1}\big\Vert_{{\sf C}^{r+|\ell|}} \lesssim \norme{f}_r.
$$
\end{lem}

\ssk

\begin{Dem}
    If $|i-j|\geq 2$ we have $\Delta_i^\ell (\Delta_j f) =\Delta_j(\Delta_i^\ell f) = 0$, so $\Delta_i^kf$ is spectrally supported in a ball $2^iB$ and 
    \begin{align*}
    (\Delta^\ell_if)(x) = \sum_{|j-i|\leq 1} \Delta_i^\ell(\Delta_jf)(x) = \sum_{|j-i|\leq 1} \int K_i(x-y)(x-y)^\ell (\Delta_jf)(y) \, dy
    \end{align*}
    Then we get
    \begin{align*}
    \big|(\Delta^\ell_if)(x)\big| &\leq \Big|\int K_i(z)z^\ell \dd z \Big| \sum_{|j-i|\leq 1}\norme{\Delta_jf}_{L^\infty} 
    \leq 2^{-ir}\Big|\int K_i(z)z^\ell \dd z \Big| \norme{f}_r
    \leq 2^{-i(r+|\ell|)}\norme{f}_r,
    \end{align*}
    using the scaling property of the kernel $K_i$ for the last inequality.
\end{Dem}

\ssk

Note that the sequence $(\Delta_i^\ell f)_{i\geq -1}$ does not represent the Littlewood-Paley blocks of any distribution as $\sum_i \Delta_i^\ell f = 0$ for any $\ell \neq 0$.

\medskip

\noindent \textbf{\textit{Proof of Lemma \ref{LemResteTaylor}}.} Pick ${\sf f}  = (f_i)_{i\geq -1} \in {\sf C}^r$ and $o>0$ with integer part $\lfloor o\rfloor$. If $f_i$ is spectrally supported in a ball $2^iB$, then $f_i(\cdot+h) - \sum_{|k|<o} \partial^kf_i\frac{h^k}{k!}$ is spectrally supported in the same ball $2^iB$. From Taylor Young inequality applied to $f_i$ at order $\lfloor o\rfloor+1$ and Bernstein inequality we have for any $x\in{\bbR^{d_0}}$
\begin{align*}
\bigg| f_i(x+h) - \sum_{|k|<o} \partial^kf_i(x)\frac{h^k}{k!} \bigg| &\lesssim |h|^{\lfloor o\rfloor+1} \norme{D^{\lfloor o\rfloor+1} f_i}_{L^\infty}
\\
&\lesssim |h|^{\lfloor o\rfloor+1} 2^{i(\lfloor o\rfloor+1)}\norme{f_i}_{L^\infty}.
\end{align*}
Similarly Taylor-Young inequality at order $\lfloor o\rfloor$ gives
\begin{align*}
\bigg| f_i(x+h) - \sum_{|k| < \lfloor o\rfloor} \partial^k f_i(x) \frac{h^k}{k!}\bigg| \lesssim |h|^{\lfloor o\rfloor} 2^{i\lfloor o\rfloor}\norme{f_i}_{L^\infty},
\end{align*}
from which we see that
\begin{align*}
\bigg| f_i(x+h) - \sum_{|k| < o} \partial^kf_i(x)\frac{h^k}{k!} \bigg| &\leq \bigg| f_i(x+h) - \sum_{|k| < \lfloor o\rfloor} \partial^k f_i(x) \frac{h^k}{k!} \bigg| + \bigg| \sum_{|k|=\lfloor o\rfloor} \partial^kf_i(x)\frac{h^k}{k!} \bigg|   \\
&\lesssim |h|^{\lfloor o\rfloor} 2^{i\lfloor o\rfloor} \norme{f_i}_{L^\infty} + |h|^{\lfloor o\rfloor} \norme{D^{\lfloor o\rfloor} f_i}_{L^\infty} \lesssim |h|^{\lfloor o\rfloor} 2^{i\lfloor o\rfloor} \norme{f_i}_{L^\infty}.
\end{align*}
We conclude by interpolation that we have
\begin{align*}
\bigg| f_i(x+h) - \sum_{|k| < o} \partial^kf_i(x)\frac{h^k}{k!} \bigg| \lesssim |h|^o 2^{io} \norme{f_i}_{L^\infty} \lesssim |h|^o \, 2^{-i(r-o)} \, \norme{{\sf f}}_r.
\end{align*}

\ssk

\noindent \textbf{\textit{Proof of Lemma \ref{lem_seuilopti}}.}  Let $\delta>0$ such that the estimate holds for $\theta\in [\gamma-\delta,\gamma+\delta]$. For $x,y\in\bfR^d$ with $|y-x|\leq 1$, one has for any integer $N$
\begin{align*}
\sum_{i\leq N} X_{yx}^i &\leq C\sum_{i\leq N} |y-x|^{\gamma+\eta}2^{i\eta} \lesssim C2^{N\eta} |y-x|^{\gamma+\eta},
\\
\sum_{i > N} X_{yx}^i &\leq C\sum_{i > N} |y-x|^{\gamma-\eta}2^{-i\eta} \lesssim C 2^{-N\eta} |y-x|^{\gamma-\eta},
\end{align*}
Choosing $N$ such that $|y-x|\simeq 2^{-N} $ gives the required bound.

\ssk

\noindent \textbf{\textit{Proof of Lemma \ref{LemProofPointB}.}} Using the definition of $\triangle_{h,r'}{\sf P}_<$ we have
\begin{align*}
\sum_{\sigma_e\prec\cdots\prec \tau_n(\bbfl )} \big(\triangle_{h,r'} &\widetilde{\sf P}_{<}^{|\tau_n(\bbfl )/\boldsymbol{\sigma}^+|_{\boldsymbol{\alpha'}}}\big) \Big([\tau_n(\bbfl )/\sigma_1]^{{{\sf M}^r}'},\dots,[\sigma_e]^{{{\sf M}^r}'},[n]^{{{\sf M}^r}'} \Big)(x)   \\
  &=\sum_{e\geq 0}\hspace{2pt}\sum_{\sigma_e\prec\cdots\prec \tau_n(\bbfl )} \widetilde{\sf P}_{<}^{|\tau_n(\bbfl )/\boldsymbol{\sigma}^+|_{\boldsymbol{\alpha'}}}\Big([\tau_n(\bbfl )/\sigma_1]^{{{\sf M}^r}'},\dots,[\sigma_e]^{{{\sf M}^r}'},[n]^{{{\sf M}^r}'} \Big)(x+h)   \\
  &\quad- \underset{\sigma_e\prec\cdots\prec \tau_n(\bbfl )}{\sum_{e\geq 0}} \sum_{|k|<r'} \partial_\star^k{\sf P}_{<}\Big([\tau_n(\bbfl )/\sigma_1]^{{{\sf M}^r}'},\dots,[\sigma_e]^{{{\sf M}^r}'},[n]^{{{\sf M}^r}'} \Big)(x) \, h^p
\end{align*}
\begin{align*}
  &\quad- \underset{\sigma_e\prec\cdots\prec \tau_n(\bbfl )}{\sum_{e\geq 0}} \; \underset{|k|<|\tau_n(\bbfl )/\sigma_m|}{\sum_{1\leq m\leq e}} \partial_\star^k{\sf P}_<\Big([\tau_n(\bbfl )/\sigma_1]^{{{\sf M}^r}'},\dots,[\sigma_{m-1}/\sigma_{m}]^{{{\sf M}^r}'} \Big) 
  \\
  &\hspace{4cm} \times \frac{h^k}{k!} \, \big(\triangle_{h,|\sigma_m|_{\alpha}+\alpha_n'}{\sf P}_<\big)\Big([\sigma_m/\sigma_{m-1}]^{{{\sf M}^r}'},\dots ,[n]^{{{\sf M}^r}'}\Big),
\end{align*}
where we use the shorthand
\begin{align*}
&\partial_\star^k{\sf P}_{<}\Big([\tau_n(\bbfl )/\sigma_1]^{{{\sf M}^r}'},\dots,[\sigma_{m-1}/\sigma_{m}]^{{{\sf M}^r}'} \Big) 
\\
&\hspace{3.2cm} = \sum_{\bbfk\in\mcP_m(k)}\binom{k}{\bbfk} \, \widetilde{\sf P}_{<}^{D^{\bbfk}|\tau_n(\bbfl )/\boldsymbol{\sigma}_{\leq m}|_{\alpha}}\Big( \partial^{k_1}[\tau_n(\bbfl )/\sigma_1]^{{{\sf M}^r}'},\dots, \partial^{k_m}[\sigma_{m-1}/\sigma_{m}]^{{{\sf M}^r}'} \Big) .
\end{align*}
The first line of the right hand side gives ${\sf P}_{\bbfl }\big( f_1,\dots, \Delta_r f_n\big)(x+h)$. As 
$$
\int_{{\bbR^{d_0}}} K_{<i}(h)h^p \,dh=0
$$ 
for $p\ne 0$, the second line of the right hand side gives a zero contribution when integrated against $K_{<i}$ except for $p=0$, in which case it gives ${\sf P}_{\bbfl }\big(f_1,\dots,\Delta_rf_n\big)(x)$. Then
\begin{align*}
  &\int_{{\bbR^{d_0}}}\sum_{\sigma_e\prec\cdots\prec \tau_n(\bbfl )} 
  \big(\triangle_{h,r'} \widetilde{\sf P}_{<}^{|\tau_n(\bbfl )/\boldsymbol{\sigma}^+|_{\boldsymbol{\alpha'}}}\big) \Big([\tau_n(\bbfl )/\sigma_1]^{{{\sf M}^r}'},\dots,[\sigma_e]^{{{\sf M}^r}'},[n]^{{{\sf M}^r}'} \Big)(x) \,dh   \\
  &=  \int_{{\bbR^{d_0}}} K_{<i}(h) \, {\sf P}_{\bbfl }\big(f_1,\dots,\Delta_r f_n\big)(x+h) \,dh-  {\sf P}_{\bbfl }\big(f_1,\dots,\Delta_r f_n\big)(x)   \\
  &\quad-\sum_{\sigma\prec \tau_n({\bbfl })} \sum_{|k|<|\tau_n(l)/\sigma|} \underset{\sigma\prec \sigma_{e_1}\prec\cdots\prec \tau_n(\bbfl )}{\sum_{e_1\geq 0}} \partial_\star^k{\sf P}_<\Big([\tau_n(l)/\sigma_1]^{{{\sf M}^r}'}, \dots, [\sigma]^{{{\sf M}^r}'}\Big)(x)
\end{align*}
\begin{align*}
  &\hspace{3.5cm}  \times \int_{{\bbR^{d_0}}} K_{<i}(h) \, \frac{h^k}{k!}\sum_{e_2\geq 0}\sum_{\nu_{e_2}\prec\cdots\prec \sigma} \triangle_{h,|\sigma|+\alpha_n'}{\sf P}_<\Big([\sigma/\nu_1]^{{{\sf M}^r}'},\dots,[n]^{{{\sf M}^r}'}\Big) \,dh.
\end{align*}
The $\sigma\in \mcB$ such that $\sigma\prec\tau_n(\bbfl )$ have form $\sigma = \llbracket m+1,n-1\rrbracket_{\bbfl +p_1}X^{p_2+s_2}$ and $\tau/\sigma=\llbracket 1,m\rrbracket_{\bbfl }^pX^{s_1}$ where $p=s_1+p_2$ and $l_{n-1}=s_1+s_2$. For such $\sigma$, using Lemma \ref{LemRelStarDerive} the sum 

$$
\sum_{\sigma\prec\cdots\prec \tau_n(\bbfl )}\partial_\star^k{\sf P}_<\Big([\tau_n(\bbfl )/\sigma_1]^{{{\sf M}^r}'},\dots,[\sigma_{m-1}/\sigma]^{{{\sf M}^r}'} \Big)(x)
$$ 
is $0$ if $s_1\ne 0 $, otherwise ($s_1=0$) this sum is equal to
$$
\sum_{\bbfp \in\mcP_m(p)}\binom{p}{\bbfp} \, \widetilde{\sf P}_{\bbfl _{\leq m}}^{D^{\bbfk+\bbfp }\alpha}\big(\partial^{k_1+p_1}f_1,\dots,\partial^{k_m+p_m}f_m\big)(x)={\sf g}^{\bm{f}}\big(\llbracket 1,m\rrbracket_{\bbfl _{\leq m}}^{\bbfk+p}  \big).
$$
On the other hand for $\sigma=\llbracket m+1,n\rrbracket_{\bbfl +p_1}X^{p_2+l_{n-1}}$, we have from the induction hypothesis 
\begin{align*}
&\int_{{\bbR^{d_0}}} K_{<i}(h) \, \frac{h^k}{k!}\sum_{e_2\geq 0}\sum_{\nu_{e_2}\prec\cdots\prec \sigma} \big(\triangle_{h,|\sigma|+\alpha_n'}{\sf P}_<\big)\Big([\sigma/\nu_1]^{{{\sf M}^r}'},\dots,[n]^{{{\sf M}^r}'}\Big) \,dh 
\\
&\hspace{7.5cm}=  \int_{{\bbR^{d_0}}} K_{<i}(h) \, {{\sf \Pi}_x^r}' \big( \llbracket m+1,n\rrbracket_{\bbfl +p}X^k \big)(x+h).
\end{align*}
One finally gets
\begin{align*}
&\int_{{\bbR^{d_0}}}\sum_{\sigma_e\prec\cdots\prec \tau_n(\bbfl )} \Big(\triangle_{h,r'} \widetilde{\sf P}_{<}^{|\tau_n(\bbfl )/\boldsymbol{\sigma}^+|_{\boldsymbol{\alpha'}}}\Big) \Big([\tau_n(\bbfl )/\sigma_1]^{{{\sf M}^r}'},\dots,[\sigma_e]^{{{\sf M}^r}'},[n]^{{{\sf M}^r}'} \Big)(x) \,dh
\end{align*}
\begin{align*}
  &=  \int_{{\bbR^{d_0}}} K_{<i}(h) \, \Big\{  {\sf P}_{\bbfl }\big(f_1,\dots,\Delta_r f_n\big)(x+h) -  {\sf P}_{\bbfl }\big(f_1,\dots,\Delta_r f_n\big)(x) 
  \\
  &\hspace{3cm} - \sum_{m,k,p} {{\sf g}^r_x}' \big(\llbracket 1,m\rrbracket_{\bbfl }^{k+p}\big) {{\sf \Pi}_x^r}'\big( X^k\llbracket m+1,n\rrbracket_{\bbfl +p}  \big)(x+h) \,dh \Big\}
  \\
  &= \int_{{\bbR^{d_0}} } K_{<i}(h) \, \Big\{ {{\sf \Pi}^r}'\big(\llbracket 1,n\rrbracket_{\bbfl  }\big)(x+h) - \sum_{\sigma<\llbracket 1,n\rrbracket_{\bbfl }}  {{\sf g}_x^r}' \big(\llbracket 1,n\rrbracket_{\bbfl }/\sigma\big)  {{\sf \Pi}_x^r}'(\sigma)(x+h)    \Big\} \,dh
  \\
  &= \int_{{\bbR^{d_0}} } K_{<i}(h) \, {{\sf \Pi}_x^r}' \big(\llbracket 1,n\rrbracket_{\bbfl  }\big)(x+h) \, dh,
\end{align*}
so we have indeed \eqref{eq_triangletopi}.

\medskip

\subsection{Proof of some algebraic lemmas.}
\label{SectionAppendixAlgebraicLemmas}

\hspace{0.1cm} We prove in this section a number of algebraic results that were used in the main body of the text. We start Section \ref{SectionAlgebraicProofs} by proving the inductive relation (Lemma \ref{rectilde2}) on the $\widetilde{\sf P}_<^{\beta^1\hspace{-0.05cm},\beta^2}$ that lead us in Section \ref{subsection_devtwisted} to the local expansion property satisfied by the $\widetilde{\sf P}_<^{\alpha}$ stated in Proposition \ref{PropExpansionTildePAlphaMinus}. The operators $\widetilde{\sf P}_<^{\beta^1\hspace{-0.05cm},\beta^2}$ have an analogue $\widetilde{\sf P}^{\beta^1\hspace{-0.05cm},\beta^2}_{\bbfl}$ defined from the (true) iterated paraproduct operator. The remainder of Section \ref{SectionAlgebraicProofs} is dedicated to proving Lemma \ref{lem_relhighcorr}, which is the analogue of Lemma \ref{rectilde2} for the operators $\widetilde{\sf P}^{\beta^1\hspace{-0.05cm},\beta^2}_{\bbfl}$. Lemma \ref{lem_relhighcorr} plays a crucial role in our proof of Lemma \ref{LemRelStarDerive}. The later is the main ingredient of our proof of Theorem \ref{ThmMain1}. The proof of Lemma \ref{LemRelStarDerive} occupies all of Section \ref{SectionProofLemmaRelStarDeriv}.

\subsubsection{Algebraic properties of the $\widetilde{\sf P}_<^{\beta^1\hspace{-0.05cm},\beta^2}$.}
\label{SectionAlgebraicProofs}

\hspace{0.1cm} We start with the 

\ssk

\noindent {\it Proof of Lemma \ref{rectilde2}.} The proof is very similar to the proof of Lemma \ref{rectilde1}. From \textbf{\textit{Assumption (A)}} we have the following partition of ${\sf MultiCut}(\beta^1\hspace{-0.05cm},\beta^2)$
$$
{\sf MultiCut}\big(\beta^1\hspace{-0.05cm},\beta^2\big) = {\sf MultiCut}\big(\beta^1\big) \sqcup \bigsqcup_{m\in {\sf Cut}(\beta^2)\backslash {\sf Cut}(\beta^1)} {\sf MultiCut}(\beta^2)\big[\beta^1,m\big],
$$
where
$$
{\sf MultiCut}(\beta^2)\big[\beta^1,m\big] \defeq \Big\{ \mathbf{i} \in {\sf MultiCut} \big(\beta^1, \beta^2 \big) \,;\, m\in \mathbf{i}, \quad \sum_{s=1}^m\beta_s^2 = \min_{j\in\mathbf{i},\enskip j\notin {\sf Cut}(\beta^1) } \sum_{s=1}^j\beta_s^2 \Big\}.
$$
We thus have

\begin{equation*} \begin{split}
\widetilde{\sf P}_<^{\beta^1 ,\beta^2}({\sf h}_1,\dots {\sf h}_n) &= {\sf P}_<^{\beta^1} \big({\sf h}_1,\dots,{\sf h}_n\big )   \\
&\qquad+\sum_{m\in {\sf Cut}(\beta^2 )\backslash {\sf Cut}(\beta^1)}  \sum_{\mathbf{i} \in {\sf MultiCut}(\beta^2)[\beta^1,m]}\hspace{-0.3cm}(-1)^{n(\bm{d})+1} \prod_{k=1}^{n(\bm{d})}{\sf P}_<\big({\sf h}_{i_{k-1}+1},\dots,{\sf h}_{i_k}\big).
\end{split} \end{equation*}
so it suffices to show that for any $m\in {\sf Cut}(\beta^2) \backslash  {\sf Cut}(\beta^1)$ we have
\begin{equation*} \begin{split}
\sum_{\mathbf{i} \in {\sf MultiCut}(\beta^2)[\beta^1,m]} (-1)^{n(\bm{d})+1} &\prod_{k=1}^{n(\bm{d})}{\sf P}_<\big({\sf h}_{i_{k-1}+1},\dots,{\sf h}_{i_k}\big)   \\
&= -\widetilde{\sf P}^{\beta^1_{\leq m},\beta^2_{\leq m}}_< \big({\sf h}_1,\dots, {\sf h}_m \big)\widetilde{\sf P}^{\beta^1_{> m},\beta^2_{> m}}_<\big({\sf h}_{m+1},\dots, {\sf h}_n\big).
\end{split} \end{equation*}
Pick $m\in {\sf Cut}(\beta^2) \backslash  {\sf Cut}(\beta^1)$. We prove that: For $1<j<m$ we have 
$$
\Big\{\exists \, \mathbf{i} \in {\sf MultiCut}(\beta^2)[\beta^1,m] \,;\,  j\in \mathbf{i} \Big\} \Leftrightarrow \Big\{j \in {\sf Cut}\big(\beta^1_{\leq m})\cup {\sf Cut}( \beta^2_{\leq m}\big) \Big\},
$$
and for $m<j<n$ we have 
$$
\Big\{\exists \, \mathbf{i} \in {\sf MultiCut}(\beta^2)[\beta^1,m] \,;\,  j\in \mathbf{i} \Big\} \Leftrightarrow \Big\{j-m \in {\sf Cut}\big(\beta^1_{> m})\cup {\sf Cut}( \beta^2_{> m} \Big\}.
$$ 
The proof of Lemma \ref{rectilde2} follows from these equivalences as in the proof of Lemma \ref{rectilde1}. 

As a preliminary remark we note that for $m\in {\sf Cut}(\beta^2) \backslash  {\sf Cut}(\beta^1)$ we have $\sum_{s=1}^m\beta^1_s > 0 $ and  $\sum_{s=m+1}^n\beta^1_s > 0 $. We prove now the first equivalence relation. Suppose $\boldsymbol{i}\in {\sf MultiCut}(\beta^2)[\beta^1,m] $ and $1<j<m$ is such that $j\in\boldsymbol{i}$. If $j\in {\sf Cut}(\beta^1)$ then $j\in {\sf Cut}(\beta^1_{\leq m})$. Otherwise $j\in {\sf Cut}(\beta^2) \backslash {\sf Cut}(\beta^1) $ and $j\in {\sf Cut}(\beta^2_{\leq m})$. Reciprocally if $j\in {\sf Cut}(\beta^1_{\leq m}) \cup {\sf Cut}(\beta^2_{\leq m})$ then necessarily $\sum_{s=j+1}^m\beta^2_s<0$ and $j\in {\sf Cut}(\beta^2)$ and $\sum_{s=j+1}^n\beta^2_s<\sum_{s=m+1}^n\beta^2_s$. 

The second equivalence relation is proved in the same way.   \hfill $\rhd$

\ssk

The remainder of this section is dedicated to stating and proving an analogue of Lemma \ref{rectilde2} for some operator $\widetilde{\sf P}_\ell^{\beta^1,\beta_2}$ that we can associate to the iterated paraproduct operators $\widetilde{\sf P}_\ell$. We first need an ad hoc setting to introduce these operators. It is very close to the setting of Section \ref{RSIPP}.

Fix $n\geq 1$. Define the set of symbols 
\begin{equation*} \begin{split}
\widehat\mcB \defeq \;&\Big\{ \llbracket a,b\rrbracket_{\bbfl }^{\bbfk} X^m\,;\, 1\leq a\leq b\leq n, \; \ell,k\in\bbN^{d_0}, \bbfl \in \mcP_{b-a}(\ell),\bbfk\in \mcP_{b-a+1}(k), m\in\bbN^{d_0}\Big\}   \\
&\cup\{X^m\}_{m\in\bbN^{d_0}}.
\end{split} \end{equation*}
Given $\beta\in\bbR^n$ and $\tau = \llbracket a,b\rrbracket_{\bbfl }^{\bbfk}\in\widehat\mcB$ we set
$$
|\tau|_{\beta} =  \sum_{j=a}^b \beta_j - |k| + |\ell|.
$$ 
We denote by $\widehat T$ the vector space freely spanned by the elements of $\widehat\mcB$, and for $\tau=\llbracket a,b\rrbracket_{\bbfl }^{\bbfk}$ we set
$$
\oplus(\tau) \defeq \Big\{ \llbracket a,c\rrbracket_{\bbfl _{< c-a}}^{\bbfk+\bbfp} \,;\, a<c<b, \, \ell_{c-a}= 0 
\Big\} ,
$$
for $\sigma = \llbracket a,c\rrbracket_{\bbfl _{< c-a}}^{\bbfk+\bbfp} \in\oplus(\tau)$ define an element of $\widehat{T}$ setting
$$
(\tau\backslash\sigma) \defeq  \sum_{k=p_1+p_2} \binom{k}{p_1} \llbracket c+1,b\rrbracket_{\bbfl _{>j_1-j}+p_1} \, X^{p_2}     
$$
Finally we define a coproduct $\widehat{\Delta}: \widehat T \rightarrow \widehat T\otimes \widehat T$ setting
$$
\widehat\Delta(\tau) = \sum_{\sigma\in \oplus(\tau)}(\tau\backslash\sigma) \otimes \sigma.
$$
Proceeding as in the proof of Proposition \ref{PropRS} one can see that $\Delta$ is co-associative. We note in particular that all the elements of $T$ in the sum defining $(\tau\backslash\sigma)$ have the same homogeneity. Re-indexing the sum defining $\widehat\Delta(\tau)$ we can write
$$
\widehat{\Delta}(\tau) \eqdef \sum_{\nu\widehat{<}\tau}\nu\otimes \tau/\nu,
$$
with $\nu$ running in the basis $\mcB$ of $T$ and $\tau/\nu$ defined by this identity. The element $\tau/\nu$ of $\widehat{T}$ is a sum of terms with the same $\vert\cdot\vert_\beta$-homogeneity, so we can abuse notations and write $\vert\tau/\nu\vert_\beta$.

\ssk

For $\tau\in\widehat\mcB$ we define the set of cuts
$$
\widehat{\sf Cut}(\tau,\beta) \defeq \Big\{ \sigma \widehat< \tau \,;\, |\sigma|_{\beta}<0 \enskip\text{and} \enskip |\tau/\sigma|_{\beta}>0 \Big\},
$$
and the set of multicuts
$$
\widehat{\sf MultiCut}(\tau,\beta) \defeq \Big\{\boldsymbol{\sigma}= (\sigma_1,\cdots,\sigma_{e(\bm{\sigma})}) \in \widehat{\sf Cut}(\tau)^{e(\boldsymbol{\sigma})} \,;\, {e(\bm{\sigma})}\geq1,\enskip \sigma_{e(\boldsymbol{\sigma})}\widehat<\cdots\widehat<\sigma_1\widehat<\tau           \Big\}.
$$
For a fixed tuple $\bm{g}=(g_1,\dots,g_n)$ of distributions, for $\sigma=\llbracket a,b\rrbracket_{\bbfl }^{\bbfk}\in\widehat\mcB$ we set
$$
\Upsilon_{\bm{g}}(\sigma)\defeq{\sf P}_{\bbfl }\big( \partial^{k_a}g_a,\dots,\partial^{k_b} g_b\big).
$$
We note that for any $\bbfp\in(\bbN^{d_0})^n$, setting $\partial^{\bbfp}\bm{g} = (\partial^{p_1}g_1,\dots,\partial^{p_n}g_n)$, one has
\begin{equation}\label{eq_commutupsilder}
\Upsilon_{\partial^\bbfp\mathbf{g}}\big(\llbracket a,b\rrbracket_{\bbfl }^{\bbfk}\big)=\Upsilon_{\bm{g}}\big(\llbracket a,b\rrbracket_{\bbfl }^{\bbfk + \bbfp }\big).
\end{equation}

\ssk

\begin{lem}\label{lemdef2highcorr}
For any $\ell\in\bbN^{d_0}$ and $\bbfl \in\mcP_{n-1}(\ell)$, letting $\tau= \llbracket 1,n\rrbracket_{\bbfl }^0$, we have
\begin{equation} \label{EqInductionTildePEllBeta}
\widetilde{\sf P}^{\beta}_{\bbfl }\big(g_1,\cdots,g_n\big) = {\sf P}_{\bbfl }\big(g_1,\cdots,g_n\big) - \hspace{-0.15cm}\sum_{ \boldsymbol{\sigma} \in \widehat{\sf MultiCut}(\tau,\beta)} (-1)^{e(\boldsymbol{\sigma})+1} \, \Upsilon_{\bm{g}}(\tau/\sigma_1)\Upsilon_{\bm{g}}(\sigma_1/\sigma_2) \cdots \Upsilon_{\bm{g}}(\sigma_{e(\boldsymbol{\sigma})}).
\end{equation}
\end{lem}

\ssk

\begin{Dem}
We prove \eqref{EqInductionTildePEllBeta} by induction on $n$. The result is true for $n=1$. We prove that the right hand side of \eqref{EqInductionTildePEllBeta} satisfies the same recursive relation as $\widetilde{\sf P}^{\beta}_{\bbfl }\big(g_1,\dots,g_n\big)$. The proof is analogous to the proof of Lemma \ref{rectilde1}. 

From \textit{\textbf{Assumption (A)}} we have a partition 
    $$
    \widehat{\sf MultiCut}(\tau,\beta) = \bigsqcup_{\sigma\widehat<\tau} \widehat{\sf MultiCut}(\tau,\beta)[\nu],
    $$
    where
    $$
    \widehat{\sf MultiCut}(\tau,\beta)[\nu] \defeq \bigg\{ \bm{\sigma} = (\sigma_1,\cdots,\sigma_{e(\bm{\sigma})})\in \widehat{\sf MultiCut}(\tau,\beta) \,;\,
    \nu\in\bm{\sigma}, \, |\tau/\nu|=\min_{1\leq j \leq e( \bm{\sigma})}|\tau/\sigma_j| \bigg\}.
    $$
For any $\nu\in \widehat{\sf Cut}(\tau,\beta)$ and $\mu\widehat<\nu$ we have the equivalence
    $$
    \Big\{\exists \, \boldsymbol\sigma\in \widehat{\sf MultiCut}(\tau,\beta)[\nu] \,;\, \mu\in\boldsymbol\sigma \Big\} \Leftrightarrow \Big\{ \mu \in \widehat{\sf Cut}(\nu,\beta) \Big\}.
    $$
Likewise, for $\nu\widehat<\mu\widehat<\tau$ we have
    $$
    \Big\{\exists \, \boldsymbol\sigma\in \widehat{\sf MultiCut}(\tau,\beta)[\nu] \,;\, \mu\in\boldsymbol\sigma \Big\} \Leftrightarrow \Big\{ \mu/\nu \in \widehat{\sf Cut}(\tau/\nu,\beta) \Big\}.
    $$
    Define
    $$
    \overline{\Upsilon}_{\bm{g}}(\tau,\beta) \defeq \sum_{ \boldsymbol{\sigma} \in \widehat{\sf MultiCut}(\tau,\beta)} (-1)^{{e(\boldsymbol{\sigma})}+1} \Upsilon_{\bm{g}}(\tau/\sigma_1)\Upsilon_{\bm{g}}(\sigma_1/\sigma_2)\cdots \Upsilon_{\bm{g}}(\sigma_{e(\boldsymbol{\sigma})}).
    $$
    Using the two equivalence relations above, the same computation as in the proof of Lemma \ref{rectilde1} gives that
    \begin{align*}
        \overline{\Upsilon}_{\bm{g}}(\tau,\beta) = -\sum_{\sigma\in \widehat{\sf Cut}(\tau,\beta)} \big({\Upsilon}_{\bm{g}}(\sigma) - \overline{\Upsilon}_{\bm{g}}(\sigma,\beta) \big) \, \Big({\Upsilon}_{\bm{g}}(\tau/\sigma) -\overline{\Upsilon}_{\bm{g}}(\tau/\sigma,\beta)\Big).
    \end{align*}
    From the induction hypothesis, for $\sigma = \llbracket c+1,n\rrbracket_{\bbfl +p} \in \widehat{\sf Cut}(\tau,\beta)$ we have 
    $$
    {\Upsilon}_{\bm{g}}(\sigma) -\overline{\Upsilon}_{\bm{g}}(\sigma,\beta)  = \widetilde{\sf P}_{\bbfl +p}^{\beta_{>m}}(g_{m+1},\dots,g_n).
    $$
    Likewise, for $\tau/\sigma = \llbracket 1,c\rrbracket^p_{\bbfl }$, using \eqref{eq_commutupsilder} we have
    \begin{align*}
    {\Upsilon}_{\bm{g}}(\tau/\sigma) -\overline{\Upsilon}_{\bm{g}}(\tau/\sigma,\beta)  
    &=\sum_{\bbfp \in\mcP_m(p)}\binom{p}{\bbfp }\Big\{ {\Upsilon}_{\bm{g}} \big(\llbracket 1,c\rrbracket_{\bbfl }^{\bbfp }\big) -\overline{\Upsilon}_{\bm{g}}\Big(\llbracket 1,c\rrbracket_{\bbfl }^{\bbfp },\beta\Big) \Big\}
    \\
     &=\sum_{\bbfp \in\mcP_m(p)}\binom{p}{\bbfp }\Big\{{\Upsilon}_{\partial^\bbfp\bm{g}}\big(\llbracket 1,c\rrbracket_{\bbfl }^{0}\big) - \overline{\Upsilon}_{\partial^\bbfp\bm{g}}\Big(\llbracket 1,c\rrbracket_{\bbfl }^0 , \beta - \vert\bbfp\vert \Big) \Big\}
    \\
    &= \sum_{\bbfp \in\mcP_m(p)}\binom{p}{\bbfp} \, \widetilde{\sf P}_{\bbfl }^{\beta_{\leq c}-\vert\bbfp\vert}\big(\partial^{p_1}g_{1},\dots,\partial^{p_m}g_m\big).
    \end{align*}
This closes the induction step.
\end{Dem}

\ssk

Define 
\begin{equation*} \begin{split}
&\widehat{\sf MultiCut}(\tau,\beta^1\hspace{-0.05cm},\beta^2)   \\
&\defeq \Big\{ \boldsymbol{\sigma}=(\sigma_1,\cdots,\sigma_{e(\bm{\sigma})})\in \big(\widehat{\sf Cut}(\tau,\beta^1) \cup \widehat{\sf Cut}(\tau,\beta^2)\big)^{e(\bm{\sigma})} \,;\, e(\bm{\sigma})\geq1,\enskip \sigma_{e(\bm{\sigma})}\widehat<\cdots\widehat<\sigma_1\widehat<\tau           \Big\},
\end{split} \end{equation*}
and set 
$$
\widetilde{\sf P}^{\beta^1\hspace{-0.05cm},\beta^2}_{\bbfl }\hspace{-0.07cm}\big(g_1,\dots,g_n\big) \hspace{-0.07cm} \defeq {\sf P}_{\bbfl }\big(g_1,\dots,g_n\big) \hspace{0.15cm} - \hspace{-0.5cm}\sum_{ \boldsymbol{\sigma} \in {\sf MultiCut}(\tau,\beta^1\hspace{-0.05cm},\beta^2)} \hspace{-0.15cm} (-1)^{e(\bm{\sigma})+1} \Upsilon_{\bm{g}}(\tau/\sigma_1) \, \Upsilon_{\bm{g}}(\sigma_1/\sigma_2) \cdots \Upsilon_{\bm{g}}(\sigma_{e(\bm{\sigma})}).
$$

\ssk

\begin{lem} \label{lem_relhighcorr}
Suppose $\beta^1\hspace{-0.05cm},\beta^2$ are two tuples of real numbers such that $\beta^1_s\geq\beta^2_s$ for every $1\leq s\leq n$. Then we have
$$
\widetilde{\sf P}_{\bbfl }^{\beta^1\hspace{-0.05cm},\beta^2}\big(g_1,\dots,g_n\big) = \widetilde{\sf P}_{\bbfl }^{\beta^1}\big(g_1,\dots,g_n\big) - \sum \widetilde{\sf P}_{\bbfl }^{\beta^1-\vert\bbfk\vert, \beta^2-\vert\bbfk\vert}\big(\partial^{k_1}g_1,\dots,\partial^{k_c}g_c\big) \, \widetilde{\sf P}_{\bbfl  +k}^{\beta^1\hspace{-0.05cm},\beta^2}\big(g_{c+1},\dots,g_n\big),
$$
for a sum over $\llbracket c+1,n\rrbracket^{\bbfk}_{\bbfl }\in \widehat{\sf Cut}(\tau,\beta^2)\backslash \widehat{\sf Cut}(\tau,\beta^1)$.    
\end{lem}

\ssk

\begin{Dem}
The proof is the same as the proof of Lemma \ref{rectilde2}. Using \textit{\textbf{Assumption (A)}} we can partition of $\widehat{\sf MultiCut}(\tau,\beta^1\hspace{-0.05cm},\beta^2)$ as
$$
\widehat{\sf MultiCut}\big(\tau,\beta^1\hspace{-0.05cm},\beta^2\big) = \widehat{\sf MultiCut}(\tau,\beta^1) \sqcup \bigsqcup_{\nu\in \widehat{\sf Cut}(\tau,\beta^2)\backslash \widehat{\sf Cut}(\tau,\beta^1)} \widehat{\sf MultiCut}(\tau,\beta^2)\big[\beta^1,\nu\big],
$$
where
$$
\widehat {\sf MultiCut}(\tau,\beta^2)\big[\beta^1,\nu\big] \defeq \bigg\{ \bm{\sigma} \in \widehat{\sf MultiCut} \big(\tau,\beta^1, \beta^2 \big) \,;\, \nu\in \bm{\sigma}, \, |\tau/\nu|_{\beta_2} = \min_{\sigma\in\bm{\sigma}, \, \sigma\notin \widehat{\sf Cut}(\tau,\beta^1) } |\tau/\sigma|_{\beta_2} \bigg\}.
$$
Then we have
\begin{align*}
\widetilde{\sf P}_{\bbfl }^{\beta^1 ,\beta^2}(g_1,\cdots,g_n  ) = &\widetilde{\sf P}_{\bbfl }^{\beta^1} \big(g_1,\cdots,g_n\big )   \\
&\qquad+ \sum_{\substack{\nu\in \widehat{\sf Cut}(\tau,\beta^2 )\backslash \widehat{\sf Cut}(\tau,\beta^1)\\ \bm{\sigma} \in \widehat{\sf MultiCut}(\tau,\beta^2)[\beta^1,\nu] } } (-1)^{e(\bm{\sigma})+1} \,  \Upsilon_{\bm{g}}(\tau/\sigma_1)\Upsilon_{\bm{g}}(\sigma_1/\sigma_2)\cdots 
    \Upsilon_{\bm{g}}(\sigma_{e(\bm{\sigma})}).
\end{align*}
It suffices then to show for any $\nu=\llbracket 1,c\rrbracket_{\bbfl }^{\bbfk} \in \widehat{\sf Cut}(\tau,\beta^2) \backslash  \widehat{\sf Cut}(\tau,\beta^1)$ we have
\begin{align*}
\sum_{\bm{\sigma} \in \widehat{\sf MultiCut}(\tau,\beta^2)[\beta^1,\nu]} (-1&)^{e(\bm{\sigma})+1}  \Upsilon_{\bm{g}}(\tau/\sigma_1)\cdots \Upsilon_{\bm{g}}(\sigma_{e(\bm{\sigma})})   \\
&= -\widetilde{\sf P}_{\bbfl _{\leq m}}^{D^{\bbfk}\beta^1_{\leq m},D^{\bbfk}\beta^2_{\leq m}} \big(\partial^{k_1}f_1,\cdots,\partial^{k_m}f_m \big) \widetilde{\sf P}^{\beta^1_{> m},\beta^2_{> m}}_{\bbfl +k}\big(f_{m+1},\cdots,f_n \big).
\end{align*}
For such a $\nu$, we show below that for $\mu\widehat<\nu$ we have 
\begin{equation} \label{EqFirstEquivalence}
\Big\{ \exists \, \bm{\sigma} \in \widehat{\sf MultiCut}(\tau,\beta^2)[\beta^1,\nu],\enskip \mu\in \boldsymbol{\sigma} \Big\} \Leftrightarrow \Big\{ \mu \in \widehat{\sf Cut}(\nu,\beta^1)\cup\widehat{\sf Cut}(\nu, \beta^2) \Big\},
\end{equation}
and that for $ \nu\widehat<\mu\widehat<\tau $ we have 
\begin{equation} \label{EqSecondEquivalence}
\Big\{\exists \, \bm{\sigma} \in \widehat{\sf MultiCut}(\tau,\beta^2)[\beta^1,\nu],\enskip \mu\in \boldsymbol{\sigma} \Big\} \Leftrightarrow \Big\{ \tau/\mu \in \widehat{\sf Cut}(\tau/\nu,\beta^1)\cup\widehat{\sf Cut}(\tau/\nu, \beta^2) \Big\}.
\end{equation}
We can then conclude the proof of our lemma in the same way as in the proof of Lemma \ref{rectilde1}. 

A basic observation we will use is that for $\nu=\llbracket 1,c\rrbracket_{\bbfl }^{\bbfk} \in \widehat{\sf Cut}(\tau,\beta^2) \backslash  \widehat{\sf Cut}(\tau,\beta^1)$ we necessarily have $|\nu|_{\beta^1} > 0 $ and  $|\tau/\nu|_{\beta^1} > 0 $. We now prove \eqref{EqFirstEquivalence}. Suppose $\bm{\sigma}\in \widehat{\sf MultiCut}(\tau,\beta^2)[\beta^1,\nu] $ and $\mu\widehat<\nu$ such that $\mu\in\bm{\sigma}$. If $\mu\in \widehat{\sf Cut}(\tau,\beta^1)$, then $\mu\in \widehat{\sf Cut}(\nu,\beta^1)$ and otherwise $\mu\in \widehat{\sf Cut}(\tau,\beta^2) \backslash \widehat{\sf Cut}(\tau,\beta^1) $, then $\mu\in \widehat{\sf Cut}(\mu,\beta^2)$. Reciprocally if $\mu\in \widehat{\sf Cut}(\nu,\beta^1, \beta^2)$, then necessarily $|\nu/\mu|_{\beta^2}<0$ and $\mu\in \widehat{\sf Cut}(\tau,\beta^2)$ and $|\tau/\mu|_{\beta^2}<|\tau/\nu|_{\beta^2}$. We proceed similarly to prove the equivalence \eqref{EqSecondEquivalence}.
\end{Dem}

\medskip

\subsubsection{Proof of Lemma \ref{LemRelStarDerive}.}
\label{SectionProofLemmaRelStarDeriv}

\textit{\textbf{We first prove point \textit{\textbf{(i)}} }} by induction. From the definition of the operator $\widetilde{\sf P}_<$ we have
\begin{align*}
\widetilde{\sf P}_<^{|\tau/\boldsymbol{\sigma}|_{\alpha-\vert\bbfk\vert }}\Big( [\tau/\sigma_1]^{\partial^\bbfk \bm{f}},\dots,[\sigma_e]^{\partial^\bbfk \bm{f}}  \Big) 
&=  {\sf P}_<\Big( [\tau/\sigma_1]^{\partial^\bbfk \bm{f}},\dots,[\sigma_e]^{\partial^\bbfk \bm{f}}  \Big)
\\
&\quad- \sum_c  \widetilde{\sf P}_<^{|\tau/\boldsymbol{\sigma}_{\leq c}|_{\alpha-\vert\bbfk\vert}} \Big( [\tau/\sigma_1]^{\partial^\bbfk \bm{f}},\cdots,[\sigma_{c-1}/\sigma_c]^{\partial^\bbfk \bm{f}}  \Big)
\\
&\hspace{1.5cm}\times\widetilde{\sf P}_<^{|\sigma_c/\bm{\sigma}_{>c}|_{\alpha-\vert\bbfk\vert}} \Big( [\sigma_c/\sigma_{c+1}]^{\partial^\bbfk \bm{f}},\cdots,[\sigma_e]^{\partial^\bbfk \bm{f}}  \Big),
\end{align*}
with a sum over the set of integers $c\in \llbracket 1,e\rrbracket$ such that $|\tau/\sigma_c|_{\alpha-\vert\bbfk\vert }>0$ and $|\sigma_c|_{\alpha-\vert\bbfk\vert }<0$. Summing over the set of descending sequences $\sigma_e\prec\cdots\prec\sigma_1\prec \tau = \llbracket 1,n\rrbracket_{\bbfl }X^m$, we obtain that 

$$
\sum_{e\geq 0} \hspace{2pt}\sum_{\sigma_e\prec\cdots\prec\sigma_1\prec \tau} \widetilde{\sf P}_<^{|\tau/\boldsymbol{\sigma}|_{\alpha-\vert\bbfk\vert }}\Big( [\tau/\sigma_1]^{\partial^\bbfk \bm{f}},\dots,[\sigma_e]^{\partial^\bbfk \bm{f}}  \Big)
$$
is equal to 

\begin{equation} \label{eqlem30b} \begin{split}
    &\sum_{e\geq 0} \hspace{2pt}\sum_{\sigma_e\prec\cdots\prec\sigma_1\prec \tau} {\sf P}_<\Big( [\tau/\sigma_1]^{\partial^\bbfk \bm{f}},\cdots,[\sigma_e]^{\partial^\bbfk \bm{f}}  \Big) - \\
    &\sum_{\substack{\sigma\prec \tau \\ \sigma\in \widehat{\sf Cut}(\tau,\alpha-\vert\bbfk\vert )}} \, \underset{\sigma\prec\sigma_{e_1}\prec\cdots\prec \tau}{\sum_{e_1\geq 0}} \widetilde{\sf P}_<^{|\tau/\boldsymbol{\sigma}|_{\alpha-\vert\bbfk\vert }}\Big( [\tau/\sigma_1]^{\partial^\bbfk \bm{f}},\dots,[\sigma_{e_1}/\sigma]^{\partial^\bbfk \bm{f}}  \Big) 
    \\
    &\hspace{5cm}\times\sum_{e_2\geq 0}\sum_{\nu_{e_2}\prec\cdots\prec\sigma} \widetilde{\sf P}_<^{|\sigma/\boldsymbol{\nu}|_{\alpha-\vert\bbfk\vert }}\Big( [\sigma/\nu_{1}]^{\partial^\bbfk \bm{f}},\dots,[\nu_{e_2}]^{\partial^\bbfk \bm{f}}  \Big) 
\end{split} \end{equation}
From Lemma \ref{PropoRepresentationProposition} the first sum in \eqref{eqlem30b} is equal to ${\sf P}_{\bbfl }\big(\partial^{k_1}f_1,\cdots,\partial^{k_n}f_n \big)$ if $m=0$ and $0$ otherwise. 

\ssk

For the second double sum in the right hand side of \eqref{eqlem30b}, note first that all the homogeneities in the tuple $|\tau/\boldsymbol{\sigma}|_{\alpha-\vert\bbfk\vert }$ are positive. It follows that we have
\begin{align*}
\widetilde{\sf P}_<^{|\tau/\boldsymbol{\sigma}|_{\alpha-\vert\bbfk\vert }}\Big( [\tau/\sigma_1]^{\partial^\bbfk \bm{f}},\dots,[\sigma_{c-1}/\sigma]^{\partial^\bbfk \bm{f}}  \Big) = {\sf P}_<\Big( [\tau/\sigma_1]^{\partial^\bbfk \bm{f}},\dots,[\sigma_{c-1}/\sigma]^{\partial^\bbfk \bm{f}}  \Big).
\end{align*} 
Now, the elements $\sigma\prec\tau$ have the form $\sigma=\llbracket c+1,n\rrbracket_{\bbfl _{>c}+p_1}X^{p_2+m_1}$, and $\tau/\sigma=\llbracket 1,c\rrbracket_{\bbfl _{<c}}^{p}X^{m_2}$ with $m=m_1+m_2$ and $p=p_1+p_2$, so it follows from Lemma \ref{PropoRepresentationProposition} that
\begin{align*}
&\sum_{e_1\geq 0}\sum_{\sigma\prec\sigma_{e_1}\prec\cdots\prec \tau} {\sf P}_<\Big( [\tau/\sigma_1]^{\partial^\bbfk \bm{f}},\dots,[\sigma_{m-1}/\sigma]^{\partial^\bbfk \bm{f}}  \Big) 
\\
&\hspace{4cm}= \mathbf{1}_{m_2=0}\sum_{\bbfp \in\mcP_c(p)} \binom{p}{\bbfp} \, \widetilde{\sf P}_{\bbfl _{<c}}^{\alpha_{\leq c} - \vert\bbfk+\bbfp\vert} \big( \partial^{k_1+p_1}f_1,\dots,\partial^{k_c+p_c}f_c\big)
\end{align*}
Also, the induction hypothesis gives
$$
\sum_{e_2\geq 0}\sum_{\nu_{e_2}\prec\cdots\prec\sigma} \widetilde{\sf P}_<^{|\sigma/\boldsymbol{\nu}|_{\alpha-\vert\bbfk\vert }}\Big( [\sigma/\nu_{1}]^{\partial^\bbfk \bm{f}},\dots,[\nu_{e_2}]^{\partial^\bbfk \bm{f}}  \Big)  = \mathbf{1}_{p_2+m_1=0} \, \widetilde{\sf P}_{\bbfl _{>c}+p_1}^{\alpha-\vert\bbfk\vert _{>c}}\big(\partial^{k_{c+1}}f_{c+1},\dots,\partial^{k_n}f_n\big).
$$
If $m\neq 0$ then $m_1\neq 0 $ or $m_2+p_1 \neq 0$, and then all the terms in the right hand side of \ref{eqlem30b} add up to $0$; this closes the induction in that case. If now $m=0$, the non-zero terms in the sum over $\sigma\prec\tau$ are the terms with $\sigma=\llbracket c+1,n\rrbracket_{\bbfl +p}$ and $\tau/\sigma=\llbracket 1,c\rrbracket^p_\bbfl$, and
\begin{align*}
    &\sum_{e\geq 0} \hspace{2pt}\sum_{\sigma_e\prec\cdots\prec\sigma_1\prec \tau} \widetilde{\sf P}_<^{|\tau/\boldsymbol{\sigma}|_{\alpha-\vert\bbfk\vert }}\big( [\tau/\sigma_1]^{\partial^\bbfk \bm{f}},\dots,[\sigma_e]^{\partial^\bbfk \bm{f}}  \big) 
    = {\sf P}_{\bbfl }\big(\partial^{k_1}f_1,\dots,\partial^{k_n}f_n \big)
    \\ 
    &\hspace{1.6cm}-\sum_{c,p} \sum_{\bbfp \in\mcP_c(p)}\widetilde{\sf P}_{\bbfl _{<c}}^{\alpha_{\leq c} - \vert\bbfk+\bbfp\vert}\big( \partial^{k_1+p_1}f_1,\dots,\partial^{k_c+p_c}f_c\big) \, \widetilde{\sf P}^{\alpha_{>c} - \vert\bbfk\vert}_{\bbfl _{>c}+p}\big(\partial^{k_{c+1}}f_{c+1},\dots,\partial^{k_n}f_n\big),
\end{align*}
where the sum in the right hand side runs over the paris $(c,p)$ such that $\ell_c=0$, $|\llbracket 1,c\rrbracket^{p}_{\bbfl _{<c}}|_{\alpha-\vert\bbfk\vert } > 0$ and $|\llbracket c+1,n\rrbracket_{\bbfl _{>c}+p}|_{\alpha-\vert\bbfk\vert} < 0$. It follows then from recursive definition of the correctors $\widetilde{P}_{\bbfl }^{\beta}$ that the above quantity is indeed equal to $\widetilde{\sf P}_{\bbfl }^{\alpha-\vert\bbfk\vert }\big(\partial^{k_1}f_1,\dots,\partial^{k_n}f_n\big)$.   \vspace{0.15cm}

\noindent -- \textbf{\textit{We now prove point $\mathbf{(ii)}$}} by proving the following stronger statement: For $\tau = \llbracket 1,n\rrbracket_{\bbfl}X^m\in T$, and for {\it any} $p\in\bbN^{d_0}$, one has
    \begin{align*}
    &\sum_{e\geq 0} \hspace{2pt}\sum_{\sigma\prec\sigma_e\prec\cdots\prec\tau} \partial^p_\star{\sf P}_<\Big( [\tau/\sigma_1]^{\bm{f}},\dots,[\sigma_e/\sigma]^{\bm{f}}  \Big) 
    \\
    &\hspace{2.5cm}= \mathbf{1}_{m=0} \sum_{\bbfp \in\mcP_n(p)}\sum_{\bbfk\in\mcP_n(k)} \binom{p}{\bbfp }\binom{k}{\bbfk} \, \widetilde{\sf P}_{\bbfl }^{\alpha-\vert\bbfk\vert , \alpha - \vert\bbfk+\bbfp\vert}\Big(\partial^{k_1+p_1}f_1,\dots,\partial^{k_n+p_n}f_n \Big),
    \end{align*}
where
$$
\partial^p_\star{\sf P}_<\Big( [\tau/\sigma_1]^{\bm{f}},\dots,[\sigma_e/\sigma]^{\bm{f}}  \Big) = \sum_{\bbfp \in\mcP_{e+1}(p)} \binom{p}{\bbfp} \, \widetilde{\sf P}^{|\tau/{\bm{\sigma}}|_{\alpha} - \vert\bbfp\vert}_<\Big( \partial^{p_1}[\tau/\sigma_1]^{\bm{f}},\dots,\partial^{p_{e+1}}[\sigma_e/\sigma]^{\bm{f}} \Big).
$$
The proof of this fact relies on Lemma \ref{lem_relhighcorr} and is an induction over $n$. The result is true for $n=1$; suppose it holds true for $(n-1)$. From the definition of $\partial^p_\star{\sf P}_<$ and the recursive relation of Lemma \ref{rectilde1}, for any descending sequence $\sigma\prec\sigma_e\prec \cdots\prec \tau$, the distribution $\partial^p_\star{\sf P}_<\big( [\tau/\sigma_1]^{\bm{f}},\cdots,[\sigma_e/\sigma]^{\bm{f}}  \big) $ is equal to
\begin{align*}
         \sum_{\bbfp  \in \mcP_e(p)} \binom{p}{\bbfp } \widetilde{\sf P}_<^{|\tau/\boldsymbol \sigma| - \vert\bbfp\vert} &\Big( \partial^{p_1}[\tau/\sigma_1]^{\bm{f}},\dots,\partial^{p_q}[\sigma_e/\sigma]^{\bm{f}}  \Big) 
	\\
	&=  \sum_{\bbfp  \in \mcP_e(p)} \binom{p}{\bbfp }\bigg\{ {\sf P}_<\Big( \partial^{p_1}[\tau/\sigma_1]^{\bm{f}},\dots,\partial^{p_e}[\sigma_e/\sigma]^{\bm{f}}  \Big) \quad  (\cdots)
\end{align*}
\begin{align*}
        &\qquad- \sum_{c\in {{\sf Cut}(|\tau/\bm{\sigma}| -\vert\bbfp\vert)}} \widetilde{\sf P}_<^{|\tau/\bm\sigma_{\leq c}| - \vert\bbfp _{\leq c}\vert}\Big( \partial^{p_1}[\tau/\sigma_1]^{\bm{f}},\dots,\partial^{p_c}[\sigma_{c-1}/\sigma_c]^{\bm{f}}  \Big)
        \\
        &\hspace{5cm} \times \widetilde{\sf P}_<^{|\sigma_c/\bm\sigma_{>c}| - \vert\bbfp _{>c}\vert} \Big( \partial^{p_{c+1}}[\sigma_c/\sigma_{c+1}]^{\bm{f}},\dots,\partial^{p_e}[\sigma_e/\sigma]^{\bm{f}}  \Big) \bigg\}.
\end{align*}
Then, summing over descending sequences $\sigma\prec\sigma_e\prec \cdots\prec \tau$ and inverting the sums over $\bbfp $ and $c$, we obtain
\begin{align*}
&\sum_{e\geq 0} \hspace{2pt}\sum_{\sigma\prec\sigma_e\prec\cdots\prec\tau} \partial^p_\star{\sf P}_<\Big( [\tau/\sigma_1]^{\bm{f}},\dots,[\sigma_e/\sigma]^{\bm{f}}  \Big) 
\\
&=\partial^p \Big\{ \sum_{e\geq 0} \sum_{\sigma\prec\sigma_e\prec\cdots\prec\tau} {\sf P}_<\Big( [\tau/\sigma_1]^{\bm{f}},\dots,[\sigma_e/\sigma]^{\bm{f}}  \Big) \Big\}  
\\
& \quad -\sum_{\nu \prec\tau}\sum_{ \substack{ p=a+b \\ (a,b)\in {\sf C}(\tau,\nu,\sigma) }  }\binom{p}{a}\bigg\{ \sum_{\substack{e_1\geq 0 \\ \nu \prec\nu_{e_1}\prec\cdots\prec\tau}} \sum_{\bbfa\in \mcP_{e_1}(a)} \binom{a}{\bbfa} \, \widetilde{\sf P}_<^{|\tau/\boldsymbol{\nu}| - \vert\bbfa\vert} \Big( \partial^{a_1}[\tau/\nu_1]^{\bm{f}},\dots,\partial^{a_{e_1}}[\nu_{e_1}/\nu]^{\bm{f}}  \Big) 
\\
&\hspace{3.9cm} \times \sum_{\substack{e_2\geq 0 \\ \sigma\prec\sigma_{e_2}\prec\cdots\prec\nu}} \sum_{\bbfb \in \mcP_{e_2}(b)}\binom{b}{\bbfb} \widetilde{\sf P}_<^{|\nu/\bm{\sigma}| - \vert\bbfb\vert}\Big( \partial^{b_1}[\nu/\sigma_1]^{\bm{f}},\dots,\partial^{b_{e_2}}[\sigma_e/\sigma]^{\bm{f}}  \Big) 
\bigg\},
\end{align*}
where
$$
{\sf C}(\tau,\sigma,\nu)\defeq \Big\{   (a,b)\in(\bbN^{d_0})^2, \quad  |\tau/\nu|>|a|\enskip \text{and}\enskip |\nu/\sigma|<|b|   \Big\}. 
$$
From Lemma \ref{PropoRepresentationProposition} the first line of the right hand side of the last equality is equal  $0$ if $m\neq 0$, as $\sum_{i\geq -1} \Delta_i^m(g) = 0$ for any function $g$; it is equal to $\partial^p {\sf g}^{\bm{f}}(\tau)$ if $m=0$. 

We are able to use induction hypothesis for the remaining terms. Suppose first that $m\neq 0$. For any $\nu\prec\tau$ the elements $\tau/\nu,\nu/\sigma\in T^+$ have the form $\tau/\nu=\llbracket 1,c\rrbracket_{\bbfl }^{k_1'}X^{m_1}$ and $\nu/\sigma=\llbracket c+1,n\rrbracket^{k_2'}_{\bbfl +v}X^{m_2}$ where $m_1,m_2$ cannot be both equal to $0$. The induction assumption ensures in that case that the second term on the right hand side is $0$, which closes the induction step.

Suppose now that $m=0$. In this case for $\nu\prec\tau$, the elements $\tau/\nu,\nu/\sigma\in T^+$ have form $\tau/\nu = \llbracket 1,c\rrbracket_{\bbfl }^{k_a+v}$ and $\nu/\sigma = \llbracket c+1 , n\rrbracket^{k_b}_{\bbfl +v_1}X^{v_2}$ with $k=k_a+k_b$ and $v=v_1+v_2$. For $v_2\neq 0$ the induction assumption ensures that the sum over $e_2$ is null. We are thus left with the $\nu$ for which $v_2=0$. This leads to the equality
\begin{align*}
    &\sum_{e\geq 0} \hspace{2pt}\sum_{\sigma\prec\sigma_e\prec\cdots\prec\tau} \partial^p_\star{\sf P}_<\Big( [\tau/\sigma_1]^{\bm{f}},\dots,[\sigma_e/\sigma]^{\bm{f}}  \Big) 
    \\
    &=   \sum_{\bbfk\in\mcP_n(k)}\binom{k}{\bbfk} \bigg\{ \partial^p\widetilde{\sf P}_{\bbfl }^{\alpha-\vert\bbfk\vert  }\Big(\partial^{k_1}f_1,\dots,\partial^{k_n}f_n\Big)
    \\
    &\quad- \sum_{\substack{p=q+q'\\ c,v}} \binom{p}{q} \sum_{\bbfq\in\mcP_c(q)} \sum_{\bbfv\in\mcP_c(v)} \binom{q}{\bbfq} \binom{v}{\bbfv} \, {\sf P}_{\bbfl _{\leq m}}^{\alpha-\vert\bbfk\vert _{\leq c} , \alpha_{\leq m} - \vert\bbfk + \bbfq + \bbfv\vert} \Big(\partial^{k_1+q_1+v_1}f_1,\dots,\partial^{k_c+q_c+v_c}f_c\Big)
    \\
    &\hspace{3cm}\times\sum_{\bbfq'\in \mcP_{n-c}(q')}\binom{q}{\bbfq'} \widetilde{\sf P}_{\bbfl _{>c}+v}^{\alpha-\vert\bbfk\vert _{>c} , \alpha_{>c} - \vert\bbfk+\bbfq\vert}\Big(\partial^{k_{c+1}+q_{c+1}}f_{c+1}, \dots,\partial^{k_n+q_n}f_n\Big)\bigg\},  
    \end{align*}
where the sum over $q,q$ in $\bbN^{d_0}$ subject to $q+q'=p$, and $c,v$ runs over the indices such that $\llbracket 1,c\rrbracket_{\bbfl }^{\bbfk_{\leq c}+v}, \llbracket c+1,n\rrbracket_{\bbfl  +v}^{\bbfk_{>c}}\in T^+$ and $\ell_c=0$ and 
$$
|q| < { \big| \llbracket 1,c\rrbracket_{\bbfl }^{\bbfk_{\leq c}+v} \big|_{\alpha-\bbfq}}, \; \textrm{ and } \;  |q'| > { \big| \llbracket c+1,n\rrbracket_{\bbfl  +v}^{\bbfk_{>c}} \big|_{\alpha-\bbfq}}.
$$
This gives then

    \begin{align*}
    &\sum_{\bbfp \in \mcP_n(p)}\sum_{\bbfk\in\mcP_n(k)} \binom{p}{\bbfp } \binom{k}{\bbfk} \bigg\{ \widetilde{\sf P}_{\bbfl }^{\alpha-\vert\bbfk\vert }\Big( \partial^{k_1+p_1}f_1,\dots,\partial^{k_n+p_n}f_n \Big)
    \\
    &\hspace{2.5cm}- \sum_{\substack{c,v \\ \bbfv\in\mcP_c(v)}} \widetilde{\sf P}_{\bbfl _{\leq c}}^{\alpha-\vert\bbfk\vert _{\leq c} , \alpha_{\leq c} - \vert\bbfk+\bbfp\vert} \Big(\partial^{k_1+p_1+v_1}f_1,\dots,\partial^{k_c+p_c+v_c}f_c\Big)
    \\
    &\hspace{6cm}\times \widetilde{\sf P}_{\bbfl _{>c}+v}^{\alpha-\vert\bbfk\vert _{>c} , \alpha_{>c} - \vert\bbfk+\bbfb\vert} \Big(\partial^{k_{c+1}+p_{c+1}} f_{c+1}, \dots,\partial^{k_n+p_n} f_n\Big)   \bigg\},
\end{align*}
where the sum sum over $c,v$ runs over indexes such that $\llbracket c+1,n\rrbracket_{\bbfl +v}^{\bbfk>c} \in \widehat{\sf Cut}\big(\tau, \alpha - \vert\bbfk+\bbfp\vert\big) \backslash \widehat{\sf Cut}\big(\tau, \alpha-\vert\bbfk\vert\big)$. The result follows in that case from Lemma \ref{lem_relhighcorr}.

\bigskip
\bigskip
\bigskip

\noindent $\bullet$ \textbf{\textsf{I. Bailleul}} -- Univ Brest, CNRS UMR 6205, Laboratoire de Mathematiques de Bretagne Atlantique, F- 29200 Brest, France. {\it E-mail}: ismael.bailleul@univ-brest.fr. Partial support from the ANR-22-CE40-0017 grant.

\bigskip

\noindent $\bullet$ \textbf{\textsf{N. Moench}} -- Univ. Rennes, CNRS, IRMAR - UMR 6625, F-35000 Rennes, France. {\it E-mail: nicolas.moench@ens-rennes.fr}

\end{document}